\documentclass{amsart}
\usepackage{hyperref}
\usepackage{caption}
\usepackage{amsmath,amsfonts,amsthm}
\usepackage{mathtools}
\usepackage{float}
\usepackage{tikz}

\theoremstyle{definition}

\theoremstyle{plain}
\newtheorem{theorem}{Theorem}[section]

\newtheorem{lemma}[theorem]{Lemma}

\theoremstyle{remark}
\newtheorem*{remark}{Remark}

\newcommand{\homeo}{{\operatorname{Homeo^+}}}

\newcommand{\D}{\mathcal{D}}
\newcommand{\C}{\mathcal{C}}
\newcommand{\X}{\mathcal{X}}
\newcommand{\Y}{\mathcal{Y}}
\newcommand{\M}{\mathcal{M}}

\makeatletter
\def\subsection{\@startsection{subsection}{2}%
  \z@{.5\linespacing\@plus.7\linespacing}{.3\linespacing}%
  {\normalfont\bfseries}}
\makeatother

\title{Configuration space, moduli space and 3-fold covering space }
\author{Byung Chun Kim}
\address{Departments of Mathematics \\ Inha University \\ Incheon 402-751, Korea}
\email{wizardbc@gmail.com}

\author{Yongjin Song}
\address{Departments of Mathematics \\ Inha University \\ Incheon 402-751, Korea}
\email{yjsong@inha.ac.kr}

\thanks{The second author was supported by the Korean National Research Fund NRF-2016R1D1A1B03934531.}
\begin{document}

\maketitle

\begin{abstract}
A function from configuration space to moduli space of surface may induce a homomorphism between their fundamental groups which are braid groups and mapping class groups of surface, respectively. This map $\phi: B_k \rightarrow \Gamma_{g,b}$ is induced by 3-fold branched covering over a disk with some branch points. In this thesis we give a concrete description of this map and show that it is injective by Birman-Hilden theory. This gives us a new interesting non-geometric embedding of braid group into mapping class group. On the other hand, we show that the map on the level of classifying spaces of groups is compatible with the action of little 2-cube operad so that it induces a trivial homomorphim between stable homology group of braid groups and that of mapping class groups(Harer conjecture). We also show how the lift $\tilde{\beta_i}$ acts on the fundamental group of the surface and through this we prove that $\tilde{\beta_i}$ equals the product of two inverse Dehn twists.

\end{abstract}

\section{Introduction}

\begin{figure}[H]
  \centering

  \definecolor{c1061a3}{RGB}{16,97,163}

  \begin{tikzpicture}[y=0.80pt, x=0.80pt, yscale=-1.000000, xscale=1.000000, inner sep=0pt, outer sep=0pt]
  \path[draw=black,line join=round,line cap=round,line width=0.800pt]
    (110.5650,762.8980) .. controls (110.5650,728.2960) and (180.6920,700.2450) ..
    (267.1990,700.2450) .. controls (353.7050,700.2450) and (423.8320,728.2960) ..
    (423.8320,762.8980) .. controls (423.8320,797.5010) and (353.7050,825.5520) ..
    (267.1990,825.5520) .. controls (180.6920,825.5520) and (110.5650,797.5010) ..
    (110.5650,762.8980) -- cycle;
  \path[draw=black,line join=round,line cap=round,line width=0.800pt]
    (154.7910,762.8980) .. controls (154.7910,752.7210) and (163.0410,744.4710) ..
    (173.2180,744.4710) .. controls (183.3960,744.4710) and (191.6460,752.7210) ..
    (191.6460,762.8980) .. controls (191.6460,773.0760) and (183.3960,781.3260) ..
    (173.2180,781.3260) .. controls (163.0410,781.3260) and (154.7910,773.0760) ..
    (154.7910,762.8980) -- cycle;
  \path[draw=black,line join=round,line cap=round,line width=0.800pt]
    (235.8720,762.8980) .. controls (235.8720,752.7210) and (244.1220,744.4710) ..
    (254.2990,744.4710) .. controls (264.4770,744.4710) and (272.7270,752.7210) ..
    (272.7270,762.8980) .. controls (272.7270,773.0760) and (264.4770,781.3260) ..
    (254.2990,781.3260) .. controls (244.1220,781.3260) and (235.8720,773.0760) ..
    (235.8720,762.8980) -- cycle;
  \path[draw=black,line join=round,line cap=round,line width=0.800pt]
    (342.7510,762.8980) .. controls (342.7510,752.7210) and (351.0020,744.4710) ..
    (361.1790,744.4710) .. controls (371.3560,744.4710) and (379.6070,752.7210) ..
    (379.6070,762.8980) .. controls (379.6070,773.0760) and (371.3560,781.3260) ..
    (361.1790,781.3260) .. controls (351.0020,781.3260) and (342.7510,773.0760) ..
    (342.7510,762.8980) -- cycle;
  \path[draw=black,fill=black,line join=round,line cap=round,line width=0.800pt]
    (295.0940,762.8980) .. controls (295.0940,762.0740) and (295.7630,761.4050) ..
    (296.5870,761.4050) .. controls (297.4120,761.4050) and (298.0810,762.0740) ..
    (298.0810,762.8980) .. controls (298.0810,763.7230) and (297.4120,764.3920) ..
    (296.5870,764.3920) .. controls (295.7630,764.3920) and (295.0940,763.7230) ..
    (295.0940,762.8980) -- cycle;
  \path[draw=black,fill=black,line join=round,line cap=round,line width=0.800pt]
    (304.0550,762.8980) .. controls (304.0550,762.0740) and (304.7230,761.4050) ..
    (305.5480,761.4050) .. controls (306.3730,761.4050) and (307.0420,762.0740) ..
    (307.0420,762.8980) .. controls (307.0420,763.7230) and (306.3730,764.3920) ..
    (305.5480,764.3920) .. controls (304.7230,764.3920) and (304.0550,763.7230) ..
    (304.0550,762.8980) -- cycle;
  \path[draw=black,fill=black,line join=round,line cap=round,line width=0.800pt]
    (313.0160,762.8980) .. controls (313.0160,762.0740) and (313.6840,761.4050) ..
    (314.5090,761.4050) .. controls (315.3340,761.4050) and (316.0030,762.0740) ..
    (316.0030,762.8980) .. controls (316.0030,763.7230) and (315.3340,764.3920) ..
    (314.5090,764.3920) .. controls (313.6840,764.3920) and (313.0160,763.7230) ..
    (313.0160,762.8980) -- cycle;
  \path[draw=c1061a3,line join=round,line cap=round,line width=0.800pt]
    (110.5650,762.8980) .. controls (110.5650,777.4990) and (154.7910,777.7560) ..
    (154.7910,762.8980);
  \path[draw=c1061a3,dash pattern=on 1.60pt,line join=round,line cap=round,line
    width=0.800pt] (154.7910,764.7410) .. controls (154.7910,750.1410) and
    (110.5650,749.8840) .. (110.5650,764.7410);
  \path[draw=c1061a3,line join=round,line cap=round,line width=0.800pt]
    (150.1820,762.8980) .. controls (150.1820,750.1760) and (160.4960,739.8620) ..
    (173.2180,739.8620) .. controls (185.9410,739.8620) and (196.2550,750.1760) ..
    (196.2550,762.8980) .. controls (196.2550,775.6210) and (185.9410,785.9350) ..
    (173.2180,785.9350) .. controls (160.4960,785.9350) and (150.1820,775.6210) ..
    (150.1820,762.8980) -- cycle;
  \path[draw=c1061a3,line join=round,line cap=round,line width=0.800pt]
    (191.6460,761.9770) .. controls (191.6460,776.5780) and (235.8720,776.8340) ..
    (235.8720,761.9770);
  \path[draw=c1061a3,dash pattern=on 1.60pt,line join=round,line cap=round,line
    width=0.800pt] (235.8720,763.8200) .. controls (235.8720,749.2190) and
    (191.6460,748.9630) .. (191.6460,763.8200);
  \path[draw=c1061a3,line join=round,line cap=round,line width=0.800pt]
    (231.2630,762.8980) .. controls (231.2630,750.1760) and (241.5770,739.8620) ..
    (254.2990,739.8620) .. controls (267.0220,739.8620) and (277.3360,750.1760) ..
    (277.3360,762.8980) .. controls (277.3360,775.6210) and (267.0220,785.9350) ..
    (254.2990,785.9350) .. controls (241.5770,785.9350) and (231.2630,775.6210) ..
    (231.2630,762.8980) -- cycle;
  \path[draw=c1061a3,line join=round,line cap=round,line width=0.800pt]
    (379.6070,761.9770) .. controls (379.6070,776.5780) and (423.8320,776.8340) ..
    (423.8320,761.9770);
  \path[draw=c1061a3,dash pattern=on 1.60pt,line join=round,line cap=round,line
    width=0.800pt] (423.8320,763.8200) .. controls (423.8320,749.2190) and
    (379.6070,748.9630) .. (379.6070,763.8200);
  \path[draw=c1061a3,line join=round,line cap=round,line width=0.800pt]
    (338.1430,762.8980) .. controls (338.1430,750.1760) and (348.4560,739.8620) ..
    (361.1790,739.8620) .. controls (373.9020,739.8620) and (384.2150,750.1760) ..
    (384.2150,762.8980) .. controls (384.2150,775.6210) and (373.9020,785.9350) ..
    (361.1790,785.9350) .. controls (348.4560,785.9350) and (338.1430,775.6210) ..
    (338.1430,762.8980) -- cycle;
  \path[fill=black,line width=0.600pt] (124.2653,784.2562) node[above right]
    (text4912) {$\alpha_1$};
  \path[fill=black,line width=0.600pt] (164.6333,797.6311) node[above right]
    (text4916) {$\alpha_2$};
  \path[fill=black,line width=0.600pt] (205.4876,784.0129) node[above right]
    (text4920) {$\alpha_3$};
  \path[fill=black,line width=0.600pt] (247.3147,797.6311) node[above right]
    (text4924) {$\alpha_4$};
  \path[fill=black,line width=0.600pt] (352.5000,796.8899) node[above right]
    (text4930) {$\alpha_{2g-2}$};
  \path[fill=black,line width=0.600pt] (393.7500,784.9104) node[above right]
    (text4934) {$\alpha_{2g-1}$};

  \end{tikzpicture}

  \caption{Dehn twists.}
  \label{fig:stdDehn}
\end{figure}
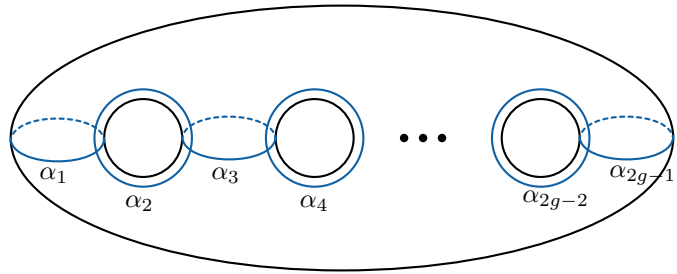

There is a well-known embedding of braid group $B_{2g}$ into mapping class group $\Gamma_g$ which maps the standard generators $\beta_i 's$ of $B_{2g}$ to consecutive standard Dehn twists $\alpha_i 's$ as in $\Gamma_g$ (see Figure \ref{fig:stdDehn}).
It was conjectured by J. Harer in 1960's that the homomorphism $\phi_*:H_*(B_\infty;\mathbb Z/2)\rightarrow H_*(\Gamma_\infty;\mathbb Z/2)$ induced by this embedding is trivial, where $B_\infty=\varinjlim B_n$ and $\Gamma_\infty=\varinjlim\Gamma_g$.
This conjecture can be proved by showing that the homomorphism $\phi_*$ preserves the Araki-Kudo-Dyer-Lashof operations which arises from the double loop space structures (of the group completion of classifying spaces of these groups), hence this problem is related with the theory of iterated loop space.
This conjecture was proved by Song-Tillmann (\cite{ST07}).

In the classical coherence theory of iterated loop spaces, there are basically two approaches.
The first approach is related with the coherence problem in category theory.
A classical result is that the group completion of classifying space of a monoidal (symmetric monoidal, braided monoidal, respectively) category is homotopy equivalent to a 1-fold loop space (infinite loop space, 2-fold loop space, respectively)  (\cite{Stasheff63},\cite{Segal74},\cite{MacLane65}).
A rather recent result is that the group completion of the classifying space of $n$-fold monoidal category is homotopy equivalent to an $n$-fold loop space (although there are several versions of definition of iterated category) (\cite{BFSV03}). The second approach is using the operad theory. If a space has a suitable action of little $n$-cube operad, then it gives rise to an $n$-fold loop space. These two approaches are not independent, rather mutually related in a natural way.

The proof of Harer conjecture of Song-Tillmann took the first approach, They lifted the embedding $\phi:B_{2g}\rightarrow\Gamma_g$ to a 2-functor between two monoidal 2-categories, tile category and surface category which are a kind of categorifications of braid groups and mapping class groups, respectively.

On the other hand, Segal and Tillmann suggested another proof of Harer conjecture taking the second approach to coherence theory of iterated loop spaces.
They lifted the embedding $\phi:B_{2g+2}\rightarrow\Gamma_{g,2}$ to a map between the classifying spaces of them and showed that this space map is compatible with the action of the framed little 2-cube operad on two spaces.

More precisely, let $\C_n=\operatorname{Conf}_n(D)/\Sigma_n$ be the configuration space of unordered $n$ distinct points on the disk $D=\{z\in\mathbb C \mid |z|\leq 1\}$.
Let $P=\{p_1,\ldots,p_n\}\subset D$. Note that
$$\C_n\simeq \operatorname{B}\homeo(D\setminus P)\simeq\operatorname{B}B_n.$$
Here $\homeo(D\setminus P)$ means the group of self-homeomorphisms of $D\setminus P$ which fix the boundary of $D$ pointwise, allowing the points in $P$ to be permuted. Let $\mathcal M_{g,2}$ be the moduli space of Riemann surface $S_{g,2}$. Note that $$\mathcal M_{g,2}\simeq\operatorname{B}\homeo(S_{g,2})\simeq \operatorname{B}\Gamma_{g,2}.$$
Segal and Tillmann considered the map $$\Phi:\C_{2g+2}\rightarrow \mathcal M_{g,2}$$
which maps $P=\{p_1,\ldots,p_{2g+2}\}$ to the part of the surface $S_g$ of the function $f_P(z)=\left((z-p_1)\cdots(z-p_{2g+2})\right)^{1/2}$ which lies over the disk $D$.

\begin{figure}[H]
  \centering

  \definecolor{c2f3b73}{RGB}{47,59,115}
  \definecolor{c74b587}{RGB}{116,181,135}

  \begin{tikzpicture}[y=0.80pt, x=0.80pt, yscale=-1.000000, xscale=1.000000, inner sep=0pt, outer sep=0pt, scale=0.8]
  \path[draw=black,line join=round,line width=0.800pt] (162.1620,410.0160) --
    (162.1620,506.4790);
  \path[draw=black,line join=round,line width=0.800pt] (75.3157,410.0160) --
    (75.3157,506.4790);
  \path[draw=black,line join=round,line cap=round,line width=0.800pt]
    (139.9630,453.0750) .. controls (139.9630,452.0890) and (140.7620,451.2890) ..
    (141.7490,451.2890) .. controls (142.7350,451.2890) and (143.5350,452.0890) ..
    (143.5350,453.0750) .. controls (143.5350,454.0620) and (142.7350,454.8610) ..
    (141.7490,454.8610) .. controls (140.7620,454.8610) and (139.9630,454.0620) ..
    (139.9630,453.0750) -- cycle;
  \path[draw=black,line join=round,line cap=round,line width=0.800pt]
    (93.5253,466.6430) .. controls (93.5253,465.6560) and (94.3250,464.8570) ..
    (95.3114,464.8570) .. controls (96.2978,464.8570) and (97.0974,465.6560) ..
    (97.0974,466.6430) .. controls (97.0974,467.6290) and (96.2978,468.4290) ..
    (95.3114,468.4290) .. controls (94.3250,468.4290) and (93.5253,467.6290) ..
    (93.5253,466.6430) -- cycle;
  \path[draw=black,line join=round,line width=0.800pt] (75.3157,410.0160) ..
    controls (75.3157,405.4340) and (94.7569,401.7190) .. (118.7390,401.7190) ..
    controls (142.7210,401.7190) and (162.1620,405.4340) .. (162.1620,410.0160) ..
    controls (162.1620,414.5980) and (142.7210,418.3120) .. (118.7390,418.3120) ..
    controls (94.7569,418.3120) and (75.3157,414.5980) .. (75.3157,410.0160) --
    cycle;
  \path[draw=black,dash pattern=on 1.60pt,line join=round,line cap=round,line
    width=0.800pt] (75.3157,506.4790) .. controls (75.3157,501.8970) and
    (94.7569,498.1830) .. (118.7390,498.1830) .. controls (142.7210,498.1830) and
    (162.1620,501.8970) .. (162.1620,506.4790);
  \path[draw=black,line join=round,line width=0.800pt] (162.1620,506.4790) ..
    controls (162.1620,511.0610) and (142.7210,514.7750) .. (118.7390,514.7750) ..
    controls (94.7569,514.7750) and (75.3157,511.0610) .. (75.3157,506.4790);
      \path[draw=c2f3b73,dash pattern=on 1.60pt,line join=round,line cap=round,line
        width=0.800pt] (136.3370,452.5460) .. controls (130.9580,452.0890) and
        (125.0040,451.8360) .. (118.7390,451.8360) .. controls (94.7569,451.8360) and
        (75.3157,455.5500) .. (75.3157,460.1320);
      \path[fill=c2f3b73] (135.4660,452.4630) -- (132.5200,455.6980) --
        (139.8210,452.8780) -- (133.1850,448.7290) -- (135.4660,452.4630) -- cycle;
  \path[draw=c2f3b73,line join=round,line cap=round,line width=0.800pt]
    (93.2432,466.8490) .. controls (82.3793,465.3410) and (75.3157,462.8940) ..
    (75.3157,460.1320);
  \path[draw=c2f3b73,line join=round,line cap=round,line width=0.800pt]
    (162.1620,460.1320) .. controls (162.1620,464.7140) and (142.7210,468.4290) ..
    (118.7390,468.4290) .. controls (110.7430,468.4290) and (103.2510,468.0160) ..
    (96.8172,467.2950);
      \path[draw=c2f3b73,dash pattern=on 1.60pt,line join=round,line cap=round,line
        width=0.800pt] (147.2320,453.8720) .. controls (156.3790,455.3920) and
        (162.1620,457.6330) .. (162.1620,460.1320);
      \path[fill=c2f3b73] (148.0970,454.0020) -- (151.2130,450.9300) --
        (143.7710,453.3520) -- (150.1740,457.8520) -- (148.0970,454.0020) -- cycle;
  \path[draw=c74b587,line join=round,draw opacity=0.798,line width=1.600pt]
    (95.3114,468.4290) -- (95.3596,513.4690);
  \path[draw=c74b587,line join=round,draw opacity=0.798,line width=1.600pt]
    (95.2766,417.0010) -- (95.3114,464.8570);
      \path[draw=black,line join=round,line cap=round,line width=0.800pt]
        (227.2090,460.4460) .. controls (223.9310,460.4640) and (183.2360,460.6870) ..
        (183.2360,460.6870);
      \path[fill=black] (224.3780,463.8480) -- (228.0280,460.8280) --
        (228.4990,460.4390) -- (228.0240,460.0540) -- (224.3400,457.0750) .. controls
        (224.1260,456.9010) and (223.8110,456.9340) .. (223.6370,457.1490) .. controls
        (223.4640,457.3630) and (223.4970,457.6780) .. (223.7110,457.8520) --
        (227.3950,460.8320) -- (227.3900,460.0580) -- (223.7400,463.0780) .. controls
        (223.5270,463.2540) and (223.4980,463.5690) .. (223.6740,463.7820) .. controls
        (223.8500,463.9950) and (224.1650,464.0240) .. (224.3780,463.8480) -- cycle;
  \path[draw=black,line join=round,line cap=round,line width=0.800pt]
    (250.5570,461.4530) .. controls (250.5570,442.7110) and (276.1460,427.5180) ..
    (307.7110,427.5180) .. controls (339.2760,427.5180) and (364.8640,442.7110) ..
    (364.8640,461.4530) .. controls (364.8640,480.1950) and (339.2760,495.3880) ..
    (307.7110,495.3880) .. controls (276.1460,495.3880) and (250.5570,480.1950) ..
    (250.5570,461.4530) -- cycle;
  \path[draw=black,line join=round,line cap=round,line width=0.800pt]
    (282.7060,461.4530) .. controls (282.7060,460.4670) and (283.5060,459.6670) ..
    (284.4920,459.6670) .. controls (285.4790,459.6670) and (286.2780,460.4670) ..
    (286.2780,461.4530) .. controls (286.2780,462.4390) and (285.4790,463.2390) ..
    (284.4920,463.2390) .. controls (283.5060,463.2390) and (282.7060,462.4390) ..
    (282.7060,461.4530) -- cycle;
  \path[draw=black,line join=round,line cap=round,line width=0.800pt]
    (329.1430,461.4530) .. controls (329.1430,460.4670) and (329.9430,459.6670) ..
    (330.9300,459.6670) .. controls (331.9160,459.6670) and (332.7160,460.4670) ..
    (332.7160,461.4530) .. controls (332.7160,462.4390) and (331.9160,463.2390) ..
    (330.9300,463.2390) .. controls (329.9430,463.2390) and (329.1430,462.4390) ..
    (329.1430,461.4530) -- cycle;
      \path[draw=c2f3b73,line join=round,line cap=round,line width=0.800pt]
        (325.6420,461.4530) .. controls (321.6180,461.4530) and (286.2780,461.4530) ..
        (286.2780,461.4530);
      \path[fill=c2f3b73] (324.7670,461.4530) -- (322.1420,464.9530) --
        (329.1420,461.4530) -- (322.1420,457.9530) -- (324.7670,461.4530) -- cycle;
  \path[draw=c74b587,line join=round,draw opacity=0.798,line width=1.600pt]
    (250.5570,461.4530) -- (282.7060,461.4530);
  \path[fill=black,line width=0.600pt] (133.2649,446.0504) node[above right]
    (text2677) {$p_{i+1}$};
  \path[fill=black,line width=0.600pt] (202.3906,455.6798) node[above right]
    (text2681) {$p$};
  \path[fill=black,line width=0.600pt] (79.9589,477.0022) node[above right]
    (text2685) {$p_i$};
  \path[fill=black,line width=0.600pt] (281.4898,455.8518) node[above right]
    (text2689) {$p_i$};
  \path[fill=black,line width=0.600pt] (327.2297,455.6798) node[above right]
    (text2693) {$p_{i+1}$};

  \end{tikzpicture}

  \caption{2-fold covering over $D_i$.}
  \label{fig:annulus}
\end{figure}
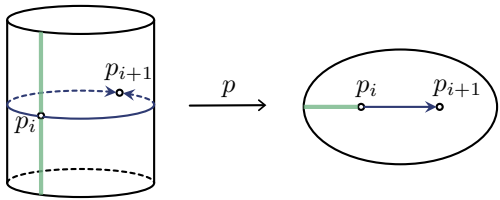

The homeomorphism $\phi=\pi_1\Phi$ on the fundamental groups induced by $\Phi$ is determined by the 2-fold branched covering $p:S_P\rightarrow D$ branched at $P=\{p_1,\ldots,p_{2g+2}\}$. Let $D_i$ be the subdisk of $D$ containing two points $p_i$ and $p_{i+1}$. Then $p^{-1}(D_i)$ is an annulus contained in $S_P$. (Figure \ref{fig:annulus})

\begin{figure}[H]
  \centering

  \definecolor{c2f3b73}{RGB}{47,59,115}
  \definecolor{c74b587}{RGB}{116,181,135}

  \begin{tikzpicture}[y=0.80pt, x=0.80pt, yscale=-1.000000, xscale=1.000000, inner sep=0pt, outer sep=0pt, scale=0.8]
  \path[draw=black,line join=round,line width=0.800pt] (165.5740,579.5490) --
    (165.5740,676.0120);
  \path[draw=black,line join=round,line width=0.800pt] (78.7279,579.5490) --
    (78.7279,676.0120);
  \path[draw=black,line join=round,line cap=round,line width=0.800pt]
    (143.3750,622.6080) .. controls (143.3750,621.6220) and (144.1750,620.8220) ..
    (145.1610,620.8220) .. controls (146.1470,620.8220) and (146.9470,621.6220) ..
    (146.9470,622.6080) .. controls (146.9470,623.5950) and (146.1470,624.3940) ..
    (145.1610,624.3940) .. controls (144.1750,624.3940) and (143.3750,623.5950) ..
    (143.3750,622.6080) -- cycle;
  \path[draw=black,line join=round,line cap=round,line width=0.800pt]
    (96.9376,636.1760) .. controls (96.9376,635.1890) and (97.7372,634.3900) ..
    (98.7236,634.3900) .. controls (99.7100,634.3900) and (100.5100,635.1890) ..
    (100.5100,636.1760) .. controls (100.5100,637.1620) and (99.7100,637.9620) ..
    (98.7236,637.9620) .. controls (97.7372,637.9620) and (96.9376,637.1620) ..
    (96.9376,636.1760) -- cycle;
  \path[draw=black,line join=round,line width=0.800pt] (78.7279,579.5490) ..
    controls (78.7279,574.9670) and (98.1692,571.2520) .. (122.1510,571.2520) ..
    controls (146.1330,571.2520) and (165.5740,574.9670) .. (165.5740,579.5490) ..
    controls (165.5740,584.1310) and (146.1330,587.8450) .. (122.1510,587.8450) ..
    controls (98.1692,587.8450) and (78.7279,584.1310) .. (78.7279,579.5490) --
    cycle;
  \path[draw=black,dash pattern=on 1.60pt,line join=round,line cap=round,line
    width=0.800pt] (78.7279,676.0120) .. controls (78.7279,671.4300) and
    (98.1692,667.7160) .. (122.1510,667.7160) .. controls (146.1330,667.7160) and
    (165.5740,671.4300) .. (165.5740,676.0120);
  \path[draw=black,line join=round,line width=0.800pt] (165.5740,676.0120) ..
    controls (165.5740,680.5940) and (146.1330,684.3080) .. (122.1510,684.3080) ..
    controls (98.1692,684.3080) and (78.7279,680.5940) .. (78.7279,676.0120);
      \path[draw=c2f3b73,dash pattern=on 1.60pt,line join=round,line cap=round,line
        width=0.800pt] (139.7490,622.0790) .. controls (134.3710,621.6220) and
        (128.4160,621.3690) .. (122.1510,621.3690) .. controls (98.1692,621.3690) and
        (78.7279,625.0830) .. (78.7279,629.6650);
      \path[fill=c2f3b73] (138.8780,621.9960) -- (135.9330,625.2310) --
        (143.2330,622.4110) -- (136.5970,618.2620) -- (138.8780,621.9960) -- cycle;
  \path[draw=c2f3b73,line join=round,line cap=round,line width=0.800pt]
    (96.6555,636.3820) .. controls (85.7915,634.8740) and (78.7279,632.4270) ..
    (78.7279,629.6650);
  \path[draw=c2f3b73,line join=round,line cap=round,line width=0.800pt]
    (165.5740,629.6650) .. controls (165.5740,634.2470) and (146.1330,637.9620) ..
    (122.1510,637.9620) .. controls (114.1550,637.9620) and (106.6640,637.5490) ..
    (100.2290,636.8280);
      \path[draw=c2f3b73,dash pattern=on 1.60pt,line join=round,line cap=round,line
        width=0.800pt] (150.6440,623.4050) .. controls (159.7910,624.9250) and
        (165.5740,627.1660) .. (165.5740,629.6650);
      \path[fill=c2f3b73] (151.5100,623.5350) -- (154.6250,620.4630) --
        (147.1830,622.8850) -- (153.5860,627.3850) -- (151.5100,623.5350) -- cycle;
  \path[draw=c74b587,line join=round,draw opacity=0.798,line width=1.600pt]
    (98.7236,637.9620) -- (98.7719,683.0020);
  \path[draw=c74b587,line join=round,draw opacity=0.798,line width=1.600pt]
    (98.6888,586.5340) -- (98.7236,634.3900);
      \path[draw=black,line join=round,line cap=round,line width=0.800pt]
        (224.5950,625.3520) .. controls (221.3180,625.3700) and (180.6220,625.5940) ..
        (180.6220,625.5940);
      \path[fill=black] (221.7640,628.7550) -- (225.4140,625.7340) --
        (225.8850,625.3450) -- (225.4100,624.9600) -- (221.7270,621.9810) .. controls
        (221.5120,621.8070) and (221.1970,621.8400) .. (221.0230,622.0550) .. controls
        (220.8500,622.2700) and (220.8830,622.5840) .. (221.0980,622.7580) --
        (224.7810,625.7380) -- (224.7770,624.9640) -- (221.1260,627.9840) .. controls
        (220.9140,628.1600) and (220.8840,628.4750) .. (221.0600,628.6880) .. controls
        (221.2360,628.9010) and (221.5510,628.9310) .. (221.7640,628.7550) -- cycle;
  \path[draw=black,line join=round,line width=0.800pt] (327.5320,579.5490) --
    (327.5320,676.0120);
  \path[draw=black,line join=round,line width=0.800pt] (240.6860,579.5490) --
    (240.6860,676.0120);
  \path[draw=black,line join=round,line cap=round,line width=0.800pt]
    (305.3320,622.6080) .. controls (305.3320,621.6220) and (306.1320,620.8220) ..
    (307.1180,620.8220) .. controls (308.1050,620.8220) and (308.9050,621.6220) ..
    (308.9050,622.6080) .. controls (308.9050,623.5950) and (308.1050,624.3940) ..
    (307.1180,624.3940) .. controls (306.1320,624.3940) and (305.3320,623.5950) ..
    (305.3320,622.6080) -- cycle;
  \path[draw=black,line join=round,line cap=round,line width=0.800pt]
    (258.8950,636.1760) .. controls (258.8950,635.1890) and (259.6950,634.3900) ..
    (260.6810,634.3900) .. controls (261.6680,634.3900) and (262.4670,635.1890) ..
    (262.4670,636.1760) .. controls (262.4670,637.1620) and (261.6680,637.9620) ..
    (260.6810,637.9620) .. controls (259.6950,637.9620) and (258.8950,637.1620) ..
    (258.8950,636.1760) -- cycle;
  \path[draw=black,line join=round,line width=0.800pt] (240.6860,579.5490) ..
    controls (240.6860,574.9670) and (260.1270,571.2530) .. (284.1090,571.2530) ..
    controls (308.0910,571.2530) and (327.5320,574.9670) .. (327.5320,579.5490) ..
    controls (327.5320,584.1310) and (308.0910,587.8450) .. (284.1090,587.8450) ..
    controls (260.1270,587.8450) and (240.6860,584.1310) .. (240.6860,579.5490) --
    cycle;
  \path[draw=black,dash pattern=on 1.60pt,line join=round,line cap=round,line
    width=0.800pt] (240.6860,676.0120) .. controls (240.6860,671.4300) and
    (260.1270,667.7160) .. (284.1090,667.7160) .. controls (308.0910,667.7160) and
    (327.5320,671.4300) .. (327.5320,676.0120);
  \path[draw=black,line join=round,line width=0.800pt] (327.5320,676.0120) ..
    controls (327.5320,680.5940) and (308.0910,684.3080) .. (284.1090,684.3080) ..
    controls (260.1270,684.3080) and (240.6860,680.5940) .. (240.6860,676.0120);
  \path[draw=c2f3b73,dash pattern=on 1.60pt,line join=round,line cap=round,line
    width=0.800pt] (240.6860,629.6650) .. controls (240.6860,625.0830) and
    (260.1270,621.3690) .. (284.1090,621.3690) .. controls (291.7600,621.3690) and
    (298.9490,621.7470) .. (305.1930,622.4110);
      \path[draw=c2f3b73,line join=round,line cap=round,line width=0.800pt]
        (255.1520,635.8480) .. controls (246.2740,634.3290) and (240.6860,632.1220) ..
        (240.6860,629.6650);
      \path[fill=c2f3b73] (254.2870,635.7150) -- (251.1600,638.7730) --
        (258.6110,636.3820) -- (252.2270,631.8550) -- (254.2870,635.7150) -- cycle;
      \path[draw=c2f3b73,line join=round,line cap=round,line width=0.800pt]
        (265.6710,637.1790) .. controls (271.2670,637.6810) and (277.5160,637.9620) ..
        (284.1090,637.9620) .. controls (308.0910,637.9620) and (327.5320,634.2470) ..
        (327.5320,629.6650);
      \path[fill=c2f3b73] (266.5420,637.2670) -- (269.5040,634.0470) --
        (262.1890,636.8290) -- (268.8030,641.0120) -- (266.5420,637.2670) -- cycle;
  \path[draw=c2f3b73,dash pattern=on 1.60pt,line join=round,line cap=round,line
    width=0.800pt] (327.5320,629.6650) .. controls (327.5320,626.8640) and
    (320.2660,624.3870) .. (309.1390,622.8850);
  \path[draw=c74b587,dash pattern=on 1.60pt,line join=round,draw
    opacity=0.798,line width=1.600pt] (307.1180,620.8220) -- (307.3420,596.0820)
    .. controls (290.0710,591.0930) and (240.6860,593.1720) ..
    (240.6860,601.0850);
  \path[draw=c74b587,line join=round,draw opacity=0.798,line width=1.600pt]
    (260.6460,586.5340) -- (260.6600,605.6080) .. controls (260.6600,605.6080) and
    (240.6860,606.2050) .. (240.6860,601.0850);
  \path[draw=c74b587,line join=round,draw opacity=0.799,line width=1.600pt]
    (327.5280,652.4740) .. controls (327.5280,665.5900) and (260.6600,662.3280) ..
    (260.6600,662.3280) -- (260.4720,682.9700);
  \path[draw=c74b587,dash pattern=on 1.60pt,line join=round,draw
    opacity=0.799,line width=1.600pt] (327.5280,652.4740) .. controls
    (327.5280,647.0970) and (307.3380,644.1780) .. (307.3380,644.1780) --
    (307.2250,623.9170);
  \path[fill=black,line width=0.600pt] (195.0000,616.8900) node[above right]
    (text3687) {$\tilde\beta_i$};

  \end{tikzpicture}

  \caption{The lift of half Dehn twist is a full Dehn twist.}
  \label{fig:lift_alpha}
\end{figure}
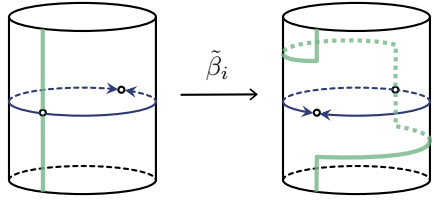

Let $\beta_i$ be a generator of $B_{2g+2}$ interchanging the $i$-th and $(i+1)$-st points.
Then $\phi:\pi_1\C_{2g+2}\rightarrow \pi_1\mathcal M_{g,2}$ maps $\beta_i$ to the Dehn twists $\alpha_i$ on the annulus $p^{-1}(D_i)$. $\alpha_i$ is described in Figure~\ref{fig:lift_alpha}.

The map $\phi$ is the well-known (Harer) embedding of braid group into mapping class group.It was shown in \cite{Segal-Tillmann} that the map $\Phi:\C_{2g+2}\rightarrow \mathcal M_{g,2}$ is compatible with the actions of the framed little 2-cube operad on configuration spaces and moduli spaces.

The main problem of this paper is an extension of that of Segal-Tillmann: what can we get in the case of 3-fold covering, instead of 2-fold covering?  This problem was, a  few years ago, suggested by Tillmann to the second author. In the case of 3-fold covering, there are mysterious difficulty and complexity. Without careful geometric observations and calculations one may easily fall into a wrong path.

Let $\Phi:\C_k\rightarrow\M_{g,b}$ map $P=\{p_1,\ldots,p_k\}$ to the part of $S_{g}$ of the function $f_P(z)=\left((z-p_1)\cdots(z-p_g)\right)^{1/3}$.
For the map $\phi:B_k\rightarrow\Gamma_{g,b}$ induced by $\Phi$, we may have the following questions.

\begin{enumerate}
  \item Over a disk with two marked points, what is the ``suitable'' 3-fold branched covering space with those two branch points?
  \item For the half Dehn twist on a disk with two marked points, what is the lift of it with respect to the 3-fold covering?
  \item In the case of more branch points, do two adjacent half Dehn twists have the lifts satisfying the braid relation so that $\Phi:\C_k\rightarrow \M_{g,b}$ induces a well-defined homomorphism $\phi:B_k\rightarrow\Gamma_{g,b}$?
  \item For $\Phi:\C_k\rightarrow\M_{g,b}$, what are $g$ and $b$ in terms of $k$?
  \item Is map $\phi:B_k\rightarrow\Gamma_{g,b}$ injective? If so, is it geometric or non-geometric? A map is said to be {\it geometric} if it sends each generator $\beta_i$ of $B_k$ to a Dehn twist.
  \item The homomorphism $\phi_*:H_*(B_\infty;R)\rightarrow H_*(\Gamma_\infty;R)$ induced by $\phi$ is zero for any constant coefficient $R$?
  \item For each lift $\tilde{\beta_i}$ of $\beta_i$, what about the action of $\tilde{\beta_i}$ on the fundamental group of the surface?

\end{enumerate}

For question 1, the suitable 3-fold branched covering over a disk with two branch points turns out to be a torus with a deleted disk (with one boundary component).

\begin{figure}[H]
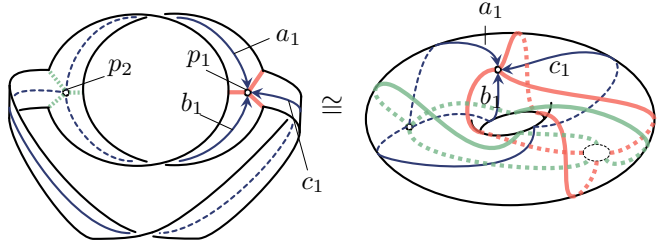

  \centering

  \definecolor{c2f3b73}{RGB}{47,59,115}
  \definecolor{cf56356}{RGB}{245,99,86}
  \definecolor{c74b587}{RGB}{116,181,135}



  \caption{3-fold branched covering over a disk with two branch points.}
  \label{fig:3f2b}
\end{figure}

The 3-fold covering over a disk with three branch points is shown to be a torus with three boundary components.

\begin{figure}[H]
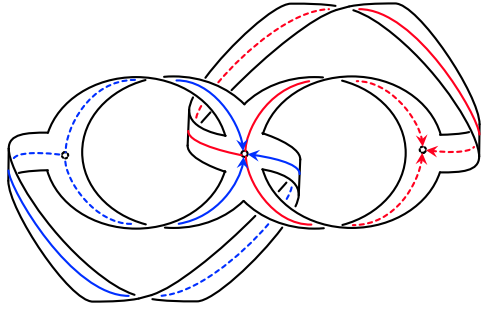

  \centering

  \definecolor{cff001f}{RGB}{255,0,31}
  \definecolor{c002fff}{RGB}{0,47,255}

  %

  \caption{3-fold covering over a disk with three branch points.}
  \label{fig:3f3b}
\end{figure}

For question 2, the lift $\tilde{\beta_i}$ of half Dehn twist $\beta_i$ on a disk turns out to be a kind of $1/6$ Dehn twist.

\begin{figure}[H]
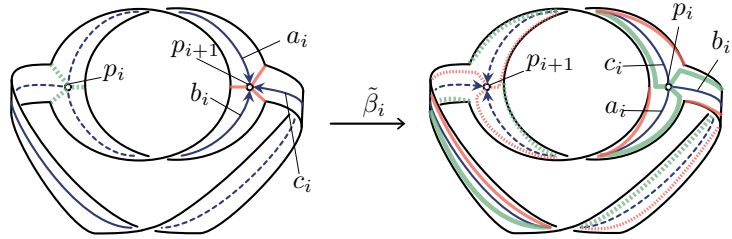

  \centering

  \definecolor{c2f3b73}{RGB}{47,59,115}
  \definecolor{cf56356}{RGB}{245,99,86}
  \definecolor{c74b587}{RGB}{116,181,135}

  %

  \caption{1/6 Dehn twist.}
\end{figure}

The answer to question 3 is affirmative which is proved in a few ways.

For question 4, we have
$$g=\begin{cases}
  k-2 & \mbox{if }g\equiv 0\quad(\mbox{mod }3) \\
  k-1 & \mbox{otherwise}
\end{cases}
$$
$$
b = \begin{cases}
  3 & \mbox{if }g\equiv 0\quad(\mbox{mod }3) \\
  1 & \mbox{otherwise}
\end{cases}
$$

From an intuitive geometric observation, it is easy to see that $\phi:B_k\rightarrow\Gamma_{g,b}$ is well-defined, that is, the braid relation
 $$\tilde{\beta_i}\tilde{\beta}_{i+1}\tilde{\beta_i}=\tilde{\beta}_{i+1}\tilde{\beta_i}\tilde{\beta}_{i+1}$$
is satisfied. By applying the Birman-Hilden theory, we can easily show that the map $\phi:B_k\rightarrow\Gamma_{g,b}$ is injective(Theorem 3.3).

In subsection 3.2, we show that $\Phi:\C_k\rightarrow \M_{g,b}$ is compatible with naturally defined actions of the framed little 2-cube operad. This implies that the homomorphism $\phi:B_k\hookrightarrow\Gamma_{g,b}\rightarrow\Gamma_{g}$ induces trivial homomorphism in stable homology with any constant coefficient(Theorem 3.4).

In section 4 by analysing the action of the lift $\tilde{\beta_i}$ on the fundamental group of the surface, we show that $\tilde{\beta_i}$ is the product(composition) of two (usual) Dehn twists along closed curves on the surface. That is, the $1/6$ Dehn twist turns out to be the product of two full Dehn twists.

Since the question of existence of non-geometric embedding of braid group into mapping class group was raised by Wajnryb (\cite{Waj06}), only a few examples have been found (\cite{Sz10},\cite{BT12},\cite{JS13}), including so called pillar switchings.
The following figure shows the pillar switchings.

\begin{figure}[H]
  \centering

  \definecolor{c92b2de}{RGB}{146,178,222}

  \begin{tikzpicture}[y=0.80pt, x=0.80pt, yscale=-1.000000, xscale=1.000000, inner sep=0pt, outer sep=0pt, scale=0.8]
  \path[draw=black,fill=black,line join=round,line width=0.800pt]
    (404.2740,550.7120) .. controls (404.2740,549.9760) and (404.8870,549.3800) ..
    (405.6440,549.3800) .. controls (406.4000,549.3800) and (407.0130,549.9760) ..
    (407.0130,550.7120) .. controls (407.0130,551.4470) and (406.4000,552.0430) ..
    (405.6440,552.0430) .. controls (404.8870,552.0430) and (404.2740,551.4470) ..
    (404.2740,550.7120) -- cycle;
  \path[draw=black,fill=black,line join=round,line width=0.800pt]
    (414.7400,550.7120) .. controls (414.7400,549.9760) and (415.3530,549.3800) ..
    (416.1100,549.3800) .. controls (416.8660,549.3800) and (417.4790,549.9760) ..
    (417.4790,550.7120) .. controls (417.4790,551.4470) and (416.8660,552.0430) ..
    (416.1100,552.0430) .. controls (415.3530,552.0430) and (414.7400,551.4470) ..
    (414.7400,550.7120) -- cycle;
  \path[draw=black,fill=black,line join=round,line width=0.800pt]
    (424.9130,550.7120) .. controls (424.9130,549.9760) and (425.5260,549.3800) ..
    (426.2820,549.3800) .. controls (427.0380,549.3800) and (427.6520,549.9760) ..
    (427.6520,550.7120) .. controls (427.6520,551.4470) and (427.0380,552.0430) ..
    (426.2820,552.0430) .. controls (425.5260,552.0430) and (424.9130,551.4470) ..
    (424.9130,550.7120) -- cycle;
  \path[draw=black,dash pattern=on 1.60pt,line join=round,line width=0.800pt]
    (246.2090,538.2300) .. controls (246.2090,536.5380) and (255.2410,535.1670) ..
    (266.3820,535.1670) .. controls (277.5230,535.1670) and (286.5550,536.5380) ..
    (286.5550,538.2300);
  \path[draw=black,dash pattern=on 1.60pt,line join=round,line width=0.800pt]
    (246.2090,559.4770) .. controls (246.2090,557.7850) and (255.2410,556.4130) ..
    (266.3820,556.4130) .. controls (277.5230,556.4130) and (286.5550,557.7850) ..
    (286.5550,559.4770);
  \path[draw=black,dash pattern=on 1.60pt,line join=round,line width=0.800pt]
    (308.5890,538.2300) .. controls (308.5890,536.5380) and (317.6210,535.1670) ..
    (328.7620,535.1670) .. controls (339.9030,535.1670) and (348.9350,536.5380) ..
    (348.9350,538.2300);
  \path[draw=black,dash pattern=on 1.60pt,line join=round,line width=0.800pt]
    (308.5890,559.4770) .. controls (308.5890,557.7850) and (317.6210,556.4130) ..
    (328.7620,556.4130) .. controls (339.9030,556.4130) and (348.9350,557.7850) ..
    (348.9350,559.4770);
  \path[draw=black,line join=round,line width=0.800pt] (245.9720,559.5000) ..
    controls (243.1910,561.9920) and (239.4340,563.5210) .. (235.3000,563.5210) ..
    controls (226.7300,563.5210) and (219.7820,556.9500) .. (219.7820,548.8450) ..
    controls (219.7820,540.7400) and (226.7300,534.1690) .. (235.3000,534.1690) ..
    controls (239.5390,534.1690) and (243.3810,535.7760) .. (246.1810,538.3810);
  \path[draw=black,line join=round,line width=0.800pt] (308.4130,559.3710) ..
    controls (305.6180,561.9390) and (301.8040,563.5210) .. (297.6000,563.5210) ..
    controls (293.4040,563.5210) and (289.5970,561.9460) .. (286.8040,559.3870);
  \path[draw=black,line join=round,line width=0.800pt] (286.3590,538.7270) ..
    controls (289.1860,535.9190) and (293.1770,534.1690) .. (297.6000,534.1690) ..
    controls (302.0740,534.1690) and (306.1060,535.9600) .. (308.9380,538.8240);
  \path[draw=black,line join=round,line width=0.800pt] (286.8040,559.3870) ..
    controls (283.8920,556.7200) and (282.0820,552.9820) .. (282.0820,548.8450) ..
    controls (282.0820,544.9230) and (283.7090,541.3600) .. (286.3590,538.7270) ..
    controls (286.3590,540.4190) and (277.3270,541.7900) .. (266.1860,541.7900) ..
    controls (255.0450,541.7900) and (246.0140,540.4190) .. (246.0140,538.7270) ..
    controls (248.8760,541.3890) and (250.6500,545.0940) .. (250.6500,549.1900) ..
    controls (250.6500,553.3860) and (248.7890,557.1700) .. (245.8040,559.8450) ..
    controls (245.8040,561.5370) and (254.8360,562.9080) .. (265.9770,562.9080) ..
    controls (277.1180,562.9080) and (286.8040,561.0790) .. (286.8040,559.3870);
  \path[draw=black,line join=round,line width=0.800pt] (348.5690,538.8990) ..
    controls (346.0190,541.5140) and (344.4620,545.0070) .. (344.4620,548.8450) ..
    controls (344.4620,553.0210) and (346.3060,556.7890) .. (349.2660,559.4610) ..
    controls (349.2660,561.1530) and (340.2340,562.5250) .. (329.0930,562.5250) ..
    controls (317.9520,562.5250) and (308.9200,561.1530) .. (308.9200,559.4610) ..
    controls (311.8220,556.7950) and (313.6250,553.0640) .. (313.6250,548.9350) ..
    controls (313.6250,545.0610) and (312.0380,541.5380) .. (309.4450,538.9150) ..
    controls (309.4450,540.6070) and (318.4770,541.9780) .. (329.6180,541.9780) ..
    controls (340.7590,541.9780) and (348.5690,540.5900) .. (348.5690,538.8990);
  \path[draw=black,line join=round,line width=0.800pt] (348.5690,538.8990) ..
    controls (351.4050,535.9910) and (355.4670,534.1690) .. (359.9800,534.1690) ..
    controls (368.5500,534.1690) and (375.4980,540.7400) .. (375.4980,548.8450) ..
    controls (375.4980,556.9500) and (368.5500,563.5210) .. (359.9800,563.5210) ..
    controls (355.8250,563.5210) and (352.0510,561.9760) .. (349.2660,559.4610);
  \path[draw=black,line join=round,line width=0.800pt] (153.4280,548.8450) ..
    controls (153.4280,540.7400) and (160.3380,534.1690) .. (168.8630,534.1690) ..
    controls (177.3870,534.1690) and (184.2970,540.7400) .. (184.2970,548.8450) ..
    controls (184.2970,556.9500) and (177.3870,563.5210) .. (168.8630,563.5210) ..
    controls (160.3380,563.5210) and (153.4280,556.9500) .. (153.4280,548.8450) --
    cycle;
  \path[draw=c92b2de,dash pattern=on 1.60pt,line join=round,draw
    opacity=0.704,line width=2.400pt] (362.0330,534.1690) .. controls
    (362.0330,525.5890) and (333.1440,516.8670) .. (297.5080,516.8670) .. controls
    (261.8710,516.8670) and (232.9820,525.5180) .. (232.9820,534.1690);
  \path[draw=c92b2de,dash pattern=on 1.60pt,line join=round,draw
    opacity=0.704,line width=2.400pt] (232.9820,563.5210) .. controls
    (232.9820,568.2590) and (261.8710,572.1000) .. (297.5080,572.1000) .. controls
    (333.1440,572.1000) and (362.0330,568.2590) .. (362.0330,563.5210);
  \path[draw=black,line join=round,line width=0.800pt] (115.2620,548.8450) ..
    controls (115.2620,517.9460) and (197.2250,492.8970) .. (298.3320,492.8970) ..
    controls (399.4390,492.8970) and (481.4020,517.9460) .. (481.4020,548.8450) ..
    controls (481.4020,579.7440) and (399.4390,604.7920) .. (298.3320,604.7920) ..
    controls (197.2250,604.7920) and (115.2620,579.7440) .. (115.2620,548.8450) --
    cycle;
  \path[draw=c92b2de,line join=round,line cap=round,draw opacity=0.704,line
    width=2.400pt] (362.0330,563.5210) .. controls (362.0330,572.1000) and
    (333.1440,580.8230) .. (297.5080,580.8230) .. controls (261.8710,580.8230) and
    (232.9820,572.1720) .. (232.9820,563.5210)(294.6460,577.1370) --
    (302.0170,580.8230) -- (294.6460,584.5080);
  \path[draw=c92b2de,line join=round,line cap=round,draw opacity=0.704,line
    width=2.400pt] (232.9820,534.1690) .. controls (232.9820,529.4310) and
    (261.8710,525.5890) .. (297.5080,525.5890) .. controls (333.1440,525.5890) and
    (362.0330,529.4310) .. (362.0330,534.1690)(294.6460,521.8320) --
    (302.0170,525.5180) -- (294.6460,529.2030);

  \end{tikzpicture}

  \caption{Pillar switchings.}
\end{figure}
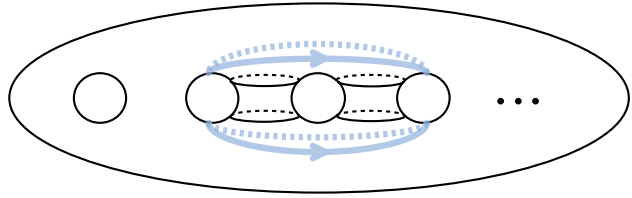

The map $\phi:B_k\rightarrow\Gamma_{g,b}$ obtained through 3-fold branched covering gives us a new interesting non-geometric embedding which is also homologically trivial.
\vskip 0.2cm
This research was partially fulfilled while the second author was staying at Dalian University of Technology(DLUT) for the Haitian Program. He is grateful Fengchun Lei for his collaboration and great hospitality. He appreciates the support from the mathematics department of DLUT. He is also grateful Ulrike Tillmann for letting him know the problem.

\section{3-fold branched covering spaces}

Let $S_{g,b,(k)}$ be the surface of genus $g$ with $b$ boundary components and set of $k$ distinct interior points $P=\{p_1,\ldots,p_k\}$.
Let $\homeo(S_{g,b,(k)})$ be the space of orientation preserving self-homemorphisms of $S_{g,b,(k)}$ which fix the boundary pointwise and preserve the set $P$.
Define $$\Gamma_{g,b,(k)}:=\pi_0\homeo(S_{g,b,(k)}),$$ called  mapping class group. Recall that braid group may be defined as a mapping class group : $B_k:=\Gamma_{0,1,(k)}$.

Let $\C_k$ be the configuration space which consists of unordered $k$-tuples of distinct points in the interior of the unit disk $D$.
Let $\M_{g,b}$ be the moduli space of connected Riemann surfaces of genus $g$ with $b$ boundary components.
The map $$\Phi:\C_k\rightarrow\M_{g,b}$$ sends $P=\{p_1,\ldots,p_k\}$ to the part of the Riemann surface $\Sigma_P$ of the function $f_P(z)=\left((z-p_1)\cdots(z-p_k)\right)^{1/3}$ which lies over the disk $D$.
Note that $\pi_1\C_k\cong B_k$ and $\pi_1\M_{g,b}\cong\Gamma_{g,b}$.

In order to determine the fundamental group homomorphism $\phi:B_k\rightarrow\Gamma_{g,b}$ induced by the map $\Phi$, we have to think of the 3-fold branched covering map $p:S_P\rightarrow D$ with the set of branch points $P=\{p_1,\ldots,p_k\}$.

\subsection{The case with two branch points}

Let us start with 2-fold branched covering over a disk with two branch points $\{p_1,p_2\}$.

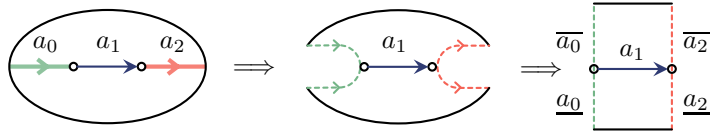
\begin{figure}[H]
  \centering

  \definecolor{cf56356}{RGB}{245,99,86}
  \definecolor{c74b587}{RGB}{116,181,135}
  \definecolor{c2f3b73}{RGB}{47,59,115}

  \begin{tikzpicture}[y=0.80pt, x=0.80pt, yscale=-1.000000, xscale=1.000000, inner sep=0pt, outer sep=0pt]
  \path[draw=cf56356,line join=round,draw opacity=0.804,line width=1.600pt]
    (105.7030,91.2613) -- (134.2800,91.2613)(116.7040,87.9392) --
    (121.0750,91.2613) -- (116.4200,94.7997);
  \path[draw=c74b587,line join=round,draw opacity=0.804,line width=1.600pt]
    (41.4056,91.2613) -- (69.9824,91.2613)(54.8175,94.6426) -- (59.2661,91.2613)
    -- (54.7279,87.8118);
  \path[draw=black,line join=round,line cap=round,line width=0.800pt]
    (69.9824,91.2613) .. controls (69.9824,90.2749) and (70.7820,89.4752) ..
    (71.7684,89.4752) .. controls (72.7548,89.4752) and (73.5545,90.2749) ..
    (73.5545,91.2613) .. controls (73.5545,92.2477) and (72.7548,93.0473) ..
    (71.7684,93.0473) .. controls (70.7820,93.0473) and (69.9824,92.2477) ..
    (69.9824,91.2613) -- cycle;
  \path[draw=black,line join=round,line cap=round,line width=0.800pt]
    (102.1310,91.2613) .. controls (102.1310,90.2749) and (102.9310,89.4752) ..
    (103.9170,89.4752) .. controls (104.9040,89.4752) and (105.7030,90.2749) ..
    (105.7030,91.2613) .. controls (105.7030,92.2477) and (104.9040,93.0473) ..
    (103.9170,93.0473) .. controls (102.9310,93.0473) and (102.1310,92.2477) ..
    (102.1310,91.2613) -- cycle;
      \path[draw=c2f3b73,line join=round,line cap=round,line width=0.800pt]
        (98.6286,91.2613) .. controls (95.8829,91.2613) and (73.5545,91.2613) ..
        (73.5545,91.2613);
      \path[fill=c2f3b73] (97.7536,91.2613) -- (95.1286,94.7613) -- (102.1290,91.2613)
        -- (95.1286,87.7613) -- (97.7536,91.2613) -- cycle;
  \path[draw=black,line join=round,line cap=round,line width=0.800pt]
    (268.2780,81.3538) .. controls (261.3930,71.5254) and (244.7720,64.5726) ..
    (225.2780,64.5726) .. controls (205.8090,64.5726) and (189.2100,71.5192) ..
    (182.3100,81.3226);
  \path[draw=cf56356,dash pattern=on 1.60pt,line join=round,line cap=round,line
    width=0.800pt] (268.3100,101.3540) -- (253.8410,101.3540) .. controls
    (248.3180,101.3540) and (243.8410,96.8767) .. (243.8410,91.3538) .. controls
    (243.8410,85.8310) and (248.3180,81.3538) .. (253.8410,81.3538) --
    (268.2780,81.3538);
  \path[draw=black,line join=round,line cap=round,line width=0.800pt]
    (182.2470,101.3230) .. controls (189.1070,111.1760) and (205.7460,118.1660) ..
    (225.2780,118.1660) .. controls (244.8110,118.1660) and (261.4500,111.2050) ..
    (268.3100,101.3540);
  \path[draw=c74b587,dash pattern=on 1.60pt,line join=round,line cap=round,line
    width=0.800pt] (182.3100,81.3226) .. controls (187.3930,81.2770) and
    (196.7890,81.1663) .. (197.3720,81.1663) .. controls (202.8950,81.1663) and
    (207.3720,85.6435) .. (207.3720,91.1663) .. controls (207.3720,96.6892) and
    (202.8950,101.1660) .. (197.3720,101.1660) .. controls (196.7020,101.1660) and
    (191.8100,101.2070) .. (188.3100,101.2600) .. controls (185.2960,101.3060) and
    (183.9490,101.3100) .. (182.2470,101.3230);
  \path[draw=black,line join=round,line cap=round,line width=0.800pt]
    (207.4250,91.3721) .. controls (207.4250,90.3857) and (208.2250,89.5860) ..
    (209.2110,89.5860) .. controls (210.1980,89.5860) and (210.9970,90.3857) ..
    (210.9970,91.3721) .. controls (210.9970,92.3585) and (210.1980,93.1581) ..
    (209.2110,93.1581) .. controls (208.2250,93.1581) and (207.4250,92.3585) ..
    (207.4250,91.3721) -- cycle;
  \path[draw=black,line join=round,line cap=round,line width=0.800pt]
    (239.5740,91.3721) .. controls (239.5740,90.3857) and (240.3740,89.5860) ..
    (241.3600,89.5860) .. controls (242.3460,89.5860) and (243.1460,90.3857) ..
    (243.1460,91.3721) .. controls (243.1460,92.3585) and (242.3460,93.1581) ..
    (241.3600,93.1581) .. controls (240.3740,93.1581) and (239.5740,92.3585) ..
    (239.5740,91.3721) -- cycle;
      \path[draw=c2f3b73,line join=round,line cap=round,line width=0.800pt]
        (236.0710,91.3721) .. controls (233.3260,91.3721) and (210.9970,91.3721) ..
        (210.9970,91.3721);
      \path[fill=c2f3b73] (235.1960,91.3721) -- (232.5710,94.8721) --
        (239.5710,91.3721) -- (232.5710,87.8721) -- (235.1960,91.3721) -- cycle;
  \path[draw=c74b587,line join=round,line cap=round,line width=0.800pt]
    (194.8230,78.8320) -- (197.3720,81.1663) -- (194.8460,83.7223);
  \path[draw=c74b587,line join=round,line cap=round,line width=0.800pt]
    (194.7790,98.7510) -- (197.3720,101.1810) -- (194.8730,103.7100);
  \path[draw=cf56356,line join=round,line cap=round,line width=0.800pt]
    (253.5270,78.9087) -- (256.0750,81.2430) -- (253.5490,83.7990);
  \path[draw=cf56356,line join=round,line cap=round,line width=0.800pt]
    (252.2520,98.7510) -- (254.8010,101.0850) -- (252.2750,103.6410);
  \path[draw=black,line join=round,line cap=round,line width=0.800pt]
    (41.4056,91.2613) .. controls (41.4056,76.4651) and (62.1963,64.4705) ..
    (87.8429,64.4705) .. controls (113.4890,64.4705) and (134.2800,76.4651) ..
    (134.2800,91.2613) .. controls (134.2800,106.0570) and (113.4890,118.0520) ..
    (87.8429,118.0520) .. controls (62.1963,118.0520) and (41.4056,106.0570) ..
    (41.4056,91.2613) -- cycle;
  \path[draw=cf56356,dash pattern=on 1.60pt,line join=round,line cap=round,line
    width=0.800pt] (354.6860,121.0080) -- (354.6860,92.3911) --
    (354.6860,61.4461);
  \path[draw=c74b587,dash pattern=on 1.60pt,line join=round,line cap=round,line
    width=0.800pt] (317.8310,121.0080) -- (317.8310,92.3911) --
    (317.8310,61.4461);
  \path[draw=black,line join=round,line cap=round,line width=0.800pt]
    (317.8310,61.4461) -- (354.6860,61.4461);
  \path[draw=black,line join=round,line cap=round,line width=0.800pt]
    (317.8310,121.0080) -- (354.6860,121.0080);
  \path[draw=black,line join=round,line cap=round,line width=0.800pt]
    (315.9880,92.3911) .. controls (315.9880,91.3734) and (316.8130,90.5483) ..
    (317.8310,90.5483) .. controls (318.8490,90.5483) and (319.6740,91.3734) ..
    (319.6740,92.3911) .. controls (319.6740,93.4088) and (318.8490,94.2338) ..
    (317.8310,94.2338) .. controls (316.8130,94.2338) and (315.9880,93.4088) ..
    (315.9880,92.3911) -- cycle;
  \path[draw=black,line join=round,line cap=round,line width=0.800pt]
    (352.8430,92.3911) .. controls (352.8430,91.3734) and (353.6680,90.5483) ..
    (354.6860,90.5483) .. controls (355.7040,90.5483) and (356.5290,91.3734) ..
    (356.5290,92.3911) .. controls (356.5290,93.4088) and (355.7040,94.2338) ..
    (354.6860,94.2338) .. controls (353.6680,94.2338) and (352.8430,93.4088) ..
    (352.8430,92.3911) -- cycle;
      \path[draw=c2f3b73,line join=round,line cap=round,line width=0.800pt]
        (349.3430,92.3911) .. controls (346.1780,92.3911) and (319.8870,92.3911) ..
        (319.8870,92.3911);
      \path[fill=c2f3b73] (348.4680,92.3911) -- (345.8430,95.8911) --
        (352.8430,92.3911) -- (345.8430,88.8911) -- (348.4680,92.3911) -- cycle;
  \path[fill=black,line width=0.600pt] (143.6938,94.9757) node[above right]
    (text5946) {$\implies$};
  \path[fill=black,line width=0.600pt] (279.3244,95.6063) node[above right]
    (text5950) {$\implies$};
  \path[fill=black,line width=0.600pt] (52.5000,84.3900) node[above right]
    (text5954) {$a_0$};
  \path[fill=black,line width=0.600pt] (82.5000,84.3900) node[above right]
    (text5958) {$a_1$};
  \path[fill=black,line width=0.600pt] (112.5000,84.3900) node[above right]
    (text5962) {$a_2$};
  \path[fill=black,line width=0.600pt] (300.0000,84.3900) node[above right]
    (text5970) {$\overline{a_0}$};
  \path[fill=black,line width=0.600pt] (300.0000,114.3900) node[above right]
    (text5974) {$\underline{a_0}$};
  \path[fill=black,line width=0.600pt] (360.0000,84.3900) node[above right]
    (text5978) {$\overline{a_2}$};
  \path[fill=black,line width=0.600pt] (360.0000,114.3900) node[above right]
    (text5982) {$\underline{a_2}$};
  \path[fill=black,line width=0.600pt] (217.5000,84.3900) node[above right]
    (text5986) {$a_1$};
  \path[fill=black,line width=0.600pt] (329.9212,88.1063) node[above right]
    (text5990) {$a_1$};

  \end{tikzpicture}

  \caption{Disk $a$ obtained by cutting the paths $a_0$ and $a_2$.}
  \label{fig:disk2pt}
\end{figure}

For a disk `$a$' with two branch points $p_1, p_2$, denote by $a_0, a_1$ and $a_2$ three paths in Figure \ref{fig:disk2pt}, one from the boundary to $p_1$, one from $p_1$ to $p_2$, and one from $p_2$ to the boundary, respectively. After being cut, the upper parts and the lower parts of the paths are denoted by overline and underline, respectively.

It is easy to see that the 2-fold branched covering space of a disk is an annulus, and the only covering transformation is the hyper-elliptic involution which interchanges the two boundary components. The 2-fold branched covering space can be constructed by cut-and-paste of two sheets of  disks; one is called $a$ and another one is called $b$. Cut the disks $a$ and $b$ along the paths $a_0, a_2, b_0, b_2$. Let  $\overline{a_{i}},\underline{a_{i}}$ and $\overline{b_{i}},\underline{b_{i}}$ be the corresponding upper paths and lower paths, for $i=0,2$. Here of course, the upper and lower path are identified to be the same after the paste.

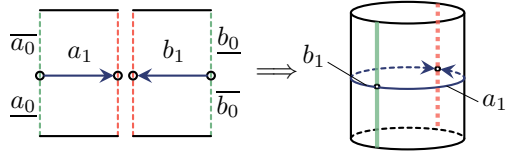
\begin{figure}[H]
  \centering

  \definecolor{cf56356}{RGB}{245,99,86}
  \definecolor{c74b587}{RGB}{116,181,135}
  \definecolor{c2f3b73}{RGB}{47,59,115}

  \begin{tikzpicture}[y=0.80pt, x=0.80pt, yscale=-1.000000, xscale=1.000000, inner sep=0pt, outer sep=0pt]
  \path[draw=cf56356,dash pattern=on 1.60pt,line join=round,line cap=round,line
    width=0.800pt] (80.8225,239.3730) -- (80.8225,210.7560) -- (80.8225,179.8110);
  \path[draw=c74b587,dash pattern=on 1.60pt,line join=round,line cap=round,line
    width=0.800pt] (43.9675,239.3730) -- (43.9675,210.7560) -- (43.9675,179.8110);
  \path[draw=black,line join=round,line cap=round,line width=0.800pt]
    (43.9675,179.8110) -- (80.8225,179.8110);
  \path[draw=black,line join=round,line cap=round,line width=0.800pt]
    (43.9675,239.3730) -- (80.8225,239.3730);
  \path[draw=black,line join=round,line cap=round,line width=0.800pt]
    (42.1247,210.7560) .. controls (42.1247,209.7380) and (42.9497,208.9130) ..
    (43.9675,208.9130) .. controls (44.9852,208.9130) and (45.8102,209.7380) ..
    (45.8102,210.7560) .. controls (45.8102,211.7730) and (44.9852,212.5990) ..
    (43.9675,212.5990) .. controls (42.9497,212.5990) and (42.1247,211.7730) ..
    (42.1247,210.7560) -- cycle;
  \path[draw=black,line join=round,line cap=round,line width=0.800pt]
    (78.9797,210.7560) .. controls (78.9797,209.7380) and (79.8048,208.9130) ..
    (80.8225,208.9130) .. controls (81.8402,208.9130) and (82.6652,209.7380) ..
    (82.6652,210.7560) .. controls (82.6652,211.7730) and (81.8402,212.5990) ..
    (80.8225,212.5990) .. controls (79.8048,212.5990) and (78.9797,211.7730) ..
    (78.9797,210.7560) -- cycle;
      \path[draw=c2f3b73,line join=round,line cap=round,line width=0.800pt]
        (75.4793,210.7560) .. controls (72.3144,210.7560) and (46.0240,210.7560) ..
        (46.0240,210.7560);
      \path[fill=c2f3b73] (74.6043,210.7560) -- (71.9792,214.2560) --
        (78.9793,210.7560) -- (71.9793,207.2560) -- (74.6043,210.7560) -- cycle;
  \path[draw=c74b587,dash pattern=on 1.60pt,line join=round,line cap=round,line
    width=0.800pt] (124.8500,239.3730) -- (124.8500,210.7560) --
    (124.8500,179.8110);
  \path[draw=cf56356,dash pattern=on 1.60pt,line join=round,line cap=round,line
    width=0.800pt] (87.9950,239.3730) -- (87.9950,210.7560) -- (87.9950,179.8110);
  \path[draw=black,line join=round,line cap=round,line width=0.800pt]
    (87.9950,179.8110) -- (124.8500,179.8110);
  \path[draw=black,line join=round,line cap=round,line width=0.800pt]
    (87.9950,239.3730) -- (124.8500,239.3730);
  \path[draw=black,line join=round,line cap=round,line width=0.800pt]
    (86.1523,210.7560) .. controls (86.1523,209.7380) and (86.9773,208.9130) ..
    (87.9950,208.9130) .. controls (89.0127,208.9130) and (89.8378,209.7380) ..
    (89.8378,210.7560) .. controls (89.8378,211.7730) and (89.0127,212.5990) ..
    (87.9950,212.5990) .. controls (86.9773,212.5990) and (86.1523,211.7730) ..
    (86.1523,210.7560) -- cycle;
  \path[draw=black,line join=round,line cap=round,line width=0.800pt]
    (123.0070,210.7560) .. controls (123.0070,209.7380) and (123.8320,208.9130) ..
    (124.8500,208.9130) .. controls (125.8680,208.9130) and (126.6930,209.7380) ..
    (126.6930,210.7560) .. controls (126.6930,211.7730) and (125.8680,212.5990) ..
    (124.8500,212.5990) .. controls (123.8320,212.5990) and (123.0070,211.7730) ..
    (123.0070,210.7560) -- cycle;
      \path[draw=c2f3b73,line join=round,line cap=round,line width=0.800pt]
        (93.5520,210.7560) .. controls (96.7168,210.7560) and (123.0070,210.7560) ..
        (123.0070,210.7560);
      \path[fill=c2f3b73] (94.4270,210.7560) -- (97.0520,207.2560) --
        (90.0520,210.7560) -- (97.0520,214.2560) -- (94.4270,210.7560) -- cycle;
  \path[draw=black,line join=round,line width=0.800pt] (244.6330,180.9750) --
    (244.6330,240.5370);
  \path[draw=black,line join=round,line width=0.800pt] (191.0090,180.9750) --
    (191.0090,240.5370);
  \path[draw=black,line join=round,line cap=round,line width=0.800pt]
    (230.9260,207.5620) .. controls (230.9260,206.9530) and (231.4200,206.4590) ..
    (232.0290,206.4590) .. controls (232.6380,206.4590) and (233.1320,206.9530) ..
    (233.1320,207.5620) .. controls (233.1320,208.1710) and (232.6380,208.6650) ..
    (232.0290,208.6650) .. controls (231.4200,208.6650) and (230.9260,208.1710) ..
    (230.9260,207.5620) -- cycle;
  \path[draw=black,line join=round,line cap=round,line width=0.800pt]
    (202.2530,215.9390) .. controls (202.2530,215.3300) and (202.7470,214.8370) ..
    (203.3560,214.8370) .. controls (203.9650,214.8370) and (204.4590,215.3300) ..
    (204.4590,215.9390) .. controls (204.4590,216.5490) and (203.9650,217.0420) ..
    (203.3560,217.0420) .. controls (202.7470,217.0420) and (202.2530,216.5490) ..
    (202.2530,215.9390) -- cycle;
  \path[draw=black,line join=round,line width=0.800pt] (191.0090,180.9750) ..
    controls (191.0090,178.1460) and (203.0130,175.8520) .. (217.8210,175.8520) ..
    controls (232.6290,175.8520) and (244.6330,178.1460) .. (244.6330,180.9750) ..
    controls (244.6330,183.8040) and (232.6290,186.0970) .. (217.8210,186.0970) ..
    controls (203.0130,186.0970) and (191.0090,183.8040) .. (191.0090,180.9750) --
    cycle;
  \path[draw=black,dash pattern=on 1.60pt,line join=round,line cap=round,line
    width=0.800pt] (191.0090,240.5370) .. controls (191.0090,237.7080) and
    (203.0130,235.4140) .. (217.8210,235.4140) .. controls (232.6290,235.4140) and
    (244.6330,237.7080) .. (244.6330,240.5370);
  \path[draw=black,line join=round,line width=0.800pt] (244.6330,240.5370) ..
    controls (244.6330,243.3660) and (232.6290,245.6600) .. (217.8210,245.6600) ..
    controls (203.0130,245.6600) and (191.0090,243.3660) .. (191.0090,240.5370);
      \path[draw=c2f3b73,dash pattern=on 1.60pt,line join=round,line cap=round,line
        width=0.800pt] (227.3540,207.1300) .. controls (224.3920,206.9150) and
        (221.1790,206.7970) .. (217.8210,206.7970) .. controls (203.0130,206.7970) and
        (191.0090,209.0900) .. (191.0090,211.9200);
      \path[fill=c2f3b73] (226.4820,207.0530) -- (223.5570,210.3060) --
        (230.8400,207.4400) -- (224.1770,203.3340) -- (226.4820,207.0530) -- cycle;
  \path[draw=c2f3b73,line join=round,line cap=round,line width=0.800pt]
    (202.0790,216.0670) .. controls (195.3710,215.1360) and (191.0090,213.6250) ..
    (191.0090,211.9200);
  \path[draw=c2f3b73,line join=round,line cap=round,line width=0.800pt]
    (244.6330,211.9200) .. controls (244.6330,214.7490) and (232.6290,217.0420) ..
    (217.8210,217.0420) .. controls (212.8840,217.0420) and (208.2590,216.7870) ..
    (204.2860,216.3430);
      \path[draw=c2f3b73,dash pattern=on 1.60pt,line join=round,line cap=round,line
        width=0.800pt] (236.7340,208.2890) .. controls (241.6130,209.2160) and
        (244.6330,210.5000) .. (244.6330,211.9200);
      \path[fill=c2f3b73] (237.5980,208.4270) -- (240.7450,205.3880) --
        (233.2780,207.7330) -- (239.6350,212.2990) -- (237.5980,208.4270) -- cycle;
  \path[draw=c74b587,line join=round,draw opacity=0.798,line width=1.600pt]
    (203.3560,217.0420) -- (203.3860,244.8530);
  \path[draw=c74b587,line join=round,draw opacity=0.798,line width=1.600pt]
    (203.3340,185.2880) -- (203.3560,214.8370);
  \path[draw=cf56356,dash pattern=on 1.60pt,line join=round,draw
    opacity=0.799,line width=1.600pt] (232.0290,206.4590) -- (232.0970,176.6390);
  \path[draw=cf56356,dash pattern=on 1.60pt,line join=round,draw
    opacity=0.799,line width=1.600pt] (232.0290,208.6650) -- (232.0300,236.1930);
  \path[fill=black,line width=0.600pt] (142.5000,211.8900) node[above right]
    (text6862) {$\implies$};
  \path[fill=black,line width=0.600pt] (57.0000,204.3900) node[above right]
    (text6866) {$a_1$};
  \path[fill=black,line width=0.600pt] (101.8469,204.3900) node[above right]
    (text6870) {$b_1$};
  \path[fill=black,line width=0.600pt] (30.0000,200.6400) node[above right]
    (text6874) {$\overline{a_0}$};
  \path[fill=black,line width=0.600pt] (30.0000,230.6400) node[above right]
    (text6878) {$\underline{a_0}$};
  \path[fill=black,line width=0.600pt] (127.7816,200.6400) node[above right]
    (text6882) {$\underline{b_0}$};
  \path[fill=black,line width=0.600pt] (127.4326,230.6400) node[above right]
    (text6886) {$\overline{b_0}$};
  \path[fill=black,line width=0.600pt] (168.0857,206.8094) node[above right]
    (text6892) {$b_1$};
  \path[draw=black,line join=miter,line cap=butt,miter limit=4.00,line
    width=0.300pt] (236.2500,215.6400) -- (251.2500,223.1400);
  \path[fill=black,line width=0.600pt] (253.0366,226.8900) node[above right]
    (text6898) {$a_1$};
  \path[draw=black,line join=miter,line cap=butt,miter limit=4.00,line
    width=0.300pt] (180.0000,204.3900) -- (198.4544,215.4429);

  \end{tikzpicture}

  \caption{Cut-and-paste of two sheets of disks.}
  \label{fig:cut_paste_2}
\end{figure}

And paste upper paths and lower paths as follows: $\overline{a_{2}}\sim \underline{b_{2}}, \overline{b_{2}}\sim \underline{a_{2}}, \overline{a_{0}}\sim \underline{b_{0}},\overline{b_0} \sim\underline{a_{0}}$.
 As a result we get an annuls as a 2-fold covering.  The identity maps of $a$ and $b$ to $S_{0,1,(2)}$ give us the branched covering map and the only nontrivial covering transformation is interchanging $a$ and $b$, which means interchanging two boundary components.

For 3-fold branched covering space, we start with three sheets of disks, called $a,b$ and $c$, and similarly labelled paths.
First, cut the disks along the paths $a_0, b_0, c_0, a_2, b_2, c_2$.

\begin{figure}[H]
  \centering

  \definecolor{cf56356}{RGB}{245,99,86}
  \definecolor{c74b587}{RGB}{116,181,135}
  \definecolor{c2f3b73}{RGB}{47,59,115}

  \begin{tikzpicture}[y=0.80pt, x=0.80pt, yscale=-1.000000, xscale=1.000000, inner sep=0pt, outer sep=0pt]
  \path[draw=cf56356,dash pattern=on 1.60pt,line join=round,line cap=round,draw
    opacity=0.799,line width=0.800pt] (79.6756,498.3360) -- (70.5803,513.0060) --
    (54.0984,513.0910);
  \path[draw=c74b587,dash pattern=on 1.60pt,line join=round,line cap=round,line
    width=0.800pt] (121.9438,531.4070) -- (129.3148,516.6650) --
    (121.9438,501.9230);
  \path[draw=black,dash pattern=on 0.80pt,line join=round,line cap=round,line
    width=0.800pt] (131.1578,516.6650) .. controls (131.1578,515.6470) and
    (130.3328,514.8220) .. (129.3148,514.8220) .. controls (128.2968,514.8220) and
    (127.4718,515.6470) .. (127.4718,516.6650) .. controls (127.4718,517.6830) and
    (128.2968,518.5080) .. (129.3148,518.5080) .. controls (130.3328,518.5080) and
    (131.1578,517.6830) .. (131.1578,516.6650) -- cycle;
  \begin{scope}[shift={(-6.8082,159.83099)}]
      \path[draw=c2f3b73,line join=round,line cap=round,line width=0.800pt]
        (89.8713,356.8340) .. controls (94.2957,356.8340) and (134.2800,356.8340) ..
        (134.2800,356.8340);
      \path[fill=c2f3b73] (90.7463,356.8340) -- (93.3712,353.3340) --
        (86.3713,356.8340) -- (93.3713,360.3340) -- cycle;
  \end{scope}
  \path[draw=c74b587,dash pattern=on 1.60pt,line join=round,line cap=round,line
    width=0.800pt] (63.7500,466.8900) -- (47.2682,466.9750) -- (37.4470,480.2120);
  \path[draw=black,line join=round,line cap=round,line width=0.800pt]
    (53.8374,520.5660) -- (37.3220,553.5140);
  \path[draw=black,line join=round,line cap=round,line width=0.800pt]
    (37.4470,480.2120) -- (54.0984,513.0910);
  \path[draw=black,line join=round,line cap=round,line width=0.800pt]
    (86.8828,502.1320) -- (121.9438,501.9230);
  \path[draw=black,line join=round,line cap=round,line width=0.800pt]
    (63.7500,466.8900) -- (79.7449,498.4730);
  \path[draw=black,dash pattern=on 0.80pt,line join=round,line cap=round,line
    width=0.800pt] (46.4356,465.3310) .. controls (45.5277,465.7910) and
    (45.1644,466.9000) .. (45.6243,467.8080) .. controls (46.0841,468.7160) and
    (47.1929,469.0790) .. (48.1008,468.6190) .. controls (49.0087,468.1600) and
    (49.3720,467.0510) .. (48.9122,466.1430) .. controls (48.4523,465.2350) and
    (47.3435,464.8720) .. (46.4356,465.3310) -- cycle;
  \path[draw=black,dash pattern=on 0.80pt,line join=round,line cap=round,line
    width=0.800pt] (69.4867,518.8370) .. controls (68.5787,519.2970) and
    (68.2155,520.4060) .. (68.6753,521.3140) .. controls (69.1351,522.2220) and
    (70.2439,522.5850) .. (71.1518,522.1250) .. controls (72.0597,521.6650) and
    (72.4230,520.5570) .. (71.9632,519.6490) .. controls (71.5034,518.7410) and
    (70.3946,518.3780) .. (69.4867,518.8370) -- cycle;
  \begin{scope}[shift={(-6.8082,159.83099)}]
      \path[draw=c2f3b73,line join=round,line cap=round,line width=0.800pt]
        (74.9741,348.4070) .. controls (73.4310,345.3600) and (60.4190,319.6680) ..
        (60.4190,319.6680) -- (54.9090,308.7880);
      \path[fill=c2f3b73] (74.5788,347.6270) -- (70.2704,346.8660) --
        (76.5554,351.5300) -- (76.5152,343.7040) -- cycle;
  \end{scope}
  \path[draw=c74b587,dash pattern=on 1.60pt,line join=round,line cap=round,line
    width=0.800pt] (37.1593,553.7890) -- (47.0352,566.9840) -- (63.5172,567.0010);
  \path[draw=black,line join=round,line cap=round,line width=0.800pt]
    (86.8134,531.4070) -- (121.9438,531.4070);
  \path[draw=black,line join=round,line cap=round,line width=0.800pt]
    (63.5172,567.0010) -- (79.4145,535.2230);
  \path[draw=black,dash pattern=on 0.80pt,line join=round,line cap=round,line
    width=0.800pt] (46.2094,568.6320) .. controls (47.1192,569.0880) and
    (48.2265,568.7200) .. (48.6826,567.8100) .. controls (49.1386,566.9000) and
    (48.7708,565.7930) .. (47.8610,565.3370) .. controls (46.9512,564.8810) and
    (45.8439,565.2490) .. (45.3878,566.1590) .. controls (44.9317,567.0680) and
    (45.2996,568.1760) .. (46.2094,568.6320) -- cycle;
  \begin{scope}[shift={(-6.8082,159.83099)}]
      \path[draw=c2f3b73,line join=round,line cap=round,line width=0.800pt]
        (74.5696,365.8060) .. controls (72.5869,369.7610) and (54.6692,405.5060) ..
        (54.6692,405.5060);
      \path[fill=c2f3b73] (74.1774,366.5880) -- (76.1300,370.5030) --
        (76.1380,362.6770) -- (69.8722,367.3660) -- cycle;
  \end{scope}
  \path[draw=cf56356,dash pattern=on 1.60pt,line join=round,line cap=round,draw
    opacity=0.799,line width=0.800pt] (86.8828,501.7160) --
    (77.7875,516.3860)(86.8828,531.4070) -- (77.7875,516.6650)(70.3754,520.3700)
    -- (53.8935,520.4550);
  \path[draw=black,dash pattern=on 0.80pt,line join=round,line cap=round,line
    width=0.800pt] (76.8848,515.0480) .. controls (75.9768,515.5080) and
    (75.6136,516.6160) .. (76.0734,517.5240) .. controls (76.5332,518.4320) and
    (77.6420,518.7950) .. (78.5499,518.3360) .. controls (79.4578,517.8760) and
    (79.8211,516.7670) .. (79.3613,515.8590) .. controls (78.9015,514.9510) and
    (77.7927,514.5880) .. (76.8848,515.0480) -- cycle;
  \path[draw=black,dash pattern=on 0.80pt,line join=round,line cap=round,line
    width=0.800pt] (69.0718,511.3350) .. controls (68.1638,511.7950) and
    (67.8006,512.9040) .. (68.2604,513.8120) .. controls (68.7202,514.7200) and
    (69.8290,515.0830) .. (70.7369,514.6230) .. controls (71.6448,514.1630) and
    (72.0081,513.0550) .. (71.5483,512.1470) .. controls (71.0885,511.2390) and
    (69.9797,510.8750) .. (69.0718,511.3350) -- cycle;
  \path[draw=cf56356,dash pattern=on 1.60pt,line join=round,line cap=round,draw
    opacity=0.799,line width=0.800pt] (79.4145,535.1970) -- (70.3192,520.4550);
  \path[draw=black,line join=round,line cap=round,line width=0.800pt]
    (211.9106,480.5380);
  \path[draw=black,line join=round,line cap=round,line width=0.800pt]
    (360.4563,517.6260) -- (361.3673,496.2900);
  \path[draw=black,line join=round,line cap=round,line width=0.800pt]
    (361.3553,496.6010) .. controls (361.9693,487.9880) and (338.9403,480.8460) ..
    (338.9403,480.8460);
  \path[draw=black,line join=round,line cap=round,line width=0.800pt]
    (328.9513,492.3030) .. controls (328.0433,492.7630) and (327.6793,493.8710) ..
    (328.1393,494.7790) .. controls (328.5993,495.6870) and (329.7083,496.0500) ..
    (330.6163,495.5910) .. controls (331.5243,495.1310) and (331.8873,494.0220) ..
    (331.4273,493.1140) .. controls (330.9673,492.2060) and (329.8583,491.8430) ..
    (328.9513,492.3030) -- cycle;
  \path[draw=black,dash pattern=on 0.80pt,line join=round,line cap=round,line
    width=0.800pt] (221.5209,492.3030) .. controls (220.6130,492.7630) and
    (220.2497,493.8710) .. (220.7095,494.7790) .. controls (221.1693,495.6870) and
    (222.2781,496.0500) .. (223.1860,495.5910) .. controls (224.0939,495.1310) and
    (224.4572,494.0220) .. (223.9974,493.1140) .. controls (223.5376,492.2060) and
    (222.4288,491.8430) .. (221.5209,492.3030) -- cycle;
  \begin{scope}[shift={(170.14629,-229.96301)}]
      \path[draw=c2f3b73,line join=round,line cap=round,line width=0.800pt]
        (164.9540,724.1900) .. controls (172.7880,725.2270) and (191.1310,728.3690) ..
        (190.8460,735.0300);
      \path[fill=c2f3b73] (165.8230,724.2920) -- (168.8370,721.1220) --
        (161.4770,723.7830) -- (168.0220,728.0740) -- cycle;
  \end{scope}
  \path[draw=c2f3b73,dash pattern=on 1.60pt,line join=round,line cap=round,line
    width=0.800pt] (192.1001,495.5680) .. controls (196.5625,490.2740) and
    (220.5316,493.6690) .. (220.5316,493.6690);
  \path[draw=c2f3b73,dash pattern=on 1.60pt,line join=round,line cap=round,line
    width=0.800pt] (263.6965,449.3330) .. controls (236.3008,449.3330) and
    (222.6732,478.4390) .. (222.6732,492.1320);
  \path[draw=c2f3b73,dash pattern=on 1.60pt,line join=round,line cap=round,line
    width=0.800pt] (261.3948,534.5050) .. controls (240.6180,534.5860) and
    (222.3535,509.4220) .. (222.3909,495.7890);
  \begin{scope}[shift={(170.14629,-229.96301)}]
      \path[draw=c2f3b73,line join=round,line cap=round,line width=0.800pt]
        (159.2010,729.2320) .. controls (156.7260,740.9160) and (142.7120,760.3680) ..
        (119.3740,764.4270);
      \path[fill=c2f3b73] (159.1020,730.1010) -- (162.2830,733.1050) --
        (159.5960,725.7540) -- (155.3270,732.3140) -- cycle;
  \end{scope}
  \begin{scope}[shift={(170.14629,-229.96301)}]
      \path[draw=c2f3b73,line join=round,line cap=round,line width=0.800pt]
        (159.3560,718.5810) .. controls (156.7210,704.0110) and (139.1190,679.3220) ..
        (118.0620,679.6260);
      \path[fill=c2f3b73] (159.2710,717.7100) -- (155.5320,715.4380) --
        (159.6970,722.0640) -- (162.4990,714.7570) -- cycle;
  \end{scope}
  \path[draw=cf56356,line join=round,draw opacity=0.799,line width=1.600pt]
    (327.9473,493.7910) -- (317.8453,493.7750);
  \path[draw=cf56356,line join=round,draw opacity=0.799,line width=1.600pt]
    (331.0293,495.3040) -- (338.9403,502.6600);
  \path[draw=cf56356,line join=round,draw opacity=0.799,line width=1.600pt]
    (330.9253,492.5000) -- (338.7983,481.0580);
  \path[draw=c74b587,dash pattern=on 0.80pt,line join=round,draw
    opacity=0.799,line width=1.600pt] (221.0744,495.2730) -- (211.7129,503.0190);
  \path[draw=c74b587,dash pattern=on 0.80pt,line join=round,draw
    opacity=0.799,line width=1.600pt] (221.5209,492.3030) -- (211.9106,480.5380);
  \path[draw=c74b587,dash pattern=on 0.80pt,line join=round,draw
    opacity=0.799,line width=1.600pt] (224.1966,493.9480) -- (232.7098,493.9470);
  \path[draw=black,line join=round,line cap=round,line width=0.800pt]
    (339.0473,502.3860) .. controls (333.4443,522.1760) and (312.0343,536.9080) ..
    (286.4783,536.9080) .. controls (284.7523,536.9080) and (283.0443,536.8410) ..
    (281.3593,536.7100);
  \path[draw=black,line join=round,line cap=round,line width=0.800pt]
    (270.0838,534.8060) .. controls (248.2767,529.0460) and (232.4797,512.1020) ..
    (232.4797,492.0970) .. controls (232.4797,472.1470) and (248.1885,455.2420) ..
    (269.9010,449.4360);
  \path[draw=black,line join=round,line cap=round,line width=0.800pt]
    (338.8843,481.2470);
  \path[draw=black,line join=round,line cap=round,line width=0.800pt]
    (212.2106,480.7280) .. controls (218.4453,461.6900) and (239.4354,447.6880) ..
    (264.3694,447.6880) .. controls (294.2163,447.6880) and (318.4123,467.7510) ..
    (318.4123,492.5000) .. controls (318.4123,517.2490) and (294.2163,537.3120) ..
    (264.3694,537.3120) .. controls (238.8681,537.3120) and (217.4922,522.6660) ..
    (211.8081,502.9650);
  \path[draw=black,line join=round,line cap=round,line width=0.800pt]
    (338.8843,481.2470) .. controls (333.0393,461.7400) and (311.7933,447.2850) ..
    (286.4783,447.2850) .. controls (284.5853,447.2850) and (282.7153,447.3660) ..
    (280.8733,447.5230);
  \path[draw=black,line join=round,line cap=round,line width=0.800pt]
    (211.8081,502.9650);
  \path[draw=black,line join=round,line cap=round,line width=0.800pt]
    (189.0724,489.9000) .. controls (189.0082,494.6820) and (188.8156,509.0270) ..
    (188.8156,509.0270);
  \path[draw=black,line join=round,line cap=round,line width=0.800pt]
    (195.0000,504.3900) .. controls (199.8500,503.1790) and (206.6811,502.8860) ..
    (212.2106,502.7780);
  \path[draw=black,line join=round,line cap=round,line width=0.800pt]
    (273.9513,579.0210) .. controls (281.1513,580.1400) and (289.5103,580.3850) ..
    (300.3723,580.3320) .. controls (313.6983,580.2680) and (359.7903,533.2370) ..
    (360.4563,517.6260) .. controls (360.9923,505.0670) and (338.8843,502.2590) ..
    (338.8843,502.2590);
  \path[draw=black,line join=round,line cap=round,line width=0.800pt]
    (189.0724,489.9000) .. controls (189.0724,502.6690) and (210.4271,537.8940) ..
    (238.7616,561.3460) .. controls (248.0684,569.0490) and (255.3808,573.7020) ..
    (263.4262,576.4720);
  \path[draw=black,line join=round,line cap=round,line width=0.800pt]
    (211.8081,480.5410) .. controls (199.6033,480.4870) and (189.0724,481.1200) ..
    (189.0724,489.9000);
  \path[draw=black,line join=round,line cap=round,line width=0.800pt]
    (265.8978,578.6820) .. controls (258.6594,580.1340) and (250.3179,580.4780) ..
    (238.7616,580.4780) .. controls (225.0643,580.4780) and (188.9440,526.4350) ..
    (188.8156,509.0270);
  \path[draw=black,line join=round,line cap=round,line width=0.800pt]
    (353.1873,510.0740) .. controls (353.1873,510.0740) and (328.5463,538.9780) ..
    (300.3723,560.4130) .. controls (289.3323,568.8120) and (281.6343,573.7540) ..
    (273.4323,576.6310);
  \path[draw=c2f3b73,line join=round,line cap=round,line width=0.800pt]
    (188.8976,502.9200) .. controls (189.5575,514.0650) and (226.2493,577.0940) ..
    (259.9478,577.0900);
  \path[draw=c2f3b73,dash pattern=on 1.60pt,line join=round,line cap=round,line
    width=0.800pt] (356.2623,512.9290) .. controls (353.1873,527.8350) and
    (308.2823,576.4010) .. (277.8233,577.0900);
  \path[draw=black,line join=round,line cap=round,line width=0.800pt]
    (265.8978,578.6820) -- (273.4323,576.6310);
  \path[draw=black,line join=miter,line cap=butt,miter limit=4.00,line
    width=0.300pt] (322.8233,472.0900) -- (342.6463,461.9270);
  \path[fill=black,line width=0.600pt] (346.3963,465.6770) node[above right]
    (text17612) {$a_1$};
  \path[draw=black,line join=miter,line cap=butt,miter limit=4.00,line
    width=0.300pt] (319.0733,517.0900) -- (301.3963,506.9270);
  \path[fill=black,line width=0.600pt] (290.1463,506.9270) node[above right]
    (text17620) {$b_1$};
  \path[draw=black,line join=miter,line cap=butt,miter limit=4.00,line
    width=0.300pt] (352.8233,498.3400) -- (365.1463,540.6770) --
    (365.1463,540.6770);
  \path[fill=black,line width=0.600pt] (361.5789,553.0821) node[above right]
    (text17626) {$c_1$};
  \path[fill=black,line width=0.600pt] (146.2500,518.2649) node[above right]
    (text16721) {$\implies$};
  \path[fill=black,line width=0.600pt] (61.4418,491.7210) node[above right]
    (text29460) {$a_1$};
  \path[fill=black,line width=0.600pt] (45.6918,544.2210) node[above right]
    (text29468) {$b_1$};
  \path[fill=black,line width=0.600pt] (97.5000,526.8900) node[above right]
    (text16725) {$c_1$};

  \end{tikzpicture}

  \caption{3-fold branched covering of disk with two branch points.}
  \label{fig:S_1_1_2}
\end{figure}
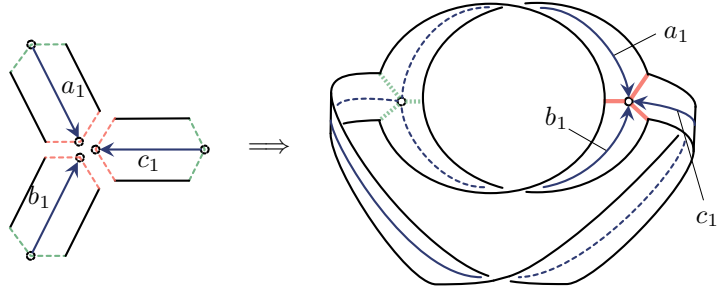

Then paste upper paths and lower paths which are labelled by ``$2$'' as follows:
$\overline{a_2}\sim\underline{c_2}, \overline{c_2} \sim \underline{b_2}, \overline{b_2} \sim \underline{a_2}$.

Now we have two choices of pasting the paths corresponding to the path ``$0$'', one gives us a pair-of-pants ($S_{0,3,(2)}$) and the other gives us a torus with one boundary component. In the case of pair-of-pants, we do not get a desired branched covering modeled by $y^3=(x-p_1)(x-p_2)$, because a small (counterclockwise oriented) loop around the point $p_1$ travels the surface in the order of $a\rightarrow c \rightarrow b$.

The correct one is obtained by pasting upper paths and lower paths which are labelled by ``$0$'' as follows:
$\overline{a_0}\sim\underline{b_0}, \overline{b_0}\sim\underline{c_0}, \overline{c_0}\sim\underline{a_0}$.
Then the identity maps on each of $a,b$ and $c$ give us the 3-fold branched covering map.


\subsection{More branch points}

For three branch points case, similar to two branch points case, first cut the disk $a$ along the paths joining the boundary and the branch points. Then paste them in a proper order.

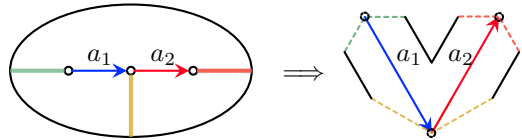
\begin{figure}[H]
  \centering

  \definecolor{c002fff}{RGB}{0,47,255}
  \definecolor{cff001f}{RGB}{255,0,31}
  \definecolor{c74b587}{RGB}{116,181,135}
  \definecolor{cf56356}{RGB}{245,99,86}
  \definecolor{cd9ac32}{RGB}{217,172,50}

  \begin{tikzpicture}[y=0.80pt, x=0.80pt, yscale=-1.000000, xscale=1.000000, inner sep=0pt, outer sep=0pt]
  \path[draw=black,line join=round,line cap=round,line width=0.800pt]
    (53.5433,94.6406) .. controls (53.5433,77.3393) and (79.1191,63.3139) ..
    (110.6690,63.3139) .. controls (142.2180,63.3139) and (167.7940,77.3393) ..
    (167.7940,94.6406) .. controls (167.7940,111.9420) and (142.2180,125.9670) ..
    (110.6690,125.9670) .. controls (79.1191,125.9670) and (53.5433,111.9420) ..
    (53.5433,94.6406) -- cycle;
  \path[draw=black,line join=round,line cap=round,line width=0.800pt]
    (108.8260,94.6406) .. controls (108.8260,93.6229) and (109.6510,92.7979) ..
    (110.6690,92.7979) .. controls (111.6860,92.7979) and (112.5110,93.6229) ..
    (112.5110,94.6406) .. controls (112.5110,95.6583) and (111.6860,96.4834) ..
    (110.6690,96.4834) .. controls (109.6510,96.4834) and (108.8260,95.6583) ..
    (108.8260,94.6406) -- cycle;
  \path[draw=black,line join=round,line cap=round,line width=0.800pt]
    (79.3418,94.6406) .. controls (79.3418,93.6229) and (80.1668,92.7979) ..
    (81.1845,92.7979) .. controls (82.2022,92.7979) and (83.0273,93.6229) ..
    (83.0273,94.6406) .. controls (83.0273,95.6583) and (82.2022,96.4834) ..
    (81.1845,96.4834) .. controls (80.1668,96.4834) and (79.3418,95.6583) ..
    (79.3418,94.6406) -- cycle;
  \path[draw=black,line join=round,line cap=round,line width=0.800pt]
    (138.3100,94.6406) .. controls (138.3100,93.6229) and (139.1350,92.7979) ..
    (140.1520,92.7979) .. controls (141.1700,92.7979) and (141.9950,93.6229) ..
    (141.9950,94.6406) .. controls (141.9950,95.6583) and (141.1700,96.4834) ..
    (140.1520,96.4834) .. controls (139.1350,96.4834) and (138.3100,95.6583) ..
    (138.3100,94.6406) -- cycle;
      \path[draw=c002fff,line join=round,line cap=round,line width=0.800pt]
        (105.3240,94.6406) .. controls (102.8600,94.6406) and (83.0273,94.6406) ..
        (83.0273,94.6406);
      \path[fill=c002fff] (104.4490,94.6406) -- (101.8240,98.1406) --
        (108.8240,94.6406) -- (101.8240,91.1406) -- (104.4490,94.6406) -- cycle;
      \path[draw=cff001f,line join=round,line cap=round,line width=0.800pt]
        (134.8080,94.6406) .. controls (132.3440,94.6406) and (112.5110,94.6406) ..
        (112.5110,94.6406);
      \path[fill=cff001f] (133.9330,94.6406) -- (131.3080,98.1406) --
        (138.3080,94.6406) -- (131.3080,91.1406) -- (133.9330,94.6406) -- cycle;
  \path[draw=c74b587,line join=round,draw opacity=0.799,line width=1.560pt]
    (53.5433,94.6406) -- (79.3418,94.6406);
  \path[draw=cf56356,line join=round,line width=1.560pt] (141.9950,94.6406) --
    (167.7940,94.6406);
  \path[draw=cd9ac32,line join=round,draw opacity=0.803,line width=1.560pt]
    (110.6690,96.4834) -- (110.6690,125.9670);
  \path[draw=cd9ac32,line join=round,line cap=round,line width=0.800pt]
    (262.5100,107.5020) -- (252.9350,124.0870);
      \path[draw=c002fff,line join=round,line cap=round,line width=0.800pt]
        (251.1840,121.0520) .. controls (248.4150,116.2550) and (221.2490,69.1861) ..
        (221.2490,69.1861);
      \path[fill=c002fff] (250.7460,120.2940) -- (246.4030,119.7700) --
        (252.9330,124.0830) -- (252.4650,116.2710) -- (250.7460,120.2940) -- cycle;
  \path[draw=c74b587,dash pattern=on 1.60pt,line join=round,line cap=round,line
    width=0.800pt] (240.3990,69.1861) -- (221.2490,69.1861) -- (211.6730,85.7708);
      \path[draw=cff001f,line join=round,line cap=round,line width=0.800pt]
        (282.9300,72.0117) .. controls (280.1590,76.8223) and (252.9350,124.0870) ..
        (252.9350,124.0870);
      \path[fill=cff001f] (282.4940,72.7699) -- (284.2160,76.7915) --
        (284.6770,68.9789) -- (278.1500,73.2976) -- (282.4940,72.7699) -- cycle;
  \path[draw=cf56356,dash pattern=on 1.60pt,line join=round,line cap=round,line
    width=0.800pt] (294.2650,85.5551) -- (284.6780,68.9775) -- (265.5280,68.9916);
  \path[draw=black,line join=round,line cap=round,line width=0.800pt]
    (240.3680,69.2038) -- (252.8780,90.8511) -- (265.2650,69.3391);
  \path[draw=cd9ac32,dash pattern=on 1.60pt,line join=round,line cap=round,draw
    opacity=0.798,line width=0.800pt] (224.3590,107.7230) -- (252.9350,124.0870);
  \path[draw=cd9ac32,dash pattern=on 1.60pt,line join=round,line cap=round,draw
    opacity=0.798,line width=0.800pt] (281.7130,107.5120) -- (252.9350,124.0870);
  \path[draw=black,dash pattern=on 0.80pt,line join=round,line cap=round,line
    width=0.800pt] (254.5780,124.1250) .. controls (254.5780,123.1070) and
    (253.7530,122.2820) .. (252.7360,122.2820) .. controls (251.7180,122.2820) and
    (250.8930,123.1070) .. (250.8930,124.1250) .. controls (250.8930,125.1420) and
    (251.7180,125.9670) .. (252.7360,125.9670) .. controls (253.7530,125.9670) and
    (254.5780,125.1420) .. (254.5780,124.1250) -- cycle;
  \path[draw=black,dash pattern=on 0.80pt,line join=round,line cap=round,line
    width=0.800pt] (286.6050,69.2038) .. controls (286.6050,68.1861) and
    (285.7800,67.3610) .. (284.7620,67.3610) .. controls (283.7440,67.3610) and
    (282.9190,68.1861) .. (282.9190,69.2038) .. controls (282.9190,70.2215) and
    (283.7440,71.0465) .. (284.7620,71.0465) .. controls (285.7800,71.0465) and
    (286.6050,70.2215) .. (286.6050,69.2038) -- cycle;
  \path[draw=black,dash pattern=on 0.80pt,line join=round,line cap=round,line
    width=0.800pt] (223.0840,68.9774) .. controls (223.0840,67.9597) and
    (222.2590,67.1347) .. (221.2410,67.1347) .. controls (220.2230,67.1347) and
    (219.3980,67.9597) .. (219.3980,68.9774) .. controls (219.3980,69.9952) and
    (220.2230,70.8202) .. (221.2410,70.8202) .. controls (222.2590,70.8202) and
    (223.0840,69.9952) .. (223.0840,68.9774) -- cycle;
  \path[draw=black,line join=round,line cap=round,line width=0.800pt]
    (281.7130,107.5120) -- (294.1410,85.9296);
  \path[draw=black,line join=round,line cap=round,line width=0.800pt]
    (211.6730,85.7708) -- (224.3590,107.7230);
  \path[fill=black,line width=0.600pt] (178.8153,98.7931) node[above right]
    (text14362) {$\implies$};
  \path[fill=black,line width=0.600pt] (90.0000,91.8900) node[above right]
    (text14370) {$a_1$};
  \path[fill=black,line width=0.600pt] (120.0000,91.8900) node[above right]
    (text14374) {$a_2$};
  \path[fill=black,line width=0.600pt] (236.2500,91.8900) node[above right]
    (text14382) {$a_1$};
  \path[fill=black,line width=0.600pt] (260.6081,91.8900) node[above right]
    (text14386) {$a_2$};

  \end{tikzpicture}

  \caption{Cut a disk with three marked points.}
  \label{fig:S_0_1_3}
\end{figure}

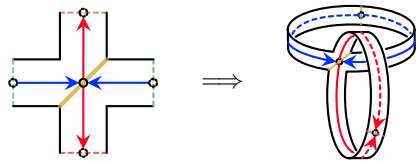
\begin{figure}[H]
  \centering

  \definecolor{cf56356}{RGB}{245,99,86}
  \definecolor{c74b587}{RGB}{116,181,135}
  \definecolor{cd9ac32}{RGB}{217,172,50}
  \definecolor{c002fff}{RGB}{0,47,255}
  \definecolor{cff001f}{RGB}{255,0,31}

  \begin{tikzpicture}[y=0.80pt, x=0.80pt, yscale=-1.000000, xscale=1.000000, inner sep=0pt, outer sep=0pt]
  \path[draw=cf56356,dash pattern=on 1.60pt,line join=round,line cap=round,draw
    opacity=0.801,line width=0.800pt] (100.1940,158.1450) -- (122.3070,158.1450);
  \path[draw=cf56356,dash pattern=on 1.60pt,line join=round,line cap=round,line
    width=0.800pt] (122.3070,224.4840) -- (100.1940,224.4840);
  \path[draw=black,line join=round,line cap=round,line width=0.800pt]
    (78.0814,180.2580) -- (100.1940,180.2580);
  \path[draw=c74b587,dash pattern=on 1.60pt,line join=round,line cap=round,draw
    opacity=0.803,line width=0.800pt] (78.0814,202.3710) -- (78.0814,180.2580);
  \path[draw=black,line join=round,line cap=round,line width=0.800pt]
    (100.1940,202.3710) -- (78.0814,202.3710);
  \path[draw=black,line join=round,line cap=round,line width=0.800pt]
    (144.4200,202.3710) -- (122.3070,202.3710);
  \path[draw=c74b587,dash pattern=on 1.60pt,line join=round,line cap=round,line
    width=0.800pt] (144.4200,180.2580) -- (144.4200,202.3710);
  \path[draw=black,line join=round,line cap=round,line width=0.800pt]
    (122.3070,180.2580) -- (144.4200,180.2580);
  \path[draw=black,line join=round,line cap=round,line width=0.800pt]
    (122.3070,202.3710) -- (122.3070,224.4840);
  \path[draw=black,line join=round,line cap=round,line width=0.800pt]
    (122.3070,158.1450) -- (122.3070,180.2580);
  \path[draw=black,line join=round,line cap=round,line width=0.800pt]
    (100.1940,180.2580) -- (100.1940,158.1450);
  \path[draw=black,line join=round,line cap=round,line width=0.800pt]
    (100.1940,224.4840) -- (100.1940,202.3710);
  \path[draw=cd9ac32,line join=round,line cap=round,draw opacity=0.796,line
    width=1.600pt] (100.1940,202.3710) -- (122.3070,180.2580);
      \path[draw=c002fff,line join=round,line cap=round,line width=0.800pt]
        (105.9050,191.3190) .. controls (103.0710,191.3220) and (79.9241,191.3520) ..
        (79.9241,191.3520);
      \path[fill=c002fff] (105.0300,191.3200) -- (102.4100,194.8230) --
        (109.4050,191.3140) -- (102.4010,187.8230) -- (105.0300,191.3200) -- cycle;
      \path[draw=cff001f,line join=round,line cap=round,line width=0.800pt]
        (111.3060,163.4900) .. controls (111.3000,166.3250) and (111.2510,189.4710) ..
        (111.2510,189.4710);
      \path[fill=cff001f] (111.3040,164.3650) -- (114.7990,166.9970) --
        (111.3140,159.9900) -- (107.7990,166.9820) -- (111.3040,164.3650) -- cycle;
      \path[draw=c002fff,line join=round,line cap=round,line width=0.800pt]
        (116.5950,191.3140) .. controls (119.4150,191.3140) and (142.4190,191.3140) ..
        (142.4190,191.3140);
      \path[fill=c002fff] (117.4700,191.3140) -- (120.0950,187.8140) --
        (113.0950,191.3140) -- (120.0950,194.8140) -- (117.4700,191.3140) -- cycle;
      \path[draw=cff001f,line join=round,line cap=round,line width=0.800pt]
        (111.2510,219.1380) .. controls (111.2510,216.3030) and (111.2510,193.1570) ..
        (111.2510,193.1570);
      \path[fill=cff001f] (111.2510,218.2630) -- (107.7510,215.6380) --
        (111.2510,222.6380) -- (114.7510,215.6380) -- (111.2510,218.2630) -- cycle;
  \path[draw=black,line join=round,line cap=round,line width=0.800pt]
    (109.4080,191.3140) .. controls (109.4080,190.2960) and (110.2330,189.4710) ..
    (111.2510,189.4710) .. controls (112.2690,189.4710) and (113.0940,190.2960) ..
    (113.0940,191.3140) .. controls (113.0940,192.3320) and (112.2690,193.1570) ..
    (111.2510,193.1570) .. controls (110.2330,193.1570) and (109.4080,192.3320) ..
    (109.4080,191.3140) -- cycle;
  \path[draw=black,dash pattern=on 0.80pt,line join=round,line cap=round,line
    width=0.800pt] (142.4190,191.3140) .. controls (142.4190,190.2960) and
    (143.2440,189.4710) .. (144.2620,189.4710) .. controls (145.2800,189.4710) and
    (146.1050,190.2960) .. (146.1050,191.3140) .. controls (146.1050,192.3320) and
    (145.2800,193.1570) .. (144.2620,193.1570) .. controls (143.2440,193.1570) and
    (142.4190,192.3320) .. (142.4190,191.3140) -- cycle;
  \path[draw=black,dash pattern=on 0.80pt,line join=round,line cap=round,line
    width=0.800pt] (109.4710,158.1450) .. controls (109.4710,157.1270) and
    (110.2960,156.3020) .. (111.3140,156.3020) .. controls (112.3310,156.3020) and
    (113.1560,157.1270) .. (113.1560,158.1450) .. controls (113.1560,159.1620) and
    (112.3310,159.9870) .. (111.3140,159.9870) .. controls (110.2960,159.9870) and
    (109.4710,159.1620) .. (109.4710,158.1450) -- cycle;
  \path[draw=black,dash pattern=on 0.80pt,line join=round,line cap=round,line
    width=0.800pt] (109.4080,224.4830) .. controls (109.4080,223.4660) and
    (110.2330,222.6410) .. (111.2510,222.6410) .. controls (112.2690,222.6410) and
    (113.0940,223.4660) .. (113.0940,224.4830) .. controls (113.0940,225.5010) and
    (112.2690,226.3260) .. (111.2510,226.3260) .. controls (110.2330,226.3260) and
    (109.4080,225.5010) .. (109.4080,224.4830) -- cycle;
  \path[draw=black,dash pattern=on 0.80pt,line join=round,line cap=round,line
    width=0.800pt] (76.2386,191.3520) .. controls (76.2386,190.3340) and
    (77.0636,189.5090) .. (78.0814,189.5090) .. controls (79.0991,189.5090) and
    (79.9241,190.3340) .. (79.9241,191.3520) .. controls (79.9241,192.3700) and
    (79.0991,193.1950) .. (78.0814,193.1950) .. controls (77.0636,193.1950) and
    (76.2386,192.3700) .. (76.2386,191.3520) -- cycle;
  \path[draw=black,line join=round,line cap=round,line width=0.800pt]
    (211.6370,169.5860) .. controls (216.6930,166.1090) and (226.6160,163.7450) ..
    (238.0300,163.7450) .. controls (249.3160,163.7450) and (259.1440,166.0560) ..
    (264.2510,169.4700);
  \path[draw=black,line join=round,line cap=round,line width=0.800pt]
    (240.3560,223.2940) .. controls (237.1460,218.2320) and (234.9690,208.4070) ..
    (234.9690,197.1170) .. controls (234.9690,193.1820) and (235.2340,189.4250) ..
    (235.7140,185.9850) .. controls (236.5190,186.0090) and (237.3320,186.0220) ..
    (238.1540,186.0220) .. controls (240.7680,186.0220) and (243.3050,185.8980) ..
    (245.7220,185.6650);
  \path[draw=black,line join=round,line cap=round,line width=0.800pt]
    (240.3860,170.8920) .. controls (239.3440,172.5220) and (238.4100,174.6540) ..
    (237.6230,177.1750) .. controls (239.3360,177.1750) and (241.0150,177.1220) ..
    (242.6500,177.0200);
  \path[draw=black,line join=round,line cap=round,line width=0.800pt]
    (236.9540,167.2990) -- (245.3130,167.2990);
  \path[draw=black,line join=round,line cap=round,line width=0.800pt]
    (236.9000,227.0310) -- (245.3130,226.9340);
  \path[draw=black,line join=round,line cap=round,line width=0.800pt]
    (245.3130,167.2990) .. controls (251.0250,167.2990) and (255.6560,180.6490) ..
    (255.6560,197.1170) .. controls (255.6560,213.5840) and (251.0250,226.9340) ..
    (245.3130,226.9340);
  \path[draw=black,line join=round,line cap=round,line width=0.800pt]
    (267.8520,166.0370) -- (267.9850,174.8360);
  \path[draw=black,line join=round,line cap=round,line width=0.800pt]
    (254.3580,175.3040) .. controls (262.4870,173.3210) and (267.8520,169.9110) ..
    (267.8520,166.0370) .. controls (267.8520,159.9120) and (254.4400,154.9470) ..
    (237.8960,154.9470) .. controls (221.3530,154.9470) and (207.9410,159.9120) ..
    (207.9410,166.0370) .. controls (207.9410,171.0700) and (216.9950,175.3190) ..
    (229.4020,176.6750) .. controls (231.2860,170.9590) and (233.9490,167.3960) ..
    (236.9000,167.3960) .. controls (242.6130,167.3960) and (247.2440,180.7460) ..
    (247.2440,197.2140) .. controls (247.2440,213.6820) and (242.6130,227.0310) ..
    (236.9000,227.0310) .. controls (231.1880,227.0310) and (226.5570,213.6820) ..
    (226.5570,197.2140) .. controls (226.5570,192.9440) and (226.8680,188.8840) ..
    (227.4290,185.2110) .. controls (216.1170,183.6270) and (208.0740,179.5790) ..
    (208.0740,174.8360) -- (207.9410,166.0370);
  \path[draw=black,line join=round,line cap=round,line width=0.800pt]
    (267.9850,174.8360) .. controls (267.9850,178.4820) and (263.2330,181.7170) ..
    (255.8970,183.7380);
  \path[draw=black,line join=round,line cap=round,line width=0.800pt]
    (230.9560,181.3590) .. controls (230.9560,180.6220) and (231.5540,180.0240) ..
    (232.2910,180.0240) .. controls (233.0290,180.0240) and (233.6270,180.6220) ..
    (233.6270,181.3590) .. controls (233.6270,182.0970) and (233.0290,182.6950) ..
    (232.2910,182.6950) .. controls (231.5540,182.6950) and (230.9560,182.0970) ..
    (230.9560,181.3590) -- cycle;
  \path[draw=black,dash pattern=on 0.80pt,line join=round,line cap=round,line
    width=0.800pt] (242.6740,158.0030) .. controls (243.4120,158.0070) and
    (244.0070,158.6090) .. (244.0030,159.3470) .. controls (243.9990,160.0840) and
    (243.3970,160.6790) .. (242.6600,160.6750) .. controls (241.9220,160.6710) and
    (241.3270,160.0690) .. (241.3310,159.3320) .. controls (241.3350,158.5940) and
    (241.9370,157.9990) .. (242.6740,158.0030) -- cycle;
  \path[draw=black,dash pattern=on 0.80pt,line join=round,line cap=round,line
    width=0.800pt] (249.3540,213.7240) .. controls (250.0910,213.7280) and
    (250.6860,214.3300) .. (250.6820,215.0680) .. controls (250.6780,215.8050) and
    (250.0760,216.4000) .. (249.3390,216.3960) .. controls (248.6010,216.3920) and
    (248.0060,215.7910) .. (248.0100,215.0530) .. controls (248.0140,214.3150) and
    (248.6160,213.7200) .. (249.3540,213.7240) -- cycle;
      \path[draw=c002fff,line join=round,line cap=round,line width=0.800pt]
        (227.5060,180.7670) .. controls (221.7820,179.6650) and (210.6880,176.9890) ..
        (208.0740,172.4230);
      \path[fill=c002fff] (226.6430,180.6190) -- (223.4640,183.6240) --
        (230.9550,181.3590) -- (224.6490,176.7250) -- (226.6430,180.6190) -- cycle;
      \path[draw=c002fff,line join=round,line cap=round,line width=0.800pt]
        (237.1270,181.4330) .. controls (239.7480,181.5030) and (242.7560,181.5710) ..
        (244.5550,181.3790);
      \path[fill=c002fff] (238.0020,181.4510) -- (240.7000,178.0070) --
        (233.6280,181.3590) -- (240.5530,185.0050) -- (238.0020,181.4510) -- cycle;
  \path[draw=c002fff,line join=round,line cap=round,line width=0.800pt]
    (267.9850,171.9690) .. controls (266.6390,176.7350) and (255.1380,179.5810) ..
    (255.1380,179.5810);
  \path[draw=cff001f,line join=round,line cap=round,line width=0.800pt]
    (232.2910,180.0240) .. controls (232.2910,180.0240) and (234.7980,170.1970) ..
    (239.0290,169.1960);
  \path[draw=cff001f,line join=round,line cap=round,line width=0.800pt]
    (232.2910,182.6950) .. controls (229.8100,195.7800) and (230.0920,219.6510) ..
    (238.9360,224.3440);
      \path[draw=cff001f,dash pattern=on 1.60pt,line join=round,line cap=round,line
        width=0.800pt] (250.3050,210.3550) .. controls (253.1870,197.7940) and
        (251.9350,175.8550) .. (242.5700,167.2990);
      \path[fill=cff001f] (250.5430,209.5130) -- (247.8880,206.0350) --
        (249.3540,213.7230) -- (254.6240,207.9380) -- (250.5430,209.5130) -- cycle;
      \path[draw=cff001f,dash pattern=on 1.60pt,line join=round,line cap=round,line
        width=0.800pt] (248.0630,219.6560) .. controls (246.8370,222.3750) and
        (245.6260,223.7450) .. (243.0090,226.9610);
      \path[fill=cff001f] (247.7450,220.4700) -- (250.0480,224.1900) --
        (249.3390,216.3960) -- (243.5290,221.6400) -- (247.7450,220.4700) -- cycle;
  \path[draw=c002fff,dash pattern=on 1.60pt,line join=round,line cap=round,line
    width=0.800pt] (244.0030,159.3470) .. controls (251.3020,159.3390) and
    (262.6650,162.9250) .. (266.0490,167.8210);
  \path[draw=c002fff,dash pattern=on 1.60pt,line join=round,line cap=round,line
    width=0.800pt] (241.3310,159.3320) .. controls (231.7830,159.3320) and
    (215.2190,159.3730) .. (210.0510,168.2450);
  \path[draw=c74b587,dash pattern=on 1.60pt,line join=round,line cap=round,draw
    opacity=0.801,line width=0.800pt] (242.6730,155.0870) -- (242.5830,163.8730);
  \path[draw=cf56356,dash pattern=on 1.60pt,line join=round,line cap=round,draw
    opacity=0.798,line width=0.800pt] (245.4110,214.1530) -- (253.1950,216.4140);
  \path[draw=cd9ac32,line join=round,line cap=round,draw opacity=0.801,line
    width=0.800pt] (227.4290,185.2110) -- (237.6230,177.1750);
  \path[fill=black,line width=0.600pt] (163.5000,193.8900) node[above right]
    (text14366) {$\implies$};

  \end{tikzpicture}

  \caption{2-fold covering with three branch points}
  \label{S_1_1_3}
\end{figure}

2-fold branched covering space of a disk with three branch points is a torus with one boundary component. It may be regarded as the union of two annuli with one branch point in common and $90^{\circ}$ twisted to each other. More generally, let $X$ be the 2-fold branched covering space of a disk with $k$ branch points, then it may be regarded as the union of $(k-1)$ sequence of annuli such that two neighboring annuli share one branch point. Let $g$ be the genus of $X$ and $b$ be the number of boundary components of $X$. Then we have the following:

\begin{table}[H]
  \centering
  \begin{tabular}{r|ccccccc}
    $k$ & 1 & 2 & 3 & 4 & 5 & 6 & 7 \\
    \hline
    $b$ & 1 & 2 & 1 & 2 & 1 & 2 & 1 \\
    \hline
    $g$ & 0 & 0 & 1 & 1 & 2 & 2 & 3
  \end{tabular}
  \caption{2-fold branched covering space with $k$ branch points.}
  \label{tab:2fold}
\end{table}

 For the case of 3-fold covering with more than two branch points, a similar argument works. If there are $k$ branch points,the 3-fold branched covering space of $S_{0,1,(k)}$ is the union of $(k-1)$ consecutive $S_{1,1,(2)}$'s and two neighboring atomic surfaces share one point.

\begin{figure}[H]
  \centering

  \definecolor{c74b587}{RGB}{116,181,135}
  \definecolor{c002fff}{RGB}{0,47,255}
  \definecolor{cf56356}{RGB}{245,99,86}
  \definecolor{cff001f}{RGB}{255,0,31}

  \begin{tikzpicture}[y=0.80pt, x=0.80pt, yscale=-1.000000, xscale=1.000000, inner sep=0pt, outer sep=0pt]
  \path[draw=c74b587,dash pattern=on 1.60pt,line join=round,line cap=round,line
    width=0.800pt] (323.5600,364.3310) -- (330.9310,349.5890) --
    (323.5600,334.8470);
  \path[draw=black,line join=round,line cap=round,line width=0.800pt]
    (332.7730,349.5890) .. controls (332.7730,348.5710) and (331.9480,347.7460) ..
    (330.9310,347.7460) .. controls (329.9130,347.7460) and (329.0880,348.5710) ..
    (329.0880,349.5890) .. controls (329.0880,350.6070) and (329.9130,351.4320) ..
    (330.9310,351.4320) .. controls (331.9480,351.4320) and (332.7730,350.6070) ..
    (332.7730,349.5890) -- cycle;
      \path[draw=c002fff,line join=round,line cap=round,line width=0.800pt]
        (284.6790,349.5890) .. controls (289.1030,349.5890) and (329.0880,349.5890) ..
        (329.0880,349.5890);
      \path[fill=c002fff] (285.5540,349.5890) -- (288.1790,346.0890) --
        (281.1790,349.5890) -- (288.1790,353.0890) -- (285.5540,349.5890) -- cycle;
  \path[draw=c74b587,dash pattern=on 1.60pt,line join=round,line cap=round,line
    width=0.800pt] (272.5040,303.4740) -- (256.0220,303.5590) --
    (246.2010,316.7950);
  \path[draw=black,line join=round,line cap=round,line width=0.800pt]
    (246.2010,316.7950) -- (262.8520,349.6740) -- (246.3370,382.6210);
  \path[draw=black,line join=round,line cap=round,line width=0.800pt]
    (272.5040,303.4740) -- (288.4990,335.0560) -- (323.5600,334.8470);
  \path[draw=black,line join=round,line cap=round,line width=0.800pt]
    (255.1890,301.9150) .. controls (254.2810,302.3750) and (253.9180,303.4830) ..
    (254.3780,304.3910) .. controls (254.8380,305.2990) and (255.9470,305.6620) ..
    (256.8550,305.2030) .. controls (257.7620,304.7430) and (258.1260,303.6340) ..
    (257.6660,302.7260) .. controls (257.2060,301.8180) and (256.0970,301.4550) ..
    (255.1890,301.9150) -- cycle;
  \path[draw=black,line join=round,line cap=round,line width=0.800pt]
    (278.5010,347.9450) .. controls (277.5940,348.4050) and (277.2300,349.5140) ..
    (277.6900,350.4220) .. controls (278.1500,351.3300) and (279.2590,351.6930) ..
    (280.1670,351.2330) .. controls (281.0750,350.7730) and (281.4380,349.6640) ..
    (280.9780,348.7560) .. controls (280.5180,347.8490) and (279.4090,347.4850) ..
    (278.5010,347.9450) -- cycle;
      \path[draw=c002fff,line join=round,line cap=round,line width=0.800pt]
        (276.9190,344.8200) .. controls (274.9200,340.8730) and (256.8550,305.2030) ..
        (256.8550,305.2030);
      \path[fill=c002fff] (276.5240,344.0400) -- (272.2150,343.2790) --
        (278.5000,347.9430) -- (278.4600,340.1170) -- (276.5240,344.0400) -- cycle;
  \path[draw=c74b587,dash pattern=on 1.60pt,line join=round,line cap=round,line
    width=0.800pt] (246.1740,382.8960) -- (256.0500,396.0920) --
    (272.5320,396.1090);
  \path[draw=black,line join=round,line cap=round,line width=0.800pt]
    (272.5320,396.1090) -- (288.4290,364.3310) -- (323.5600,364.3310);
  \path[draw=black,line join=round,line cap=round,line width=0.800pt]
    (255.2240,397.7390) .. controls (256.1340,398.1950) and (257.2410,397.8280) ..
    (257.6970,396.9180) .. controls (258.1530,396.0080) and (257.7860,394.9010) ..
    (256.8760,394.4450) .. controls (255.9660,393.9890) and (254.8590,394.3560) ..
    (254.4030,395.2660) .. controls (253.9470,396.1760) and (254.3140,397.2830) ..
    (255.2240,397.7390) -- cycle;
      \path[draw=c002fff,line join=round,line cap=round,line width=0.800pt]
        (276.7760,354.7440) .. controls (274.7930,358.7000) and (256.8760,394.4450) ..
        (256.8760,394.4450);
      \path[fill=c002fff] (276.3840,355.5260) -- (278.3370,359.4420) --
        (278.3450,351.6150) -- (272.0790,356.3050) -- (276.3840,355.5260) -- cycle;
  \path[draw=cf56356,dash pattern=on 1.60pt,line join=round,line cap=round,line
    width=0.800pt] (379.6060,364.0950) -- (372.2360,349.3530) --
    (379.6060,334.6110);
  \path[draw=black,line join=round,line cap=round,line width=0.800pt]
    (370.3930,349.3530) .. controls (370.3930,348.3350) and (371.2180,347.5100) ..
    (372.2350,347.5100) .. controls (373.2530,347.5100) and (374.0780,348.3350) ..
    (374.0780,349.3530) .. controls (374.0780,350.3700) and (373.2530,351.1950) ..
    (372.2350,351.1950) .. controls (371.2180,351.1950) and (370.3930,350.3700) ..
    (370.3930,349.3530) -- cycle;
      \path[draw=cff001f,line join=round,line cap=round,line width=0.800pt]
        (377.5810,349.3530) .. controls (382.0050,349.3530) and (421.9900,349.3530) ..
        (421.9900,349.3530);
      \path[fill=cff001f] (378.4560,349.3530) -- (381.0810,345.8530) --
        (374.0810,349.3530) -- (381.0810,352.8530) -- (378.4560,349.3530) -- cycle;
  \path[draw=cf56356,dash pattern=on 1.60pt,line join=round,line cap=round,line
    width=0.800pt] (430.6620,303.2380) -- (447.1440,303.3220) --
    (456.9650,316.5590);
  \path[draw=black,line join=round,line cap=round,line width=0.800pt]
    (456.9650,316.5590) -- (440.3140,349.4380) -- (456.8290,382.3850);
  \path[draw=black,line join=round,line cap=round,line width=0.800pt]
    (430.6620,303.2380) -- (414.6670,334.8200) -- (379.6060,334.6110);
  \path[draw=black,line join=round,line cap=round,line width=0.800pt]
    (447.9770,301.6780) .. controls (448.8850,302.1380) and (449.2480,303.2470) ..
    (448.7880,304.1550) .. controls (448.3280,305.0630) and (447.2200,305.4260) ..
    (446.3120,304.9660) .. controls (445.4040,304.5070) and (445.0400,303.3980) ..
    (445.5000,302.4900) .. controls (445.9600,301.5820) and (447.0690,301.2190) ..
    (447.9770,301.6780) -- cycle;
  \path[draw=black,line join=round,line cap=round,line width=0.800pt]
    (424.6650,347.7090) .. controls (425.5730,348.1690) and (425.9360,349.2770) ..
    (425.4760,350.1850) .. controls (425.0160,351.0930) and (423.9070,351.4560) ..
    (423.0000,350.9970) .. controls (422.0920,350.5370) and (421.7280,349.4280) ..
    (422.1880,348.5200) .. controls (422.6480,347.6120) and (423.7570,347.2490) ..
    (424.6650,347.7090) -- cycle;
      \path[draw=cff001f,line join=round,line cap=round,line width=0.800pt]
        (444.7290,308.0910) .. controls (442.7300,312.0380) and (424.6650,347.7090) ..
        (424.6650,347.7090);
      \path[fill=cff001f] (444.3340,308.8720) -- (446.2700,312.7950) --
        (446.3100,304.9690) -- (440.0250,309.6320) -- (444.3340,308.8720) -- cycle;
  \path[draw=cf56356,dash pattern=on 1.60pt,line join=round,line cap=round,line
    width=0.800pt] (456.9920,382.6600) -- (447.1160,395.8560) --
    (430.6340,395.8720);
  \path[draw=black,line join=round,line cap=round,line width=0.800pt]
    (430.6340,395.8720) -- (414.7370,364.0950) -- (379.6060,364.0950);
  \path[draw=black,line join=round,line cap=round,line width=0.800pt]
    (447.9420,397.5030) .. controls (447.0320,397.9590) and (445.9250,397.5910) ..
    (445.4690,396.6810) .. controls (445.0130,395.7720) and (445.3810,394.6640) ..
    (446.2900,394.2080) .. controls (447.2000,393.7520) and (448.3070,394.1200) ..
    (448.7640,395.0300) .. controls (449.2200,395.9400) and (448.8520,397.0470) ..
    (447.9420,397.5030) -- cycle;
      \path[draw=cff001f,line join=round,line cap=round,line width=0.800pt]
        (444.7210,391.0770) .. controls (442.7380,387.1220) and (424.8200,351.3770) ..
        (424.8200,351.3770);
      \path[fill=cff001f] (444.3290,390.2950) -- (440.0230,389.5170) --
        (446.2890,394.2060) -- (446.2810,386.3800) -- (444.3290,390.2950) -- cycle;
  \path[fill=black,line width=0.600pt] (267.1418,323.7614) node[above right]
    (text67550) {$a_1$};
  \path[fill=black,line width=0.600pt] (255.0000,380.6400) node[above right]
    (text67554) {$b_1$};
  \path[fill=black,line width=0.600pt] (296.2500,361.8900) node[above right]
    (text67558) {$c_1$};
  \path[fill=black,line width=0.600pt] (442.5000,380.6400) node[above right]
    (text67562) {$a_2$};
  \path[fill=black,line width=0.600pt] (386.2500,361.8900) node[above right]
    (text67566) {$c_2$};
  \path[fill=black,line width=0.600pt] (439.6418,328.5053) node[above right]
    (text67570) {$b_2$};

  \end{tikzpicture}

  \caption{Two $S_{1,1,(2)}$'s.}
  \label{fig:2S_1_1_2}
\end{figure}
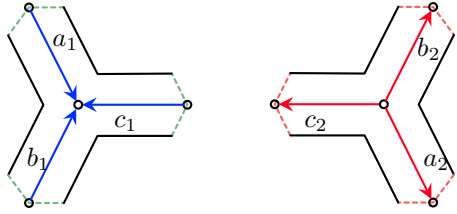

For better understanding we first consider the case with three branch points. By gluing two $S_{1,1,(2)}$'s, we have the following:

\begin{lemma}
A 3-fold branched covering space of a disk with three branch points is a torus with three boundary components.
\end{lemma}

\begin{figure}[H]
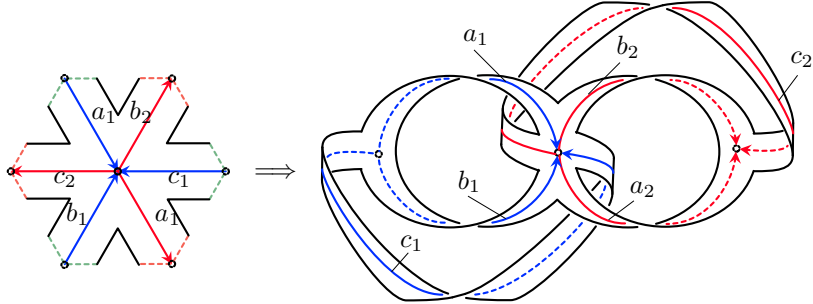

  \centering

  \definecolor{c002fff}{RGB}{0,47,255}
  \definecolor{c74b587}{RGB}{116,181,135}
  \definecolor{cff001f}{RGB}{255,0,31}
  \definecolor{cf56356}{RGB}{245,99,86}



  \caption{$S_{1,3,(3)}$ for three branch points}
  \label{fig:S_1_3_3}
\end{figure}

Figure \ref{fig:S_1_3_3} gives us more precise description of the covering space with three branch points.

Through gluing more atomic surfaces we can get the following result for the general case.

\begin{lemma}
 Let $X$ be the 3-fold branched covering space over a disk with $k$ branch points. Let $g$ be the genus and $b$ be the number of boundary components of $X$. Then we have
 $$g=\begin{cases}
  k-2 & \mbox{if }g\equiv 0\quad(\mbox{mod }3) \\
  k-1 & \mbox{otherwise}
\end{cases}
$$
$$
b = \begin{cases}
  3 & \mbox{if }g\equiv 0\quad(\mbox{mod }3) \\
  1 & \mbox{otherwise}
\end{cases}
$$

\end{lemma}

This lemma may also be expressed as in the following table.

\begin{table}[H]
  \centering
  \begin{tabular}{r|ccccccc}
    $k$ & 1 & 2 & 3 & 4 & 5 & 6 & 7 \\
    \hline
    $b$ & 1 & 1 & 3 & 1 & 1 & 3 & 1 \\
    \hline
    $g$ & 0 & 1 & 1 & 3 & 4 & 4 & 6
  \end{tabular}
  \caption{3-fold branched covering space with $k$ branch points.}
  \label{tab:3fold}
\end{table}

\section{A new non-geometric embedding and its homology triviality}

In this section we are going to construct the fundamental group homomorphism $\phi:B_k \rightarrow \Gamma_{g,b}$ induced by the space map $\Phi: \mathcal{C}_k \rightarrow \mathcal{M}_{g,b}$.  This homomorphism is defined through the lift, with respect to 3-fold branched covering, of half Dehn twist on the disk. We also show that
\begin{itemize}
  \item $\phi$ is injective.
  \item $\phi$ induces trivial homomorphism on homology groups in the stable range (Harer conjecture).
\end{itemize}
\subsection{The map $\phi:B_k \rightarrow \Gamma_{g,b}$ induced by 3-fold covering}
 First recall the case of 2-fold branched covering over a disk with two branch points, the half Dehn twist $\beta_i$ on a disk interchanging two branch points $\{p_i,p_{i+1}\}$ is lifted to the Dehn twist $\tilde\beta_i$ on an annulus(cf. \cite{Segal-Tillmann}).

\begin{figure}[H]
  \centering

  \definecolor{c2f3b73}{RGB}{47,59,115}
  \definecolor{c74b587}{RGB}{116,181,135}
  \definecolor{cf56356}{RGB}{245,99,86}
  \definecolor{c92b2de}{RGB}{146,178,222}

  \begin{tikzpicture}[y=0.80pt, x=0.80pt, yscale=-1.000000, xscale=1.000000, inner sep=0pt, outer sep=0pt, scale=0.8]
  \path[draw=black,line join=round,line cap=round,line width=0.800pt]
    (60.7257,443.5650) .. controls (60.7257,424.8240) and (86.3142,409.6300) ..
    (117.8790,409.6300) .. controls (149.4440,409.6300) and (175.0330,424.8240) ..
    (175.0330,443.5650) .. controls (175.0330,462.3070) and (149.4440,477.5000) ..
    (117.8790,477.5000) .. controls (86.3142,477.5000) and (60.7257,462.3070) ..
    (60.7257,443.5650) -- cycle;
  \path[draw=black,line join=round,line cap=round,line width=0.800pt]
    (92.8746,443.5650) .. controls (92.8746,442.5790) and (93.6743,441.7790) ..
    (94.6607,441.7790) .. controls (95.6471,441.7790) and (96.4467,442.5790) ..
    (96.4467,443.5650) .. controls (96.4467,444.5520) and (95.6471,445.3510) ..
    (94.6607,445.3510) .. controls (93.6743,445.3510) and (92.8746,444.5520) ..
    (92.8746,443.5650) -- cycle;
  \path[draw=black,line join=round,line cap=round,line width=0.800pt]
    (139.3120,443.5650) .. controls (139.3120,442.5790) and (140.1120,441.7790) ..
    (141.0980,441.7790) .. controls (142.0840,441.7790) and (142.8840,442.5790) ..
    (142.8840,443.5650) .. controls (142.8840,444.5520) and (142.0840,445.3510) ..
    (141.0980,445.3510) .. controls (140.1120,445.3510) and (139.3120,444.5520) ..
    (139.3120,443.5650) -- cycle;
      \path[draw=c2f3b73,line join=round,line cap=round,line width=0.800pt]
        (135.8110,443.5650) .. controls (131.7870,443.5650) and (96.4467,443.5650) ..
        (96.4467,443.5650);
      \path[fill=c2f3b73] (134.9360,443.5650) -- (132.3110,447.0650) --
        (139.3110,443.5650) -- (132.3110,440.0650) -- (134.9360,443.5650) -- cycle;
  \path[draw=c74b587,line join=round,miter limit=4.00,draw opacity=0.798,line
    width=1.800pt] (60.7257,443.5650) -- (92.8746,443.5650);
  \path[draw=cf56356,line join=round,miter limit=4.00,draw opacity=0.799,line
    width=1.200pt] (142.8840,443.5650) -- (175.0330,443.5650);
  \path[draw=c92b2de,line join=round,line cap=round,draw opacity=0.697,line
    width=2.427pt] (77.8194,443.5650) .. controls (77.8194,430.4290) and
    (95.7549,419.7800) .. (117.8790,419.7800) .. controls (140.0040,419.7800) and
    (157.9390,430.4290) .. (157.9390,443.5650) .. controls (157.9390,456.7020) and
    (140.0040,467.3510) .. (117.8790,467.3510) .. controls (95.7549,467.3510) and
    (77.8194,456.7020) .. (77.8194,443.5650) -- cycle(116.0930,464.4620) --
    (119.6650,467.4670) -- (116.0930,470.2400);
  \path[draw=black,line join=round,line cap=round,line width=0.800pt]
    (225.0420,443.5650) .. controls (225.0420,424.8240) and (250.6310,409.6300) ..
    (282.1960,409.6300) .. controls (313.7610,409.6300) and (339.3500,424.8240) ..
    (339.3500,443.5650) .. controls (339.3500,462.3070) and (313.7610,477.5000) ..
    (282.1960,477.5000) .. controls (250.6310,477.5000) and (225.0420,462.3070) ..
    (225.0420,443.5650) -- cycle;
  \path[draw=black,line join=round,line cap=round,line width=0.800pt]
    (257.1910,443.5650) .. controls (257.1910,442.5790) and (257.9910,441.7790) ..
    (258.9770,441.7790) .. controls (259.9640,441.7790) and (260.7630,442.5790) ..
    (260.7630,443.5650) .. controls (260.7630,444.5520) and (259.9640,445.3510) ..
    (258.9770,445.3510) .. controls (257.9910,445.3510) and (257.1910,444.5520) ..
    (257.1910,443.5650) -- cycle;
  \path[draw=black,line join=round,line cap=round,line width=0.800pt]
    (303.6290,443.5650) .. controls (303.6290,442.5790) and (304.4280,441.7790) ..
    (305.4150,441.7790) .. controls (306.4010,441.7790) and (307.2010,442.5790) ..
    (307.2010,443.5650) .. controls (307.2010,444.5520) and (306.4010,445.3510) ..
    (305.4150,445.3510) .. controls (304.4280,445.3510) and (303.6290,444.5520) ..
    (303.6290,443.5650) -- cycle;
      \path[draw=c2f3b73,line join=round,line cap=round,line width=0.800pt]
        (264.2640,443.5650) .. controls (268.2880,443.5650) and (303.6290,443.5650) ..
        (303.6290,443.5650);
      \path[fill=c2f3b73] (265.1390,443.5650) -- (267.7650,440.0650) --
        (260.7640,443.5650) -- (267.7640,447.0650) -- (265.1390,443.5650) -- cycle;
  \path[draw=c74b587,line join=round,miter limit=4.00,draw opacity=0.798,line
    width=1.800pt] (225.0420,443.5650) -- (240.4970,443.5650) .. controls
    (240.4470,477.5000) and (320.4640,477.5000) .. (320.4640,443.5650) --
    (307.2010,443.5650);
  \path[draw=cf56356,line join=round,miter limit=4.00,draw opacity=0.799,line
    width=1.200pt] (257.1910,443.5650) -- (243.7170,443.5650) .. controls
    (243.7170,409.6300) and (323.9330,409.6300) .. (323.9090,443.5650) --
    (339.3500,443.5650);
  \path[draw=black,line join=round,line width=0.800pt] (163.6460,247.3430) --
    (163.6460,343.8060);
  \path[draw=black,line join=round,line width=0.800pt] (76.8002,247.3430) --
    (76.8002,343.8060);
  \path[draw=black,line join=round,line cap=round,line width=0.800pt]
    (141.4470,290.4020) .. controls (141.4470,289.4160) and (142.2470,288.6160) ..
    (143.2330,288.6160) .. controls (144.2200,288.6160) and (145.0190,289.4160) ..
    (145.0190,290.4020) .. controls (145.0190,291.3890) and (144.2200,292.1880) ..
    (143.2330,292.1880) .. controls (142.2470,292.1880) and (141.4470,291.3890) ..
    (141.4470,290.4020) -- cycle;
  \path[draw=black,line join=round,line cap=round,line width=0.800pt]
    (95.0098,303.9700) .. controls (95.0098,302.9830) and (95.8094,302.1840) ..
    (96.7959,302.1840) .. controls (97.7823,302.1840) and (98.5819,302.9830) ..
    (98.5819,303.9700) .. controls (98.5819,304.9560) and (97.7823,305.7560) ..
    (96.7959,305.7560) .. controls (95.8094,305.7560) and (95.0098,304.9560) ..
    (95.0098,303.9700) -- cycle;
  \path[draw=black,line join=round,line width=0.800pt] (76.8002,247.3430) ..
    controls (76.8002,242.7610) and (96.2414,239.0460) .. (120.2230,239.0460) ..
    controls (144.2050,239.0460) and (163.6460,242.7610) .. (163.6460,247.3430) ..
    controls (163.6460,251.9250) and (144.2050,255.6390) .. (120.2230,255.6390) ..
    controls (96.2414,255.6390) and (76.8002,251.9250) .. (76.8002,247.3430) --
    cycle;
  \path[draw=black,dash pattern=on 1.60pt,line join=round,line cap=round,line
    width=0.800pt] (76.8002,343.8060) .. controls (76.8002,339.2240) and
    (96.2414,335.5090) .. (120.2230,335.5090) .. controls (144.2050,335.5090) and
    (163.6460,339.2240) .. (163.6460,343.8060);
  \path[draw=black,line join=round,line width=0.800pt] (163.6460,343.8060) ..
    controls (163.6460,348.3880) and (144.2050,352.1020) .. (120.2230,352.1020) ..
    controls (96.2414,352.1020) and (76.8002,348.3880) .. (76.8002,343.8060);
      \path[draw=c2f3b73,dash pattern=on 1.60pt,line join=round,line cap=round,line
        width=0.800pt] (137.8210,289.8720) .. controls (132.4430,289.4160) and
        (126.4880,289.1630) .. (120.2230,289.1630) .. controls (96.2414,289.1630) and
        (76.8002,292.8770) .. (76.8002,297.4590);
      \path[fill=c2f3b73] (136.9500,289.7890) -- (134.0050,293.0250) --
        (141.3050,290.2050) -- (134.6690,286.0560) -- (136.9500,289.7890) -- cycle;
  \path[draw=c2f3b73,line join=round,line cap=round,line width=0.800pt]
    (94.7277,304.1760) .. controls (83.8638,302.6680) and (76.8002,300.2210) ..
    (76.8002,297.4590);
  \path[draw=c2f3b73,line join=round,line cap=round,line width=0.800pt]
    (163.6460,297.4590) .. controls (163.6460,302.0410) and (144.2050,305.7560) ..
    (120.2230,305.7560) .. controls (112.2270,305.7560) and (104.7360,305.3430) ..
    (98.3017,304.6220);
      \path[draw=c2f3b73,dash pattern=on 1.60pt,line join=round,line cap=round,line
        width=0.800pt] (148.7170,291.1990) .. controls (157.8640,292.7190) and
        (163.6460,294.9600) .. (163.6460,297.4590);
      \path[fill=c2f3b73] (149.5820,291.3280) -- (152.6970,288.2570) --
        (145.2550,290.6790) -- (151.6590,295.1790) -- (149.5820,291.3280) -- cycle;
  \path[draw=c92b2de,dash pattern=on 1.60pt,line join=round,draw
    opacity=0.697,line width=2.400pt] (76.8002,268.8790) .. controls
    (76.8002,264.2970) and (96.2414,260.5830) .. (120.2230,260.5830) .. controls
    (144.2050,260.5830) and (163.6460,264.2970) .. (163.6460,268.8790);
  \path[draw=c92b2de,dash pattern=on 1.60pt,line join=round,draw
    opacity=0.697,line width=2.400pt] (76.8002,322.6420) .. controls
    (76.8002,318.0600) and (96.2414,314.3450) .. (120.2230,314.3450) .. controls
    (144.2050,314.3450) and (163.6460,318.0600) .. (163.6460,322.6420);
  \path[draw=c74b587,line join=round,miter limit=4.00,draw opacity=0.798,line
    width=1.800pt] (96.7959,305.7560) -- (96.8441,350.7950);
  \path[draw=c74b587,line join=round,miter limit=4.00,draw opacity=0.798,line
    width=1.800pt] (96.7610,254.3280) -- (96.7958,302.1840);
  \path[draw=cf56356,dash pattern=on 1.20pt,line join=round,miter limit=4.00,draw
    opacity=0.799,line width=1.200pt] (143.2330,288.6160) -- (143.3440,240.3210);
  \path[draw=cf56356,dash pattern=on 1.20pt,line join=round,miter limit=4.00,draw
    opacity=0.799,line width=1.200pt] (143.2330,292.1880) -- (143.2340,336.7710);
  \path[draw=c92b2de,line join=round,line cap=round,draw opacity=0.697,line
    width=2.427pt] (163.6460,268.8790) .. controls (163.6460,273.4610) and
    (144.2050,277.1750) .. (120.2230,277.1750) .. controls (96.2414,277.1750) and
    (76.8002,273.4610) .. (76.8002,268.8790)(120.2230,274.2860) --
    (116.6510,277.2910) -- (120.2230,280.0650);
  \path[draw=c92b2de,line join=round,line cap=round,draw opacity=0.697,line
    width=2.427pt] (163.6460,322.6420) .. controls (163.6460,327.2230) and
    (144.2050,330.9380) .. (120.2230,330.9380) .. controls (96.2414,330.9380) and
    (76.8002,327.2230) .. (76.8002,322.6420)(116.6510,328.0490) --
    (120.2230,331.0540) -- (116.6510,333.8270);
      \path[draw=black,line join=round,line cap=round,line width=0.800pt]
        (215.8810,443.1100) .. controls (213.1470,443.1290) and (182.3990,443.3470) ..
        (182.3990,443.3470);
      \path[fill=black] (213.0550,446.5170) -- (216.7000,443.4910) --
        (217.1700,443.1010) -- (216.6950,442.7170) -- (213.0070,439.7430) .. controls
        (212.7920,439.5700) and (212.4770,439.6030) .. (212.3040,439.8180) .. controls
        (212.1300,440.0330) and (212.1640,440.3480) .. (212.3790,440.5210) --
        (216.0670,443.4950) -- (216.0610,442.7210) -- (212.4160,445.7470) .. controls
        (212.2040,445.9240) and (212.1740,446.2390) .. (212.3510,446.4510) .. controls
        (212.5270,446.6640) and (212.8420,446.6930) .. (213.0550,446.5170) -- cycle;
      \path[draw=black,line join=round,line cap=round,line width=0.800pt]
        (222.6670,293.1460) .. controls (219.3900,293.1640) and (178.6940,293.3870) ..
        (178.6940,293.3870);
      \path[fill=black] (219.8360,296.5480) -- (223.4860,293.5280) --
        (223.9570,293.1390) -- (223.4820,292.7540) -- (219.7990,289.7750) .. controls
        (219.5840,289.6010) and (219.2690,289.6340) .. (219.0960,289.8490) .. controls
        (218.9220,290.0630) and (218.9550,290.3780) .. (219.1700,290.5520) --
        (222.8530,293.5320) -- (222.8490,292.7580) -- (219.1990,295.7780) .. controls
        (218.9860,295.9540) and (218.9560,296.2690) .. (219.1320,296.4820) .. controls
        (219.3080,296.6950) and (219.6230,296.7240) .. (219.8360,296.5480) -- cycle;
      \path[draw=black,line join=round,line cap=round,line width=0.800pt]
        (120.2230,393.1570) .. controls (120.2230,390.6790) and (120.2230,364.1630) ..
        (120.2230,364.1630);
      \path[fill=black] (116.8360,390.3070) -- (119.8360,393.9740) --
        (120.2230,394.4470) -- (120.6100,393.9740) -- (123.6100,390.3070) .. controls
        (123.7850,390.0930) and (123.7540,389.7780) .. (123.5400,389.6040) .. controls
        (123.3260,389.4290) and (123.0110,389.4600) .. (122.8360,389.6740) --
        (119.8360,393.3410) -- (120.6100,393.3410) -- (117.6100,389.6740) .. controls
        (117.4350,389.4600) and (117.1200,389.4290) .. (116.9070,389.6030) .. controls
        (116.6930,389.7780) and (116.6620,390.0930) .. (116.8360,390.3070) -- cycle;
      \path[draw=black,line join=round,line cap=round,line width=0.800pt]
        (282.1960,393.1580) .. controls (282.1960,390.6050) and (282.1960,362.8770) ..
        (282.1960,362.8770);
      \path[fill=black] (278.8090,390.3070) -- (281.8090,393.9740) --
        (282.1960,394.4470) -- (282.5830,393.9740) -- (285.5830,390.3070) .. controls
        (285.7580,390.0940) and (285.7260,389.7790) .. (285.5130,389.6040) .. controls
        (285.2990,389.4290) and (284.9840,389.4610) .. (284.8090,389.6740) --
        (281.8090,393.3410) -- (282.5830,393.3410) -- (279.5830,389.6740) .. controls
        (279.4080,389.4600) and (279.0930,389.4290) .. (278.8790,389.6040) .. controls
        (278.6660,389.7790) and (278.6340,390.0940) .. (278.8090,390.3070) -- cycle;
  \path[draw=black,line join=round,line width=0.800pt] (325.6040,247.3430) --
    (325.6040,343.8060);
  \path[draw=black,line join=round,line width=0.800pt] (238.7580,247.3430) --
    (238.7580,343.8060);
  \path[draw=black,line join=round,line cap=round,line width=0.800pt]
    (303.4050,290.4020) .. controls (303.4050,289.4160) and (304.2040,288.6160) ..
    (305.1910,288.6160) .. controls (306.1770,288.6160) and (306.9770,289.4160) ..
    (306.9770,290.4020) .. controls (306.9770,291.3890) and (306.1770,292.1880) ..
    (305.1910,292.1880) .. controls (304.2040,292.1880) and (303.4050,291.3890) ..
    (303.4050,290.4020) -- cycle;
  \path[draw=black,line join=round,line cap=round,line width=0.800pt]
    (256.9670,303.9700) .. controls (256.9670,302.9830) and (257.7670,302.1840) ..
    (258.7530,302.1840) .. controls (259.7400,302.1840) and (260.5390,302.9830) ..
    (260.5390,303.9700) .. controls (260.5390,304.9560) and (259.7400,305.7560) ..
    (258.7530,305.7560) .. controls (257.7670,305.7560) and (256.9670,304.9560) ..
    (256.9670,303.9700) -- cycle;
  \path[draw=black,line join=round,line width=0.800pt] (238.7580,247.3430) ..
    controls (238.7580,242.7610) and (258.1990,239.0460) .. (282.1810,239.0460) ..
    controls (306.1630,239.0460) and (325.6040,242.7610) .. (325.6040,247.3430) ..
    controls (325.6040,251.9250) and (306.1630,255.6390) .. (282.1810,255.6390) ..
    controls (258.1990,255.6390) and (238.7580,251.9250) .. (238.7580,247.3430) --
    cycle;
  \path[draw=black,dash pattern=on 1.60pt,line join=round,line cap=round,line
    width=0.800pt] (238.7580,343.8060) .. controls (238.7580,339.2240) and
    (258.1990,335.5090) .. (282.1810,335.5090) .. controls (306.1630,335.5090) and
    (325.6040,339.2240) .. (325.6040,343.8060);
  \path[draw=black,line join=round,line width=0.800pt] (325.6040,343.8060) ..
    controls (325.6040,348.3880) and (306.1630,352.1020) .. (282.1810,352.1020) ..
    controls (258.1990,352.1020) and (238.7580,348.3880) .. (238.7580,343.8060);
  \path[draw=c2f3b73,dash pattern=on 1.60pt,line join=round,line cap=round,line
    width=0.800pt] (238.7580,297.4590) .. controls (238.7580,292.8770) and
    (258.1990,289.1630) .. (282.1810,289.1630) .. controls (289.8320,289.1630) and
    (297.0220,289.5410) .. (303.2650,290.2050);
      \path[draw=c2f3b73,line join=round,line cap=round,line width=0.800pt]
        (253.2240,303.6420) .. controls (244.3460,302.1230) and (238.7580,299.9160) ..
        (238.7580,297.4590);
      \path[fill=c2f3b73] (252.3600,303.5080) -- (249.2320,306.5670) --
        (256.6840,304.1750) -- (250.2990,299.6490) -- (252.3600,303.5080) -- cycle;
      \path[draw=c2f3b73,line join=round,line cap=round,line width=0.800pt]
        (263.7440,304.9730) .. controls (269.3390,305.4750) and (275.5880,305.7560) ..
        (282.1810,305.7560) .. controls (306.1630,305.7560) and (325.6040,302.0410) ..
        (325.6040,297.4590);
      \path[fill=c2f3b73] (264.6140,305.0600) -- (267.5760,301.8410) --
        (260.2610,304.6230) -- (266.8760,308.8060) -- (264.6140,305.0600) -- cycle;
  \path[draw=c2f3b73,dash pattern=on 1.60pt,line join=round,line cap=round,line
    width=0.800pt] (325.6040,297.4590) .. controls (325.6040,294.6580) and
    (318.3380,292.1810) .. (307.2120,290.6790);
  \path[draw=cf56356,dash pattern=on 1.20pt,line join=round,miter limit=4.00,draw
    opacity=0.799,line width=1.200pt] (325.6040,268.8790) .. controls
    (325.6040,263.5020) and (305.2500,259.1400) .. (305.2500,259.1400) --
    (305.3020,240.3210);
  \path[draw=c74b587,dash pattern=on 1.80pt,line join=round,miter limit=4.00,draw
    opacity=0.798,line width=1.800pt] (305.1910,288.6160) -- (305.2500,266.6400)
    .. controls (287.9780,261.6510) and (238.7580,260.9660) ..
    (238.7580,268.8790);
  \path[draw=cf56356,line join=round,miter limit=4.00,draw opacity=0.799,line
    width=1.200pt] (325.6040,268.8790) .. controls (325.6040,281.9950) and
    (258.7500,279.3900) .. (258.7500,279.3900) -- (258.7530,302.1840);
  \path[draw=c74b587,line join=round,miter limit=4.00,draw opacity=0.798,line
    width=1.800pt] (258.7190,254.3280) -- (258.7500,272.6400) .. controls
    (258.7500,272.6400) and (238.7580,273.9990) .. (238.7580,268.8790);
  \path[draw=c74b587,line join=round,miter limit=4.00,draw opacity=0.799,line
    width=1.800pt] (325.6000,320.2680) .. controls (325.6000,333.3840) and
    (258.7500,331.1400) .. (258.7500,331.1400) -- (258.5450,350.7640);
  \path[draw=c74b587,dash pattern=on 1.80pt,line join=round,miter limit=4.00,draw
    opacity=0.799,line width=1.800pt] (325.6000,320.2680) .. controls
    (325.6000,314.8910) and (305.2500,309.3900) .. (305.2500,309.3900) --
    (305.2980,291.7110);
  \path[draw=cf56356,dash pattern=on 1.20pt,line join=round,miter limit=4.00,draw
    opacity=0.798,line width=1.200pt] (305.3273,335.6400) -- (305.2500,318.3900)
    .. controls (287.9780,313.4020) and (238.5930,315.4800) ..
    (238.5930,323.3940);
  \path[draw=cf56356,line join=round,miter limit=4.00,draw opacity=0.798,line
    width=1.200pt] (258.7340,305.7560) -- (258.7500,323.6400) .. controls
    (258.7500,323.6400) and (238.7730,325.4270) .. (238.7730,320.3070);
  \path[fill=black,line width=0.600pt] (124.8245,383.9441) node[above right]
    (text19984) {$p$};
  \path[fill=black,line width=0.600pt] (288.7500,383.9441) node[above right]
    (text19988) {$p$};
  \path[fill=black,line width=0.600pt] (195.0000,440.6400) node[above right]
    (text19992) {$\beta_i$};
  \path[fill=black,line width=0.600pt] (195.0000,290.6400) node[above right]
    (text19996) {$\tilde\beta_i$};

  \end{tikzpicture}

  \caption{half Dehn twist on disk is lifted to full Dehn twist on annulus.}
  \label{fig:halfDehn}
\end{figure}
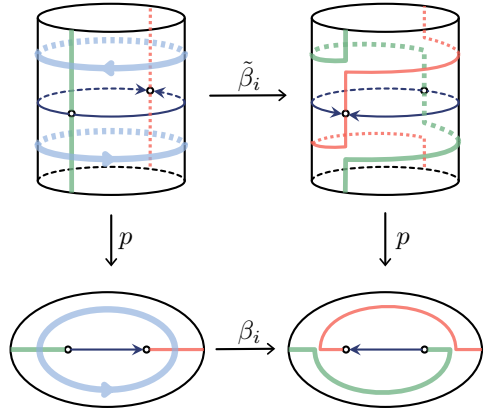

If we cut the annulus (2-fold covering space over a disk) along a middle circle into two annuli, then the lift $\tilde\beta_i$ is separated into two half Dehn twists in each of these two annuli.

 In the case of 3-fold covering, as the case of 2-fold covering, we should first start from finding the lift $\tilde\beta_i$ of half Dehn twist $\beta_i$ on a disk. We have found out that the 3-fold covering space (obtained by gluing three copies of disk) has one boundary component, that is, its the boundary forms a circle. However, if one  has one full trip along this boundary, it may be thought of as three times of circle-type trips around the two branch points. Therefore, lift $\tilde\beta_i$ of half Dehn twist $\beta_i$ may consist of three local half Dehn twists, one on each copy of three disks, That is, the lift $\tilde\beta_i$ may be regarded as a kind of 1/6 Dehn twist in the entire covering space which is described in the following figure.

\begin{figure}[H]
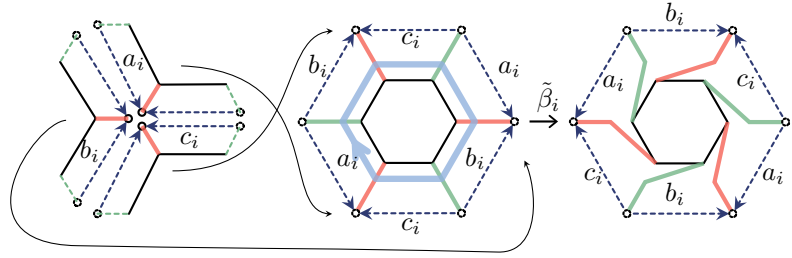

  \centering

  \definecolor{c2f3b73}{RGB}{47,59,115}
  \definecolor{cf56356}{RGB}{245,99,86}
  \definecolor{c74b587}{RGB}{116,181,135}
  \definecolor{c92b2de}{RGB}{146,178,222}



  \caption{1/6 Dehn twist}
  \label{fig:1/6Dehn}
\end{figure}

The lift $\tilde{\beta_i}$ shown in Figure \ref{fig:1/6Dehn} may be expressed as in the following theorem.

\begin{theorem}
For 3-fold branched covering over a disk with two branch points, the lift $\tilde\beta_i$ of half Dehn twist $\beta_i$ on a disk is a 1/6 Dehn twist on the covering space $S_{1,1,(2)}$ fixing the boundary pointwise. More precisely, it sends the path $a_i$ to the path $b_i^{-1}$, $b_i$ to $c_i^{-1}$, and $c_i$ to $a_i^{-1}$, where $a_i^{-1}, b_i^{-1}$ and $c_i^{-1}$ means the paths $a_i, b_i$ and $c_i$, respectively, with the reverse direction.
\end{theorem}

\begin{figure}[H]
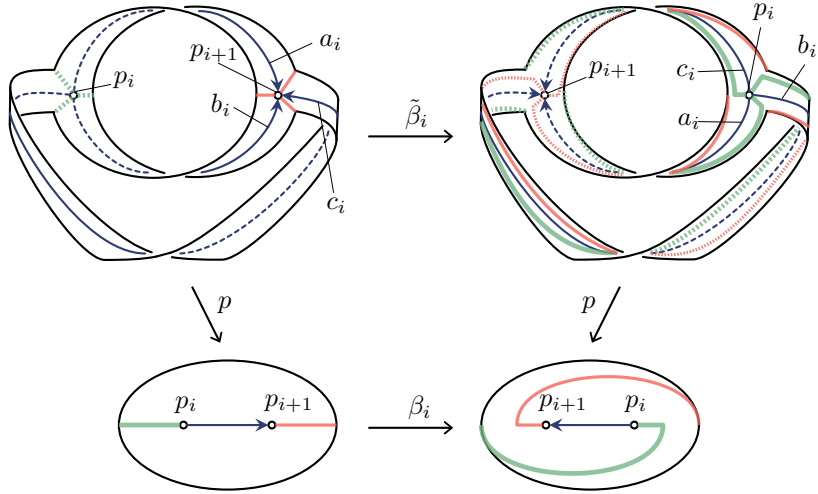

  \centering

  \definecolor{c2f3b73}{RGB}{47,59,115}
  \definecolor{cf56356}{RGB}{245,99,86}
  \definecolor{c74b587}{RGB}{116,181,135}



  \caption{Lift of the half Dehn twist}
  \label{fig:lifting}
\end{figure}

Since the self-homeomorphism $\tilde\beta_i:S_{1,1}\rightarrow S_{1,1}$ fixes the boundary, it can be extended to the entire surface via the identity map.

The map $$\phi: B_k \rightarrow \Gamma_{g,b}, \beta_i \mapsto \tilde{\beta_i}$$
is a well-defined homomorphism, that is, $\tilde{\beta_i}$'s satisfy the braid relation. Those readers who have good geometric vision can notice this fact by chasing the moves of branch points and arcs. On the other hand, the powerful Birman-Hilden theorem immediately tells us that $\phi$ is a homomorphism. Nevertheless, in the next section we will show this through explicit geometric and algebraic calculations, because of its own interest.

We now show that the map $\phi$ is indeed \textbf{injective} by exploiting Birman-Hilden theorem. See \cite{BH_Annals} and also see \cite{BH_Err}, \cite{GW16}) and the survey paper of Margalit and Winarski \cite{MW17}.

Let $p:X\rightarrow B$ be a branched covering where $B$ has the marked points as the branch points.
A homeomorphism $f:X\rightarrow X$ is said to be {\em fiber preserving} if for every pair of points $x,x'\in X$, $p(x)=p(x')$ implies $p(g(x))=p(g(x'))$; or simply, $f$ takes fibers to fibers.
If two fiber preserving homeomorphisms isotopic through fiber preserving homeomorphisms, then they are said to be {\em fiber-isotopic}.

A branched covering $p:X\rightarrow B$ is said to have the {\em Birman-Hilden property} if any two isotopic fiber preserving homeomorphisms are fiber-isotopic.

Let $L\Gamma_{0,1,(k)}\subset\Gamma_{0,1,(k)}$ be the {\em liftable mapping class group} consisting of the isotopy classes that have representatives that can be lifted to homeomorphisms of $S_{g,b}$ with respect to a branched covering $p$. And let $S\Gamma_{g,b}\subset\Gamma_{g,b}$ be the {\em symmetric mapping class group} consisting of the isotopy classes which can be represented by fiber preserving maps with respect to a branched covering $p$.

The following lemma is an immediate consequence Birman-Hilden theory (cf. Theorem 5 of \cite{BH_Annals} or Proposition 3.1 of \cite{MW17}).
\begin{lemma}\label{lemma:BH}
  Let $p:S_{g,b}\rightarrow S_{0,1,(k)}$ be a $n$-fold branched covering. Then $p$ has the Birman-Hilden property if and only if $L\Gamma_{0,1,(k)}\cong S\Gamma_{g,b}.$
\end{lemma}

Since the group of covering transformation is $\mathbb Z/3\mathbb Z$ and all branch points are fixed by the covering transformations, by Theorem 1 of \cite{BH_Annals}, $p$ has the Birman-Hilden property.
And we have seen that for 3-fold branched covering, every half Dehn twist $\beta_i$ can be lifted to $\tilde\beta_i$. Thus lemma \ref{lemma:BH} implies $$\Gamma_{0,1,(k)}=L\Gamma_{0,1,(k)}\cong S\Gamma_{g,b}\subset\Gamma_{g,b}.$$ Hence we have

\begin{theorem}
  Let $p:S_{g,b}\rightarrow S_{0,1,(k)}$ be a 3-fold branched covering. Then the assignment $\beta_i \mapsto \tilde\beta_i$ defines an injective homomorphism $B_k\hookrightarrow\Gamma_{g,b}$.
\end{theorem}

\subsection{The homology triviality of the embedding}

We show that the homology homomorphism $\phi_*:H_*(B_{\infty};\mathbb Z/p\mathbb Z)\rightarrow H_*(\Gamma_{\infty};\mathbb Z/p\mathbb Z)$ induced by $B_k \hookrightarrow \Gamma_{g,b} \hookrightarrow \Gamma_g$ is trivial. Note that by the Harer-Ivanov stability theorem the latter map of this composition induces homology isomorphism in the stable range. In the proof of the triviality, we take advantage of the fact that the map $\phi$ is defined through the fundamental group homomorphism of $\Phi: \mathcal{C}_{k}\rightarrow \M_{g,b}$. We show that the map $\Phi$ is compatible with naturally defined actions of the framed little 2-cube operad. This part is similar to that of  \cite{Segal-Tillmann}.

Let $\D=\{\D_k\}$ be the framed little 2-disks operad, with structure maps
$$\gamma:\D_k\times(\D_{m_1}\times\cdots\times\D_{m_k})\rightarrow \D_{\sum m_i}$$
given by composition of embeddings. (cf. \cite{Getzler94},\cite{SW01})

For the consistency of the number of boundary components, we may here deal with surfaces with $k\equiv 0$ (mod 3). Let $\X=\coprod_{m\geq 0}\X_m$ with $\X_m=\C_{3m}$. Then $\X$ is a $\D$-algebra defined by maps
$$\gamma_\X:\D_k\times(\X_{m_1}\times\cdots\times\X_{m_k})\rightarrow\X_{\sum m_i}$$
which $(f;P_1,\ldots,P_k)\mapsto f(P_1\cup\cdots\cup P_k)$.

\begin{figure}[H]
  \centering

  \begin{tikzpicture}[y=0.80pt, x=0.80pt, yscale=-1.000000, xscale=1.000000, inner sep=0pt, outer sep=0pt, scale=0.8]
  \path[draw=black,line join=round,line width=0.800pt] (156.5200,255.3670) ..
    controls (156.5200,246.3350) and (192.6460,246.3350) .. (192.6460,255.3670);
  \path[draw=black,line join=round,line width=0.800pt] (261.9990,190.4130) ..
    controls (261.9990,183.1750) and (277.9560,177.3070) .. (297.6400,177.3070) ..
    controls (317.3240,177.3070) and (333.2810,183.1750) .. (333.2810,190.4130) ..
    controls (333.2810,197.6510) and (317.3240,203.5180) .. (297.6400,203.5180) ..
    controls (277.9560,203.5180) and (261.9990,197.6510) .. (261.9990,190.4130) --
    cycle;
  \path[draw=black,line join=round,line width=0.800pt] (85.2382,255.3670) ..
    controls (85.2382,208.2870) and (261.9990,237.4920) .. (261.9990,190.4130);
  \path[draw=black,line join=round,line width=0.800pt] (333.2810,190.4130) ..
    controls (333.2810,231.6290) and (484.8120,214.1500) .. (484.8120,255.3670);
  \path[draw=black,line join=round,line width=0.800pt] (395.2670,248.5930) ..
    controls (404.3650,248.5600) and (413.5310,250.8170) .. (413.5310,255.3670);
  \path[draw=black,line join=round,line width=0.800pt] (156.5200,255.3670) ..
    controls (156.5200,262.6050) and (140.5630,268.4720) .. (120.8790,268.4720) ..
    controls (101.1950,268.4720) and (85.2382,262.6050) .. (85.2382,255.3670);
  \path[draw=black,dash pattern=on 1.60pt,line join=round,line width=0.800pt]
    (85.2382,255.3670) .. controls (85.2382,248.1280) and (101.1950,242.2610) ..
    (120.8790,242.2610) .. controls (140.5630,242.2610) and (156.5200,248.1280) ..
    (156.5200,255.3670);
  \path[draw=black,line join=round,line width=0.800pt] (263.9280,255.3670) ..
    controls (263.9280,262.6050) and (247.9710,268.4720) .. (228.2870,268.4720) ..
    controls (208.6030,268.4720) and (192.6460,262.6050) .. (192.6460,255.3670);
  \path[draw=black,dash pattern=on 1.60pt,line join=round,line width=0.800pt]
    (192.6460,255.3670) .. controls (192.6460,248.1280) and (208.6030,242.2610) ..
    (228.2870,242.2610) .. controls (247.9710,242.2610) and (263.9280,248.1280) ..
    (263.9280,255.3670);
  \path[draw=black,line join=round,line width=0.800pt] (484.8120,255.3670) ..
    controls (484.8120,262.6050) and (468.8550,268.4720) .. (449.1710,268.4720) ..
    controls (429.4880,268.4720) and (413.5310,262.6050) .. (413.5310,255.3670);
  \path[draw=black,dash pattern=on 1.60pt,line join=round,line width=0.800pt]
    (413.5310,255.3670) .. controls (413.5310,248.1280) and (429.4880,242.2610) ..
    (449.1710,242.2610) .. controls (468.8550,242.2610) and (484.8120,248.1280) ..
    (484.8120,255.3670);
  \path[draw=black,dash pattern=on 1.60pt,line join=round,line width=0.800pt]
    (85.2382,337.9820) .. controls (85.2382,326.4520) and (101.1950,317.1040) ..
    (120.8790,317.1040) .. controls (140.5630,317.1040) and (156.5200,326.4520) ..
    (156.5200,337.9820) .. controls (156.5200,349.5130) and (140.5630,358.8600) ..
    (120.8790,358.8600) .. controls (101.1950,358.8600) and (85.2382,349.5130) ..
    (85.2382,337.9820) -- cycle;
  \path[draw=black,dash pattern=on 1.60pt,line join=round,line width=0.800pt]
    (193.1310,337.9820) .. controls (193.1310,326.4520) and (209.0880,317.1040) ..
    (228.7720,317.1040) .. controls (248.4560,317.1040) and (264.4130,326.4520) ..
    (264.4130,337.9820) .. controls (264.4130,349.5130) and (248.4560,358.8600) ..
    (228.7720,358.8600) .. controls (209.0880,358.8600) and (193.1310,349.5130) ..
    (193.1310,337.9820) -- cycle;
  \path[draw=black,dash pattern=on 1.60pt,line join=round,line width=0.800pt]
    (413.5310,337.9820) .. controls (413.5310,326.4520) and (429.4880,317.1040) ..
    (449.1710,317.1040) .. controls (468.8550,317.1040) and (484.8120,326.4520) ..
    (484.8120,337.9820) .. controls (484.8120,349.5130) and (468.8550,358.8600) ..
    (449.1710,358.8600) .. controls (429.4880,358.8600) and (413.5310,349.5130) ..
    (413.5310,337.9820) -- cycle;
  \path[draw=black,line join=round,line width=0.800pt] (156.5200,416.8600) ..
    controls (156.5200,425.8920) and (192.6460,425.8920) .. (192.6460,416.8600);
  \path[draw=black,line join=round,line width=0.800pt] (85.2382,416.8600) ..
    controls (85.2382,463.9400) and (261.9990,434.7350) .. (261.9990,481.8140);
  \path[draw=black,line join=round,line width=0.800pt] (333.2810,481.8140) ..
    controls (333.2810,440.5980) and (484.8120,458.0770) .. (484.8120,416.8600);
  \path[draw=black,line join=round,line width=0.800pt] (396.8360,423.6170) ..
    controls (405.3890,423.4000) and (413.5310,421.1480) .. (413.5310,416.8600);
  \path[draw=black,line join=round,line width=0.800pt] (156.5200,416.8600) ..
    controls (156.5200,424.0980) and (140.5630,429.9660) .. (120.8790,429.9660) ..
    controls (101.1950,429.9660) and (85.2382,424.0980) .. (85.2382,416.8600);
  \path[draw=black,dash pattern=on 1.60pt,line join=round,line width=0.800pt]
    (85.2382,416.8600) .. controls (85.2382,409.6220) and (101.1950,403.7550) ..
    (120.8790,403.7550) .. controls (140.5630,403.7550) and (156.5200,409.6220) ..
    (156.5200,416.8600);
  \path[draw=black,line join=round,line width=0.800pt] (263.9280,416.8600) ..
    controls (263.9280,424.0980) and (247.9710,429.9660) .. (228.2870,429.9660) ..
    controls (208.6030,429.9660) and (192.6460,424.0980) .. (192.6460,416.8600);
  \path[draw=black,dash pattern=on 1.60pt,line join=round,line width=0.800pt]
    (192.6460,416.8600) .. controls (192.6460,409.6220) and (208.6030,403.7550) ..
    (228.2870,403.7550) .. controls (247.9710,403.7550) and (263.9280,409.6220) ..
    (263.9280,416.8600);
  \path[draw=black,line join=round,line width=0.800pt] (484.8120,416.8600) ..
    controls (484.8120,424.0980) and (468.8550,429.9660) .. (449.1710,429.9660) ..
    controls (429.4880,429.9660) and (413.5310,424.0980) .. (413.5310,416.8600);
  \path[draw=black,dash pattern=on 1.60pt,line join=round,line width=0.800pt]
    (413.5310,416.8600) .. controls (413.5310,409.6220) and (429.4880,403.7550) ..
    (449.1710,403.7550) .. controls (468.8550,403.7550) and (484.8120,409.6220) ..
    (484.8120,416.8600);
  \path[draw=black,dash pattern=on 1.60pt,line join=round,line width=0.800pt]
    (156.5200,255.3670);
  \path[draw=black,line join=round,line width=0.800pt] (85.2382,416.8600) ..
    controls (76.0355,377.5500) and (76.0354,294.6770) .. (85.2382,255.3670);
  \path[draw=black,line join=round,line width=0.800pt] (156.5200,255.3670) ..
    controls (165.7230,294.6770) and (165.7230,377.5500) .. (156.5200,416.8600);
  \path[draw=black,line join=round,line width=0.800pt] (192.6460,416.8600) ..
    controls (183.4430,377.5500) and (183.4430,294.6770) .. (192.6460,255.3670);
  \path[draw=black,line join=round,line width=0.800pt] (263.9280,255.3670) ..
    controls (273.1300,294.6770) and (273.1300,377.5500) .. (263.9280,416.8600);
  \path[draw=black,line join=round,line width=0.800pt] (413.5310,416.8600) ..
    controls (404.3280,377.5500) and (404.3280,294.6770) .. (413.5310,255.3670);
  \path[draw=black,line join=round,line width=0.800pt] (484.8120,255.3670) ..
    controls (494.0150,294.6770) and (494.0150,377.5500) .. (484.8120,416.8600);
  \path[draw=black,dash pattern=on 1.60pt,line join=round,line width=0.800pt]
    (280.4500,353.9400) .. controls (286.6960,362.8560) and (298.0200,368.8140) ..
    (310.9520,368.8140) .. controls (330.6360,368.8140) and (346.5930,355.0100) ..
    (346.5930,337.9820) .. controls (346.5930,320.9540) and (330.6360,307.1500) ..
    (310.9520,307.1500) .. controls (298.1720,307.1500) and (286.9630,312.9690) ..
    (280.6730,321.7100) .. controls (277.2740,326.4330) and (275.3110,332.0090) ..
    (275.3110,337.9820) .. controls (275.3110,343.8230) and (277.1880,349.2850) ..
    (280.4500,353.9400) -- cycle;
  \path[draw=black,line join=round,line width=0.800pt] (263.9280,255.3670) ..
    controls (263.9280,250.8040) and (273.1480,248.5460) .. (282.2730,248.5940);
  \path[draw=black,line join=round,line width=0.800pt] (263.9280,416.8600) ..
    controls (263.9280,421.3920) and (273.0250,423.6500) .. (282.0890,423.6340);
  \path[draw=black,fill=black,line join=round,line width=0.800pt]
    (316.9530,339.0660) .. controls (316.9530,338.0480) and (317.7780,337.2230) ..
    (318.7960,337.2230) .. controls (319.8130,337.2230) and (320.6390,338.0480) ..
    (320.6390,339.0660) .. controls (320.6390,340.0840) and (319.8130,340.9090) ..
    (318.7960,340.9090) .. controls (317.7780,340.9090) and (316.9530,340.0840) ..
    (316.9530,339.0660) -- cycle;
  \path[draw=black,fill=black,line join=round,line width=0.800pt]
    (339.0660,339.0660) .. controls (339.0660,338.0480) and (339.8910,337.2230) ..
    (340.9090,337.2230) .. controls (341.9260,337.2230) and (342.7510,338.0480) ..
    (342.7510,339.0660) .. controls (342.7510,340.0840) and (341.9260,340.9090) ..
    (340.9090,340.9090) .. controls (339.8910,340.9090) and (339.0660,340.0840) ..
    (339.0660,339.0660) -- cycle;
  \path[draw=black,fill=black,line join=round,line width=0.800pt]
    (361.1790,339.0660) .. controls (361.1790,338.0480) and (362.0040,337.2230) ..
    (363.0220,337.2230) .. controls (364.0390,337.2230) and (364.8640,338.0480) ..
    (364.8640,339.0660) .. controls (364.8640,340.0840) and (364.0390,340.9090) ..
    (363.0220,340.9090) .. controls (362.0040,340.9090) and (361.1790,340.0840) ..
    (361.1790,339.0660) -- cycle;
  \path[draw=black,dash pattern=on 1.60pt,line join=round,line width=0.800pt]
    (333.2810,481.8140) .. controls (333.2810,474.5760) and (317.3240,468.7080) ..
    (297.6400,468.7080) .. controls (277.9560,468.7080) and (261.9990,474.5760) ..
    (261.9990,481.8140);
  \path[draw=black,line join=round,line width=0.800pt] (261.9990,481.8140) ..
    controls (261.9990,489.0520) and (277.9560,494.9200) .. (297.6400,494.9200) ..
    controls (317.3240,494.9200) and (333.2810,489.0520) .. (333.2810,481.8140);
  \path[draw=black,dash pattern=on 1.60pt,line join=round,line width=0.800pt]
    (165.5570,313.1820) .. controls (121.3080,312.8440) and (81.7872,307.5350) ..
    (81.7872,340.9090) .. controls (81.7872,372.3540) and (120.9540,362.9460) ..
    (164.5630,362.2080);
  \path[draw=black,line join=round,line width=0.800pt] (164.5630,362.2080) ..
    controls (167.8820,362.1520) and (171.2270,362.1460) .. (174.5830,362.2120) ..
    controls (177.8400,362.2760) and (180.9770,362.3630) .. (184.0060,362.4710);
  \path[draw=black,dash pattern=on 1.60pt,line join=round,line width=0.800pt]
    (184.0060,362.4710) .. controls (223.7260,363.8780) and (244.8450,368.7780) ..
    (272.5950,371.4730);
  \path[draw=black,line join=round,line width=0.800pt] (272.5950,371.4730) ..
    controls (283.7590,372.5580) and (295.9970,373.2860) .. (310.9520,373.2860) ..
    controls (365.9100,373.2860) and (351.8020,372.0500) .. (395.2670,364.3420) ..
    controls (398.1510,363.8310) and (401.0610,363.4590) .. (403.9820,363.1970);
  \path[draw=black,dash pattern=on 1.60pt,line join=round,line width=0.800pt]
    (403.9820,363.1970) .. controls (445.0870,359.5170) and (488.2630,377.7720) ..
    (488.2630,339.0660) .. controls (488.2630,299.6970) and (445.1600,315.7330) ..
    (404.0870,311.8060);
  \path[draw=black,line join=round,line width=0.800pt] (404.0870,311.8060) ..
    controls (401.1310,311.5240) and (398.1850,311.1380) .. (395.2670,310.6200) ..
    controls (351.8020,302.9120) and (365.8050,302.0890) .. (310.9520,302.0890) ..
    controls (296.4490,302.0890) and (284.4850,302.8110) .. (273.5670,303.8870);
  \path[draw=black,dash pattern=on 1.60pt,line join=round,line width=0.800pt]
    (273.5670,303.8870) .. controls (245.4390,306.6590) and (224.2530,311.7820) ..
    (184.4980,312.9640);
  \path[draw=black,line join=round,line width=0.800pt] (184.4980,312.9640) ..
    controls (181.3160,313.0590) and (178.0160,313.1280) .. (174.5830,313.1690) ..
    controls (171.5600,313.2050) and (168.5480,313.2050) .. (165.5570,313.1820);
  \path[fill=black,line width=0.600pt] (95,297.0264) node[above right]
    (text1567) {$3m_1-2$};
  \path[fill=black,line width=0.600pt] (200,297.3082) node[above right]
    (text1571) {$3m_2-2$};
  \path[fill=black,line width=0.600pt] (423,298.4855) node[above right]
    (text1575) {$3m_k-2$};

  \end{tikzpicture}

  \caption{$\D$ acting on moduli spaces.}
  \label{fig:D_moduli}
\end{figure}
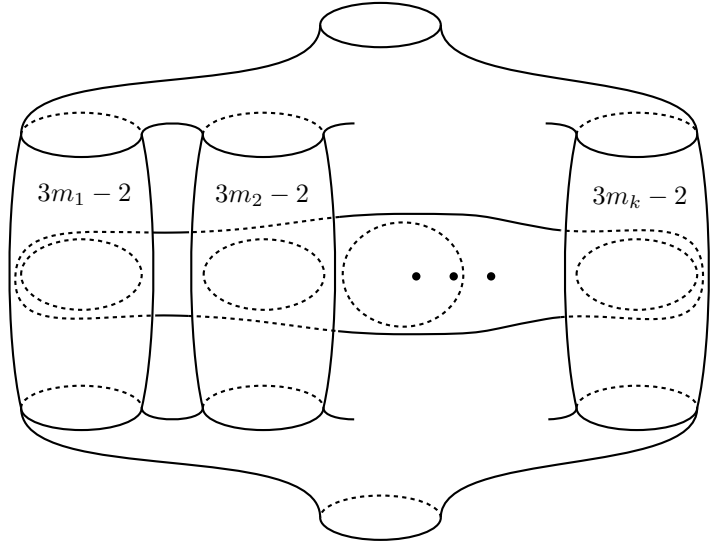

Let $\Y_m=\M_{3m-2,3}$ for $m\geq 1$ and $\Y_0=\M_{0,1}\sqcup\M{0,1}\sqcup\M{0,1}$. Each surface $\Sigma\in\Y_m$ has three parametrised boundaries $\partial_a\Sigma,\partial_b\Sigma$ and $\partial_c\Sigma$.
For $f\in\D_k$ and $j\in\{a,b,c\}$, let $(D_f)_j=D\setminus F(D\cup\cdots\cup D)=S_{0,k+1}$.
Then $\Y=\coprod_{m\geq 0}\Y_m$ is $\D$-algebra, defined by maps
$$(f;\Sigma_1,\ldots,\Sigma_k)\mapsto ((D_f)_a\cup (D_f)_b\cup (D_f)_c\cup\Sigma_1\cup\cdots\cup\Sigma_k)/\equiv$$
where the $\cup_i\partial_j\Sigma_i$ is identified with the interior boundaries of $(D_f)_j$ for $j\in\{a,b,c\}$ (see Figure \ref{fig:D_moduli}).

We have that map $\Phi$ gives us the map $\X\rightarrow\Y$ of $\D$-algebras. This gives us the proof of the following strong version of Harer conjecture(theorem).

\begin{theorem}
The homology homomorphism $$\phi_*:H_n(B_{\infty}; R)\rightarrow H_n(\Gamma_{\infty};R)$$
  is zero for all $n\geq 1$ in any constant coefficient $R$.
\end{theorem}

By Harer-Ivanov homology stability theorem, we may have a modified version of this theorem as follows:
$$\phi_*:H_n(B_{k}; R)\rightarrow H_n(\Gamma_{g};R)$$
is zero for all $g > 2n\geq 1$ in any constant coefficient $R$. Here, $g=k-1$ or $k-2$.

\section{The analysis of the embedding $\phi$}

In this section we analyse the action of the lift $\tilde{\beta_i}$ on the fundamental group of the surface and give algebraic descriptions of it. From these we could check that our construction of the map $\phi$ is indeed correct, that is, $\tilde{\beta_i}$'s  satisfy the braid relations. We believe that the calculations and the results are of their own interest. Moreover, We show that the lift $\tilde{\beta_i}$ is the product of two inverse Dehn twists, that is, $\phi$ is a non-geometric embedding.

\subsection{Action of $\tilde{\beta_i}$'s on fundamental group}
 We adopt the notations of fundamental groupoid with a set of base points. A  fundamental groupoid $\Pi_1(X,A)$ with the set of base points $A\subset X$ is the category whose the object set is $A$ and the morphisms from $x$ to $y$ are the homotopy classes of paths from $x$ to $y$ in a path connected space $X$. Note that for $x_0\in A$, $\hom_{\Pi_1(X,A)}(x_0,x_0) = \pi_1(X,x_0)$.

Let $P=\{p_1,\ldots,p_k\}$ be the marked points of $S_{0,1,(k)}$ and let $p_0, p_{k+1}$ be two distinct points on the boundary. $\mathcal G = \Pi_1\left(S_{0,1,(k)},\{p_0,p_{k+1}\}\cup P\right)$ is the groupoid generated by the following data:
$$\begin{array}{rcl}
\mbox{objects} & = & \{p_0,\ldots,p_{k+1}\}, \\
\mbox{morphisms} & \mbox{are} & \mbox{generated by the arrows }\left(p_i\xrightarrow{a_i} p_{i+1}\right)_{0\leq i\leq k}.
\end{array}$$

%
%

A half Dehn twist $\beta_i$ may be regarded as self-functor of $\mathcal{G}$ such that
$$\beta_i:\left\{
\begin{array}{rcl}
  p_i & \mapsto & p_{i+1}, \\
  p_{i+1} & \mapsto & p_i, \\
  \\
  a_{i-1} & \mapsto &  a_{i-1} a_i,\\
  a_i & \mapsto & a_i^{-1}, \\
  a_{i+1} & \mapsto & a_i a_{i+1},\\
\end{array}\right.$$
This functor fixes points and paths that do not appear in the list.

Let $p:X\rightarrow S_{0,1,(k)}$ be the 3-fold branched covering with branch points $P=\{p_1,\ldots,p_k\}$.
Then for a path from $p_{i}$ to $p_{i+1}$ in $S_{0,1,(k)}$ is lifted to three paths in $X$.
Then there is a groupoid $\tilde{\mathcal G}$ on the space $X$ which is the lift of $\mathcal{G}$ with respect to the covering map $p$. That is $\tilde{\mathcal G}$ is generated by :
$$\begin{array}{rcl}
\mbox{objects} & = & \{a_{0}(0), b_{0}(0), c_{0}(0), a_{k}(1), b_{k}(1), c_{k}(1)\}\cup P, \\
\\
\mbox{morphisms} & \mbox{are} & \mbox{generated by the morphisms} \\
\end{array}
$$
$$
\left\{\begin{array}{cl}
\left(a_0(0)\xrightarrow{a_0} p_{1}\right), \left(b_0(0)\xrightarrow{b_0} p_{1}\right), \left(c_0(0)\xrightarrow{c_0} p_{1}\right),\\

\left(p_i\xrightarrow{a_i} p_{i+1}\right),\left(p_i\xrightarrow{b_i} p_{i+1}\right),\left(p_i\xrightarrow{c_i} p_{i+1}\right), & \mbox{ for }1\leq i\leq (k-1),\\

\left(p_k\xrightarrow{a_{k}} a_{k}(1)\right), \left(p_k\xrightarrow{b_{k}} b_{k}(1)\right), \left(p_k\xrightarrow{c_{k}} c_{k}(1)\right).
\end{array}\right.
$$

For morphism $\left(p_i\xrightarrow i p_{i+1}\right)$ in $\mathcal{G}$, denote its three lifts by $a_i, b_i, c_i$ in $\tilde{\mathcal{G}}$ .

%
%
%


Let $\tilde\beta_i$ be the lift of $\beta_i$, then it may be regarded as a self-functor of $\tilde{\mathcal{G}}$ such that
$$\tilde\beta_i:\left\{
\begin{array}{rcl}
  p_i & \mapsto & p_{i+1}, \\
  p_{i+1} & \mapsto & p_i, \\
  \\
  a_{i-1} & \mapsto & a_{i-1} c_{i} \\
  b_{i-1} & \mapsto & b_{i-1} a_{i} \\
  c_{i-1} & \mapsto & c_{i-1} b_{i} \\
  \\
  a_{i} & \mapsto & b_{i}^{-1} \\
  b_{i} & \mapsto & c_{i}^{-1} \\
  c_{i} & \mapsto & a_{i}^{-1} \\
  \\
  a_{i+1} & \mapsto & c_{i} a_{i+1} \\
  b_{i+1} & \mapsto & a_{i} b_{i+1} \\
  c_{i+1} & \mapsto & b_{i} c_{i+1} \\
\end{array}\right.$$

\begin{remark}
  We may check that the lifts $\tilde{\beta_i}$'s, as self-functors of $\tilde{\mathcal{G}}$, satisfy the braid relation $$\tilde\beta_i\tilde\beta_{i+1}\tilde\beta_i = \tilde\beta_{i+1}\tilde\beta_{i}\tilde\beta_{i+1}$$
  so that the homomorphism $\phi: B_k \rightarrow \Gamma_{g,b}$ is correctly constructed.

    On the objects of $\tilde{\mathcal{G}}$, it is obvious that the braid relations are satisfied.

We are going to show that the equality $\tilde\beta_i\tilde\beta_{i+1}\tilde\beta_i(a_{j}) = \tilde\beta_{i+1}\tilde\beta_{i}\tilde\beta_{i+1}(a_{j})$ holds for each morphism $a_j$.

For the morphism $a_{i-1}$,
    $$a_{i-1}\xmapsto{\tilde\beta_{i}} a_{i-1}c_i\xmapsto{\tilde\beta_{i+1}} a_{i-1}c_i b_{i+1}\xmapsto{\tilde\beta_{i}} (a_{i-1}c_i) (a_i^{-1}) (a_ib_{i+1}) = a_{i-1}c_ib_{i+1} ,$$
    $$a_{i-1}\xmapsto{\tilde\beta_{i+1}} a_{i-1}\xmapsto{\tilde\beta_{i}} a_{i-1}c_i\xmapsto{\tilde\beta_{i+1}} a_{i-1}c_ib_{i+1}.$$

For the morphism $a_{i}$,
    $$a_{i}\xmapsto{\tilde\beta_{i}} b_{i}^{-1}\xmapsto{\tilde\beta_{i+1}} a_{i+1}^{-1}b_{i}^{-1}\xmapsto{\tilde\beta_{i}} (a_{i+1}^{-1}c_i^{-1})(c_i)=a_{i+1}^{-1},$$
    $$a_{i}\xmapsto{\tilde\beta_{i+1}} a_{i}c_{i+1}\xmapsto{\tilde\beta_{i}} b_{i}^{-1}b_{i}c_{i+1}\xmapsto{\tilde\beta_{i+1}} a_{i+1}^{-1}.$$

For the morphism $a_{i+1}$,
    $$a_{i+1}\xmapsto{\tilde\beta_{i}} c_{i}a_{i+1}\xmapsto{\tilde\beta_{i+1}} c_{i}b_{i+1}b_{i+1}^{-1}\xmapsto{\tilde\beta_{i}} a_{i}^{-1},$$
    $$a_{i+1}\xmapsto{\tilde\beta_{i+1}} b_{i+1}^{-1}\xmapsto{\tilde\beta_{i}} b_{i+1}^{-1}a_{i}^{-1}\xmapsto{\tilde\beta_{i+1}} (c_{i+1})(c_{i+1}^{-1}a_{i}^{-1})=a_{i}^{-1}.$$

For the morphism $a_{i+2}$,
    $$a_{i+2}\xmapsto{\tilde\beta_{i}} a_{i+2}\xmapsto{\tilde\beta_{i+1}} c_{i+1}a_{i+2}\xmapsto{\tilde\beta_{i}} b_{i}c_{i+1}a_{i+2},$$
    $$a_{i+2}\xmapsto{\tilde\beta_{i+1}} c_{i+1}a_{i+2}\xmapsto{\tilde\beta_{i}} b_{i}c_{i+1}a_{i+2}\xmapsto{\tilde\beta_{i+1}} (b_{i}a_{i+1})(a_{i+1}^{-1})(c_{i+1}a_{i+2})=b_{i}c_{i+1}a_{i+2}.$$

Thus we have $\tilde\beta_i\tilde\beta_{i+1}\tilde\beta_i = \tilde\beta_{i+1}\tilde\beta_{i}\tilde\beta_{i+1}$.
\end{remark}

On the other hand, we can also prove that the braid relation is satisfied with respect to the fundamental group of the surface. Note that each element in the mapping class group $\Gamma_{g,b}$ is completed determined by the action on the fundamental group of the surface. Let $x$ and $y$ be loops in $X$ with base points on the boundary given by:
$$x_i = (c_0 \cdots c_{i-1})b_i a_i^{-1}(c_0 \cdots c_{i-1})^{-1}, y_i = (c_0 \cdots c_{i-1})a_i c_i^{-1}(c_0 \cdots c_{i-1})^{-1}$$
Note that $\{x_i,y_i\mid 1\leq i\leq k-1\}$ is a generator set of $\pi_1(X)$ which is a free group.

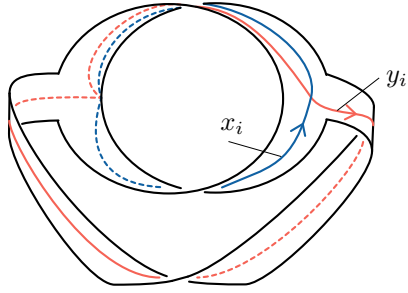
\begin{figure}[H]
  \centering

  \definecolor{cf56356}{RGB}{245,99,86}
  \definecolor{c1061a3}{RGB}{16,97,163}

  \begin{tikzpicture}[y=0.80pt, x=0.80pt, yscale=-1.000000, xscale=1.000000, inner sep=0pt, outer sep=0pt]
  \path[draw=black,line join=round,line cap=round,line width=0.800pt]
    (67.3718,129.2400);
  \path[draw=black,line join=round,line cap=round,line width=0.800pt]
    (215.9180,166.3280) -- (216.8280,144.9930);
  \path[draw=black,line join=round,line cap=round,line width=0.800pt]
    (216.8160,145.3030) .. controls (217.4300,136.6910) and (194.4010,129.5490) ..
    (194.4010,129.5490);
  \path[draw=cf56356,line join=round,line cap=round,line width=0.800pt]
    (132.5500,97.6423) .. controls (156.3540,98.1384) and (172.6500,119.6870) ..
    (187.0830,140.4930) .. controls (194.3450,150.9610) and (216.8160,145.3030) ..
    (216.4540,153.7690);
  \path[draw=cf56356,dash pattern=on 1.60pt,line join=round,line cap=round,line
    width=0.800pt] (47.5613,144.2710) .. controls (52.0237,138.9760) and
    (87.9408,140.7990) .. (87.9408,140.7990) .. controls (68.5749,126.9640) and
    (100.5970,93.9354) .. (125.3560,96.6220);
  \path[draw=c1061a3,line join=round,line cap=round,line width=0.800pt]
    (144.9810,183.1670) .. controls (167.4350,173.3160) and (178.9270,173.8140) ..
    (187.0830,140.4930) .. controls (191.4110,122.8120) and (164.2940,99.9968) ..
    (136.3340,96.2256);
  \path[draw=black,line join=round,line cap=round,line width=0.800pt]
    (194.5080,151.0890) .. controls (188.9050,170.8780) and (167.4960,185.6110) ..
    (141.9400,185.6110) .. controls (140.2130,185.6110) and (138.5050,185.5430) ..
    (136.8200,185.4120);
  \path[draw=black,line join=round,line cap=round,line width=0.800pt]
    (194.3450,129.9500);
  \path[draw=black,line join=round,line cap=round,line width=0.800pt]
    (67.6718,129.4300) .. controls (73.9065,110.3920) and (94.8966,96.3906) ..
    (119.8310,96.3906) .. controls (149.6780,96.3906) and (173.8730,116.4540) ..
    (173.8730,141.2020) .. controls (173.8730,165.9510) and (149.6780,186.0140) ..
    (119.8310,186.0140) .. controls (94.3293,186.0140) and (72.9534,171.3680) ..
    (67.2693,151.6670);
  \path[draw=black,line join=round,line cap=round,line width=0.800pt]
    (194.3450,129.9500) .. controls (188.5000,110.4430) and (167.2540,95.9871) ..
    (141.9400,95.9871) .. controls (140.0470,95.9871) and (138.1770,96.0679) ..
    (136.3340,96.2256);
  \path[draw=black,line join=round,line cap=round,line width=0.800pt]
    (67.2693,151.6670);
  \path[draw=black,line join=round,line cap=round,line width=0.800pt]
    (44.5336,138.6020) .. controls (44.4694,143.3840) and (44.2768,157.7300) ..
    (44.2768,157.7300);
  \path[draw=black,line join=round,line cap=round,line width=0.800pt]
    (50.4613,153.0920) .. controls (55.3112,151.8810) and (62.1423,151.5880) ..
    (67.6718,151.4800);
  \path[draw=black,line join=round,line cap=round,line width=0.800pt]
    (129.4120,227.7240) .. controls (136.6120,228.8420) and (144.9710,229.0870) ..
    (155.8340,229.0350) .. controls (169.1590,228.9700) and (215.2510,181.9400) ..
    (215.9180,166.3280) .. controls (216.4540,153.7690) and (194.3450,150.9610) ..
    (194.3450,150.9610);
  \path[draw=black,line join=round,line cap=round,line width=0.800pt]
    (44.5336,138.6020) .. controls (44.5336,151.3710) and (65.8883,186.5960) ..
    (94.2228,210.0490) .. controls (103.5300,217.7520) and (110.8420,222.4040) ..
    (118.8870,225.1740);
  \path[draw=black,line join=round,line cap=round,line width=0.800pt]
    (67.2693,129.2430) .. controls (55.0645,129.1900) and (44.5336,129.8220) ..
    (44.5336,138.6020);
  \path[draw=black,line join=round,line cap=round,line width=0.800pt]
    (121.3590,227.3850) .. controls (114.1210,228.8370) and (105.7790,229.1810) ..
    (94.2228,229.1810) .. controls (80.5255,229.1810) and (44.4052,175.1370) ..
    (44.2768,157.7300);
  \path[draw=black,line join=round,line cap=round,line width=0.800pt]
    (208.6490,158.7760) .. controls (208.6490,158.7760) and (184.0070,187.6800) ..
    (155.8340,209.1150) .. controls (144.7940,217.5140) and (137.0960,222.4560) ..
    (128.8940,225.3330);
  \path[draw=cf56356,line join=round,line cap=round,line width=0.800pt]
    (44.3588,151.6220) .. controls (45.0187,162.7670) and (81.7105,225.7970) ..
    (115.4090,225.7920);
  \path[draw=cf56356,dash pattern=on 1.60pt,line join=round,line cap=round,line
    width=0.800pt] (211.7240,161.6320) .. controls (208.6490,176.5380) and
    (163.7430,225.1030) .. (133.2850,225.7920);
  \path[draw=black,line join=round,line cap=round,line width=0.800pt]
    (121.3590,227.3850) -- (128.8940,225.3330);
  \path[draw=cf56356,line join=round,line cap=round,line width=0.800pt]
    (203.9090,144.2710) -- (208.2540,148.4350) -- (203.0870,150.9610);
  \path[draw=c1061a3,line join=round,line cap=round,line width=0.800pt]
    (177.5880,156.4680) -- (183.1560,153.1260) -- (184.2970,158.7760);
  \path[draw=c1061a3,dash pattern=on 1.60pt,line join=round,line cap=round,line
    width=0.800pt] (125.3620,98.1384) .. controls (97.6552,100.3270) and
    (82.0102,124.5710) .. (87.9408,140.7990) -- (87.9408,140.7990) .. controls
    (81.7105,160.3030) and (84.8911,181.8880) .. (117.1820,183.0790);
  \path[draw=black,line join=round,line cap=round,line width=0.800pt]
    (125.5450,183.5080) .. controls (103.7380,177.7480) and (87.9408,160.8040) ..
    (87.9408,140.7990) .. controls (87.9408,120.8490) and (103.6500,103.9450) ..
    (125.3620,98.1384);
  \path[draw=black,line join=miter,line cap=butt,miter limit=4.00,line
    width=0.300pt] (199.4339,146.7881) -- (219.2984,131.3379);
  \path[fill=black,line width=0.600pt] (144.3520,157.8038) node[above right]
    (text29155) {$x_i$};
  \path[fill=black,line width=0.600pt] (222.6092,135.5946) node[above right]
    (text29159) {$y_i$};
  \path[draw=black,line join=miter,line cap=butt,miter limit=4.00,line
    width=0.300pt] (148.0442,161.4217) -- (173.0155,168.5564);

  \end{tikzpicture}

  \caption{$x_i$ and $y_i$.}
  \label{fig:xy}
\end{figure}

On the fundamental group, we have
$$
\tilde\beta_{i}:\left\{\begin{array}{rcll}
    x_{i-1} & \mapsto & x_{i-1} \cdot y_{i-1} \cdot y_i \cdot y_{i-1}^{-1} & \mbox{ if } i\geq 2, \\
    y_{i-1} & \mapsto & y_{i-1} \cdot y_i^{-1} \cdot x_i^{-1} & \mbox{ if } i\geq 2, \\
    \\
    x_i & \mapsto & x_i \cdot y_i, & \\
    y_i & \mapsto & x_i^{-1}, & \\
    \\
    x_{i+1} & \mapsto & x_i \cdot y_i \cdot x_{i+1} \cdot x_i^{-1} & \mbox{ if } i\leq k-2, \\
    y_{i+1} & \mapsto & x_i \cdot y_{i+1} \cdot y_i^{-1} \cdot x_i^{-2} & \mbox{ if } i\leq k-2,\\
    \\
    x_{i+j} & \mapsto & x_i^2 \cdot y_i \cdot x_{i+j} \cdot y_i^{-1} \cdot x_i^{-2}  & \mbox{ if } j\geq 2 \mbox{ and } i+j\leq k-1, \\
    x_{i+j} & \mapsto & x_i^2 \cdot y_i \cdot y_{i+j} \cdot y_i^{-1} \cdot x_i^{-2} & \mbox{ if } j\geq 2 \mbox{ and } i+j\leq k-1, \\
\end{array}\right.
$$
From this result (action on the fundamental group) and some calculations, by hands or some computer program,  we can again confirm that the braid relations for $\tilde{\beta_i}$'s are satisfied.


\subsection{$\tilde{\beta_i}$ in terms of full Dehn twists}

Let $X=S_{g,b}$ be the branched covering space over a disk. Let $x_i, y_i$ be the loops on $X$ given as above. Let $z_i = (c_0 \cdots c_{i-1})c_i b_i^{-1}(c_0 \cdots c_{i-1})^{-1} = (x_iy_i)^{-1}$ for $1\leq i \leq k-1$. The Dehn on $X$ along the loops $x_i, y_i, z_i$ are denoted by $D_{x_i}, D_{y_i}, D_{z_i}$, respectively. The actions of these Dehn twists on the groupoid $\tilde{\mathcal{G}}$ are as follows:
$$D_{x_i}:\left\{
\begin{array}{rcl}
  p_i & \mapsto & p_{i+1}, \\
  p_{i+1} & \mapsto & p_i, \\
  \\
  a_{i-1} & \mapsto & a_{i-1} a_{i} \\
  b_{i-1} & \mapsto & b_{i-1} a_{i} \\
  c_{i-1} & \mapsto & c_{i-1} b_{i} \\
  \\
  a_{i} & \mapsto & b_{i}^{-1} \\
  b_{i} & \mapsto & a_{i}^{-1} \\
  c_{i} & \mapsto & a_{i}^{-1} c_i a_{i}^{-1} \\
  \\
  a_{i+1} & \mapsto & a_{i} a_{i+1} \\
  b_{i+1} & \mapsto & a_{i} b_{i+1} \\
  c_{i+1} & \mapsto & b_{i} c_{i+1} \\
\end{array}\right.
D_{y_i}:\left\{
\begin{array}{rcl}
  p_i & \mapsto & p_{i+1}, \\
  p_{i+1} & \mapsto & p_i, \\
  \\
  a_{i-1} & \mapsto & a_{i-1} a_{i} \\
  b_{i-1} & \mapsto & b_{i-1} c_{i} \\
  c_{i-1} & \mapsto & c_{i-1} a_{i} \\
  \\
  a_{i} & \mapsto & c_{i}^{-1} \\
  b_{i} & \mapsto & a_{i}^{-1} b_i a_{i}^{-1} \\
  c_{i} & \mapsto & a_{i}^{-1} \\
  \\
  a_{i+1} & \mapsto & a_{i} a_{i+1} \\
  b_{i+1} & \mapsto & c_{i} b_{i+1} \\
  c_{i+1} & \mapsto & a_{i} c_{i+1} \\
\end{array}\right.
D_{z_i}:\left\{
\begin{array}{rcl}
  p_i & \mapsto & p_{i+1}, \\
  p_{i+1} & \mapsto & p_i, \\
  \\
  a_{i-1} & \mapsto & a_{i-1} b_{i} \\
  b_{i-1} & \mapsto & b_{i-1} c_{i} \\
  c_{i-1} & \mapsto & c_{i-1} c_{i} \\
  \\
  a_{i} & \mapsto & c_{i}^{-1} a_i c_{i}^{-1} \\
  b_{i} & \mapsto & c_{i}^{-1} \\
  c_{i} & \mapsto & b_{i}^{-1} \\
  \\
  a_{i+1} & \mapsto & b_{i} a_{i+1} \\
  b_{i+1} & \mapsto & c_{i} b_{i+1} \\
  c_{i+1} & \mapsto & c_{i} c_{i+1} \\
\end{array}\right.$$

The actions on $\pi_1 X$ are as follows:
$$D_{x_i}:\left\{\begin{array}{rcll}
  y_{i-1} & \mapsto & y_{i-1} \cdot x_i^{-1} & \mbox{ if } i\geq 2, \\
  \\
  y_{i} & \mapsto & y_i \cdot x_i^{-1}, \\
  \\
  x_{i+1} & \mapsto & x_i \cdot x_{i+1} \cdot x_i^{-1} & \mbox{ if } i\leq k-2, \\
  y_{i+1} & \mapsto & x_i \cdot y_{i+1} \cdot y_i^{-1} \cdot x_i^{-1} \cdot y_i \cdot x_i^{-1} & \mbox{ if } i\leq k-2 \\
  \\
  x_{i+j} & \mapsto & x_i \cdot y_i^{-1} \cdot x_i \cdot y_i \cdot x_{i+j} \cdot y_i^{-1} \cdot x_i^{-1} \cdot y_i \cdot x_i^{-1} & \mbox{ if } j\geq 2 \mbox{ and } i+j \leq k-1 \\
  y_{i+j} & \mapsto & x_i \cdot y_i^{-1} \cdot x_i \cdot y_i \cdot y_{i+j} \cdot y_i^{-1} \cdot x_i^{-1} \cdot y_i \cdot x_i^{-1} & \mbox{ if } j\geq 2 \mbox{ and } i+j \leq k-1 \\
\end{array}\right.$$
$$D_{y_i}:\left\{\begin{array}{rcll}
  x_{i-1} & \mapsto & x_{i-1} \cdot y_{i-1} \cdot y_i^{-1} \cdot y_{i-1}^{-1} & \mbox{ if } i\geq 2, \\
  \\
  x_{i} & \mapsto & x_i \cdot y_i^{-1}, \\
  \\
  x_{i+1} & \mapsto & x_{i+1} \cdot y_i^{-1} & \mbox{ if } i\leq k-2, \\
  \\
  y_{i+j} & \mapsto & y_i \cdot y_{i+1} \cdot y_i^{-1} & \mbox{ if } j\geq 1 \mbox{ and } i+j\leq k-1 \\
  x_{i+j} & \mapsto & y_i \cdot x_{i+j} \cdot y_i^{-1} & \mbox{ if } j\geq 2 \mbox{ and } i+j \leq k-1 \\
\end{array}\right.$$
$$D_{z_i}:\left\{\begin{array}{rcll}
  x_{i-1} & \mapsto & x_{i-1} \cdot y_{i-1} \cdot y_i^{-1} \cdot x_i^{-1} \cdot y_{i-1}^{-1} & \mbox{ if } i\geq 2, \\
  y_{i-1} & \mapsto & y_{i-1} \cdot x_i \cdot y_i & \mbox{ if } i\geq 2, \\
  \\
  x_{i} & \mapsto & y_i^{-1}, \\
  y_{i} & \mapsto & y_i \cdot x_i \cdot y_i, \\
  \\
  x_{i+1} & \mapsto & y_i^{-1} \cdot x_i^{-1} \cdot x_{i+1} & \mbox{ if } i\leq k-2, \\
  y_{i+1} & \mapsto & y_{i+1} \cdot x_i \cdot y_i & \mbox{ if } i\leq k-2, \\
  \\
  x_{i+j} & \mapsto & y_i^{-1} \cdot x_i^{-1} \cdot x_{i+j} \cdot x_i \cdot y_i & \mbox{ if } j\geq 2 \mbox{ and } i+j \leq k-1 \\
  y_{i+j} & \mapsto & y_i^{-1} \cdot x_i^{-1} \cdot y_{i+j} \cdot x_i \cdot y_i & \mbox{ if } j\geq 2 \mbox{ and } i+j\leq k-1 \\
\end{array}\right.$$

From the result of calculations (by hands or some computer program), we can show that $\tilde{\beta_i}$ equals the product of the inverse of two Dehn twists.

\begin{theorem}
 $\tilde\beta_i$ is equal to $D_{z_i}^{-1}\cdot D_{y_i}^{-1}$ as an element of $\Gamma_{g,b}$ .
\end{theorem}


\vskip 0.5cm

\begin{thebibliography}{99}

  \bibitem{BFSV03}
    Balteanu, C. and Fiedorowicz, Z. and SchwÃ€nzl, R. and Vogt, R.
    \textit{Iterated monoidal categories}.
    Advances in Mathematics, 176.2 (2003), pp. 277--349.

  \bibitem{BH_Annals}
    Birman, Joan S. and Hilden, Hugh M.
    \textit{On Isotopies of Homeomorphisms of Riemann Surfaces}.
    Annals of Mathematics, 97.3 (1973), pp. 424--439.

  \bibitem{BH_Err}
    Birman, Joan S. and Hilden, Hugh M.
    \textit{Erratum to `Isotopies of homeomorphisms of Riemann surfaces'}.
    Annals of Mathematics, 185.1 (2017), pp. 345--345.

  \bibitem{BT12}
    {B{\"o}digheimer}, C.-F. and {Tillmann}, U.
    \textit{Embeddings of braid groups into mapping class groups and their homology}.
    ArXiv e-prints http://adsabs.harvard.edu/abs/2012arXiv1204.4310B (2012).

  \bibitem{Getzler94}
    {Getzler}, E.
    \textit{Batalin-Vilkovisky algebras and two-dimensional topological field theories}.
    Communications in Mathematical Physics, 159 (1994), pp. 265--285.

  \bibitem{GW16}
    {Ghaswala}, T. and {Winarski}, R.~R.
    \textit{Lifting Homeomorphisms and Cyclic Branched Covers of Spheres}.
    ArXiv e-prints, http://adsabs.harvard.edu/abs/2016arXiv160706060G (2016).

  \bibitem{JS13}
    JEONG, CHAN-SEOK and SONG, YONGJIN
    \textit{THE PILLAR SWITCHINGS OF MAPPING CLASS GROUPS OF SURFACES}.
    International Journal of Mathematics, 24.13 (2013), p. 1350103.

  \bibitem{MacLane65}
    MacLane, S.
    \textit{Categorical algebra}.
    Bull. Amer. Math. Soc., 71.1 (1965), pp. 40--106.

  \bibitem{MW17}
    {Margalit}, D. and {Winarski}, R.~R.
    \textit{The Birman-Hilden theory}.
    ArXiv e-prints, http://adsabs.harvard.edu/abs/2017arXiv170303448M (2017).

  \bibitem{Segal74}
    Segal, Graeme.
    \textit{Categories and cohomology theories}
    Topology, 13.3 (1974), pp. 293--312.

  \bibitem{Segal04}
    Segal, Graeme.
    \textit{The definition of conformal field theory}.
    Topology, Geometry and Quantum Field Theory: Proceedings of the 2002 Oxford Symposium in Honour of the 60th Birthday of Graeme Segal, (2004), pp.421--422.

  \bibitem{Segal-Tillmann}
    Segal, G. and Tillmann, U.
    \textit{Mapping configuration spaces to moduli spaces}.
    Groups of Diffeomorphisms. Adv. Stud. Pure Math., 52 (2008), pp. 469--477.

  \bibitem{ST07}
    Song, Yongjin and Tillmann, Ulrike.
    \textit{Braids, mapping class groups, and categorical delooping}.
    Mathematische Annalen, 339.2 (2007), pp. 377--393.

  \bibitem{Stasheff63}
    Stasheff, James Dillon.
    \textit{Homotopy Associativity of H-Spaces. I}.
    Transactions of the American Mathematical Society, 108.2 (1963), pp. 275--292.

  \bibitem{SW01}
    {Salvatore}, P. and {Wahl}, N.
    \textit{Framed discs operads and the equivariant recognition principle}.
    ArXiv Mathematics e-prints, http://adsabs.harvard.edu/abs/2001math......6242S (2001).

  \bibitem{Sz10}
    Szepietowski, Bla{\.z}ej
    \textit{EMBEDDING THE BRAID GROUP IN MAPPING CLASS GROUPS}.
    Publicacions Matem{\`a}tiques, 54.2 (2010), pp. 359--368.

  \bibitem{Waj99}
    Wajnryb, Bronislaw.
    \textit{Artin groups and geometric monodromy}.
    Inventiones mathematicae, 138.3 (1999), pp.563--571.

  \bibitem{Waj06}
    Wajnryb, Bronislaw.
    \textit{Relations in the mapping class group}.
    Proceedings of Symposia in Pure Mathematics, (2006), pp. 115--120.

\end{thebibliography}
\end{document}